%% file: main.tex
\documentclass[11pt]{article}

\usepackage[utf8]{inputenc}
\usepackage[margin=1in]{geometry}
\usepackage{authblk}

\usepackage{graphicx}
\usepackage{amsmath,amsthm,amsfonts,latexsym,amssymb,amsbsy}
\usepackage{thmtools}
\usepackage{mathtools}
\usepackage{comment}
\usepackage{multirow}

\usepackage[textwidth=0.8in,color=green]{todonotes}
\usepackage{booktabs}
\usepackage{siunitx}
\newcommand{\sta}{\texttt{STA}}
\newcommand{\rec}{\texttt{REC}}
\newcommand{\msp}{\texttt{MSP}}
\newcommand{\twosp}{\texttt{2SP}}
\newcommand{\bstgmsp}{\texttt{B-Stg-MSP}}
\newcommand{\cstgmsp}{\texttt{C-Stg-MSP}}
\newcommand{\bscenmsp}{\texttt{B-Scen-MSP}}
\newcommand{\cscenmsp}{\texttt{C-Scen-MSP}}
\newcommand{\bstgsta}{\texttt{B-Stg-STA}}
\newcommand{\cstgsta}{\texttt{C-Stg-STA}}
\newcommand{\bscensta}{\texttt{B-Scen-STA}}
\newcommand{\cscensta}{\texttt{C-Scen-STA}}

\newcommand{\cscenrec}{\texttt{C-Scen-REC}}
\newcommand{\Npairs}{\mathcal{N}^{\texttt{pairs}}}

\newcommand{\zsta}[1]{z_{#1}^{\text{\sta{}},*}}
\newcommand{\Nrec}{\mathcal{N}^{\text{\rec{}}}}
\usepackage{enumitem}
\setlist[itemize]{itemsep=5pt}

\usepackage{subcaption}

\usepackage{tikz}
\usetikzlibrary{arrows.meta,positioning,calc}

\tikzset{
  sctree/.style={x=1.5cm, y=0.75cm, line width=0.6pt,
    >={Straight Barb[length=1.3mm, width=2.2mm]}},
  stnode/.style={draw, rounded corners=2.5pt, minimum width=0.7cm,
    minimum height=0.62cm, inner sep=1.5pt, font=\footnotesize, fill=white},
  stmerged/.style={stnode, fill=black!8},
  stedge/.style={->},
}
\newcommand{\sctreebbox}{\path[use as bounding box] (-0.35,-0.45) rectangle (3.35,7.45);}

\usepackage{url}
\usepackage{xcolor,colortbl}
\usepackage{hyperref}
\hypersetup{
    colorlinks=true,
    linkcolor=blue,
    filecolor=magenta,      
    urlcolor=blue,
    citecolor = blue
}

\usepackage[numbers]{natbib}

\newtheoremstyle{break}
  {\topsep}{\topsep}%
  {\upshape}{}%
  {\bfseries}{}%
  {.5em}{}%
\theoremstyle{break}

\newtheorem{property}{Property}
\newtheorem{remark}{Remark}

\newcommand{\BI}{\begin{itemize}}
\newcommand{\EI}{\end{itemize}}

\title{Optimizing Flexibility and Commitments: Partially-Adaptive Multi-Stage Matrix Assembly Layout Design Problem}
\author[1]{Baturhan Bayraktar}
\author[2]{Merve Bodur}
\author[1]{Martin Grunow}
\author[1]{Rainer Kolisch}
\affil[1]{Department of Operations and Technology, TUM School of Management\\
Technical University of Munich, Arcisstr. 21, Munich 80333, Germany\\
\href{mailto:baturhan.bayraktar@tum.de}{\texttt{baturhan.bayraktar@tum.de}},
\href{mailto:martin.grunow@tum.de}{\texttt{martin.grunow@tum.de}},
\href{mailto:rainer.kolisch@tum.de}{\texttt{rainer.kolisch@tum.de}}}
\affil[2]{School of Mathematics and Maxwell Institute for Mathematical Sciences\\
University of Edinburgh, \href{mailto:merve.bodur@ed.ac.uk}{\texttt{merve.bodur@ed.ac.uk}}}
\date{}

\begin{document}

\numberwithin{equation}{section}

\maketitle

\begin{abstract}
Matrix assembly layouts provide manufacturing flexibility under product variety and demand volatility through locally reconfigurable production stations. Periodic redesign exploits this flexibility, yet existing models allow only annual revisions. Finer calendars expose a flexibility-commitment trade-off: advance layout commitments support coordination, supplier agreements, and discounts. We introduce the layout revision timing problem, 
%whose revision policy determines revision points and committed layout sequences under demand uncertainty, while operational production planning is embedded within those layouts. 
jointly choosing revision points and committed layout sequences under demand uncertainty while production adapts within those layouts.
We formulate it as a partially adaptive multistage stochastic program with two complementary extensions. Our cost-based extension places the economic value of advance commitment in the objective rather than prescribing revision budgets; our scenario-adaptive extension adapts revision timing to revealed information. 
For station-use minimization, we establish a cost decomposition across commitment intervals, yielding an exact dynamic program that reuses station subproblem solutions across model variants and parameter settings. As a more general alternative, we develop an integer L-shaped method strengthened by monotonicity-based cuts for budgeted formulations and enhanced cuts exploiting this decomposition.
On automotive instances, including scenario trees with nearly 10,000 nodes, our combined model’s average savings in expected station-use cost at the base-case discount reach 4.99\% versus full adaptivity and 28.65\% versus full commitment. Scenario adaptivity is most valuable when demand paths diverge and discounts justify commitment on some paths but not others, as in electrification-driven scenarios.

\vspace{1ex}
\end{abstract}

%%%%%%%%%%%%%%%%%%%%%%%%%%%%%%%%%%%%%%%%%%%%%%%%
%%%%%%%%%%%%%%%%%%%%%%%%%%%%%%%%%%%%%%%%%%%%%%%%
%% content

\section{Introduction}

\emph{Matrix assembly layouts} have emerged as a prominent response to the growing need for production systems that can accommodate high product variety, volatile demand, evolving sustainability requirements, and rapid changes in products and production technologies.
These pressures are particularly visible in automotive manufacturing, where rising customization and the transition across powertrain technologies have intensified the need for flexible production.
Consistent with this shift, the AMS/ABB Automotive Manufacturing Outlook Survey 2025 identifies the transition toward more flexible and modular manufacturing processes as one of the principal strategies being pursued by automotive manufacturers \citep{ams2025survey}.

\begin{figure}[t]
    \centering
    \includegraphics[width=0.7\linewidth]{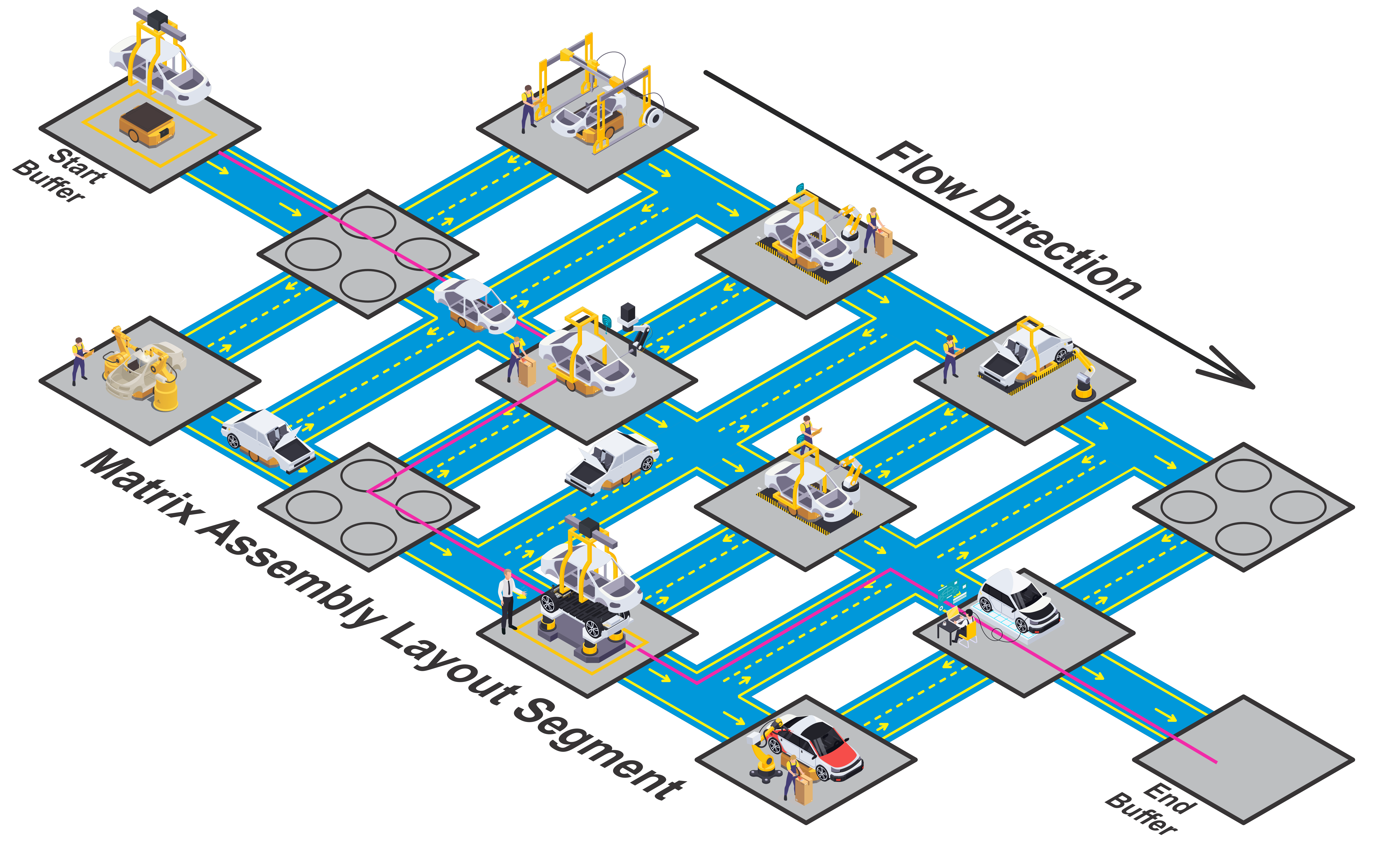}
    \caption{Representative matrix assembly layout schematic. Design icons from~\citep{freepik}.}\label{fig:matrixillustration}
\end{figure}

Matrix assembly layouts are particularly suited to production segments with high product variety and heterogeneous workloads and task requirements.
In these unpaced systems, modular workstations are arranged in a grid, and automated guided vehicles (AGVs) transport products between stations.
Rather than following a common station sequence, each product is transported along a customized route that visits only the stations needed to perform its required tasks.
A matrix segment can be integrated into a broader production system rather than operated as a standalone facility \citep{hottenrott2019flexible,li2026matrix}.
Figure~\ref{fig:matrixillustration} illustrates a representative matrix layout: the pink line traces one possible product-specific AGV route, while the buffer station represents an interface between the matrix segment and the surrounding production system.

The physical structure of a matrix layout provides several complementary forms of \emph{flexibility}.
Routing flexibility arises because a task may be assigned to alternative stations, whereas operation flexibility arises because task sequences may vary subject to precedence constraints rather than being fixed in advance.
There is also volume flexibility because production capacity can be adjusted through local station-level changes as demand changes.
The same modularity allows the layout to respond locally to disruptions without requiring the entire production line to be shut down. % \citep{bayraktar}
Together, these capabilities form the operational basis for dynamic layout planning, in which the matrix segment can be periodically redesigned in response to changes in product mix, demand conditions, or station availability over the product lifecycle.

In current practice, matrix assembly layout design is undertaken at the tactical level of the facility-planning hierarchy, after the main production segments have been identified and before operational production decisions are made. The planning horizon typically covers the product lifecycle and may, for instance, span five years. Major layout reconfigurations are generally coordinated between consecutive production years, for example during planned maintenance periods such as end-of-summer shutdowns \citep{bayraktar}. Within this horizon, the planner selects a \emph{layout configuration} for each production year, comprising the active station set and the task-to-station assignment. The configuration is determined before the year's actual demand is known; once demand is observed, the operational planning step determines production quantities within the selected configuration.

\begin{figure}[t]
    \centering
    \includegraphics[width=0.95\linewidth]{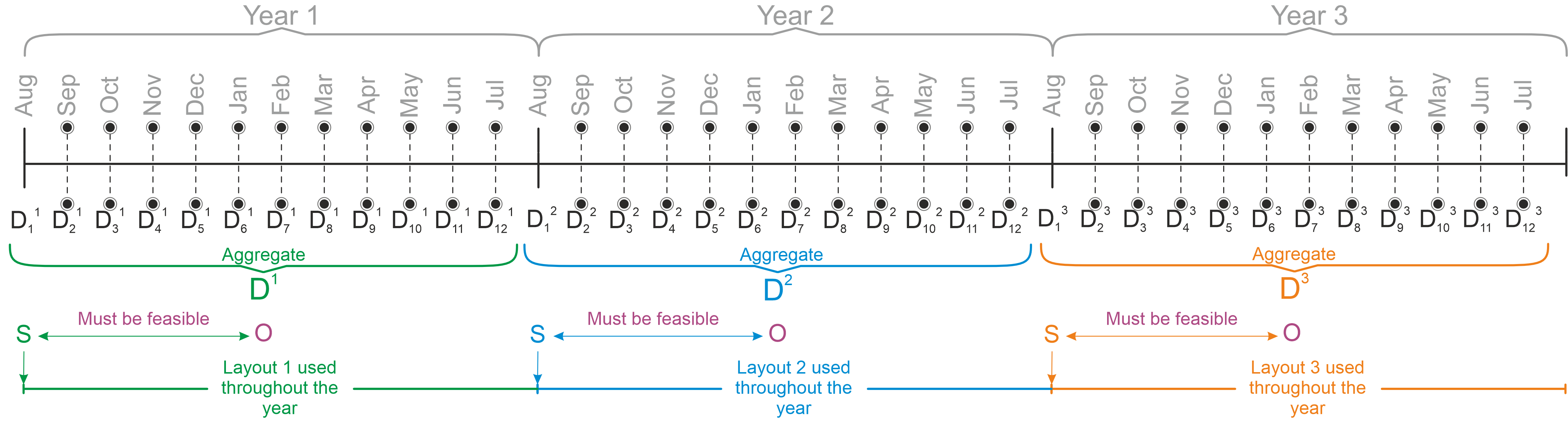}
    \caption{Planning process in existing dynamic matrix assembly layout models.}
    \label{fig:timeline-baseline}
\end{figure}

Existing \emph{dynamic matrix-layout planning} models formalize this planning process on an annual time grid and show that periodic redesign can yield efficiency gains under both deterministic and stochastic demand \citep{bayraktar,bayraktar2}. Figure~\ref{fig:timeline-baseline}, synthesized from these models, depicts the corresponding annual planning template. Each production year contains a pre-demand layout-planning substage (S), in which the layout configuration is selected based on aggregated demand forecasts, followed by a post-demand operational substage (O), in which production quantities are determined after demand is realized. We refer to the pre-demand substage as \emph{strategic} to distinguish its information timing from that of post-demand operational recourse. At this annual resolution, the stochastic formulation requires the selected configuration to permit feasible production under every demand realization represented for the year. This structure captures an important source of value: matrix layouts can be redesigned over time as demand changes. However, because admissible revision times are \emph{fixed} at annual boundaries rather than optimized, these models cannot fully exploit the local reconfigurability that distinguishes matrix layouts from traditional serial assembly lines. In particular, they cannot evaluate whether an intermediate layout revision would better adapt the layout to demand patterns that emerge within the year. Beyond these missed revision opportunities, the annual representation also limits the feasibility statement delivered by the model: it verifies feasibility for demand and capacity measured at the annual planning unit, but it does not certify feasibility for a finer operational calendar in which the same annual demand may be distributed unevenly across shorter intervals.

\begin{figure}[b]
    \centering
    \includegraphics[width=0.95\linewidth]{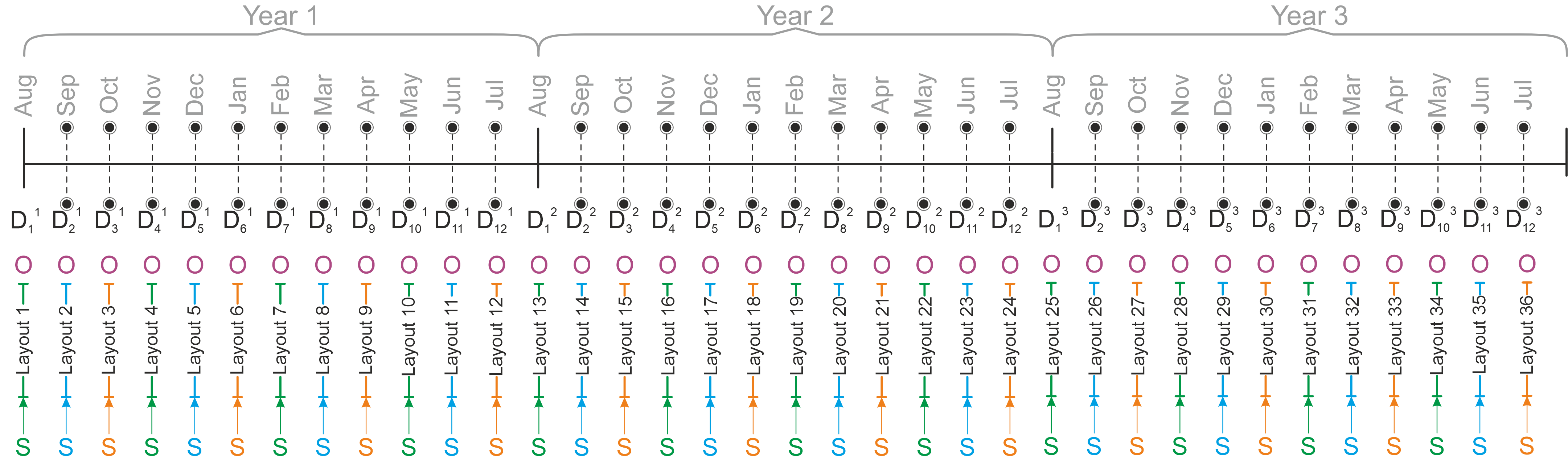}
    \caption{Planning process with increased temporal granularity (monthly).}
    \label{fig:timeline-fine}
\end{figure}

To expose such demand-responsive flexibility, the planner can refine the time grid. Dividing the horizon into monthly or quarterly intervals creates additional candidate points at which the layout plan may be revised and permits production feasibility to be assessed within the year; Figure~\ref{fig:timeline-fine} illustrates the monthly case. 
Planning at this finer resolution makes it possible to respond to seasonal or regime-driven demand patterns that are obscured by annual aggregation. Yet revising the plan at every opportunity may leave insufficient time to prepare the resulting configuration changes. Station activation may require hiring personnel, installing equipment, and arranging infrastructure, whereas task reassignment may require worker training, equipment movement, and validation. We refer to the interval between consecutive selected revision points as a \emph{commitment}: at the earlier revision point, the planner fixes the sequence of layout configurations for the intervening periods until the next revision, providing time to coordinate their implementation. Beyond this coordination benefit, longer commitments can create economic value through longer-term arrangements, such as supplier agreements, bulk purchasing terms, or priority delivery.

\begin{figure}[!htb]
    \centering
    \includegraphics[width=0.95\linewidth]{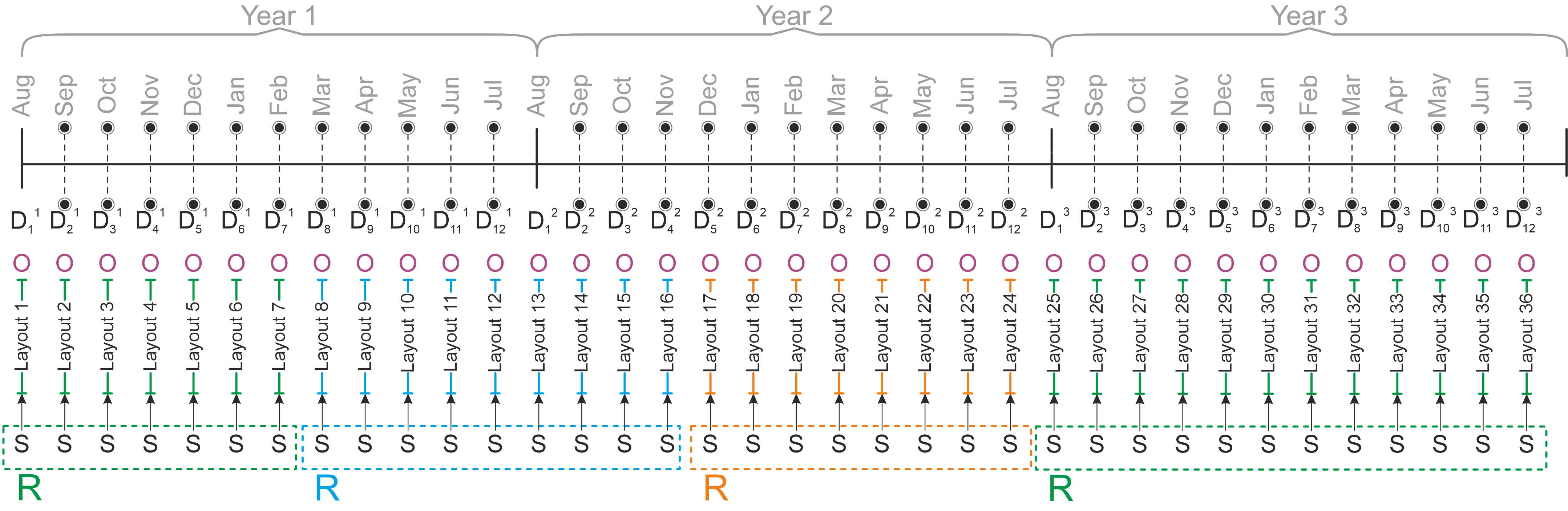}
    \caption{Planning process under an optimized revision policy with layout commitments.}
    \label{fig:timeline-finer}
\end{figure}

To optimize this flexibility-commitment trade-off under demand uncertainty, we formalize the \emph{layout revision timing problem}. Given a fine-grained planning horizon with candidate revision opportunities, the problem aims to find a \emph{revision policy}: consistently with the information available at each decision point, the policy determines which opportunities are used as revision points and which sequence of layout configurations is committed between successive revisions. Along each realization of uncertainty, the selected revision points partition the horizon into commitment intervals. Within each interval, configurations may differ across periods, but the committed sequence remains fixed until another revision is selected, while production quantities are determined operationally within those configurations. Figure~\ref{fig:timeline-finer} illustrates this structure where (R) marks the revision points selected by the policy and the dashed boxes indicate the resulting commitment intervals.

We address this problem with a partially adaptive multistage stochastic programming (MSP) formulation reflecting the strategic-operational timing of matrix-layout planning. Revision timing and layout configurations are optimized jointly, while operational production decisions are embedded to certify production feasibility under the committed layouts. Because station activation dominates layout cost, we solve the strategic design problem hierarchically: the station problem first determines an optimal revision policy and station-use profile, after which the reconfiguration problem minimizes task-to-station assignment changes while preserving that profile. 
To model the flexibility-commitment trade-off, we build on the budgeted, stage-wise partially adaptive MSP framework \citep{sezenece} and extend it in two complementary directions. The cost-based extension places the economic value of advance commitment in the objective, while the scenario-adaptive extension allows revision timing to vary across scenario paths. Together, they yield a coherent family of four models.
To solve this family, we develop a dynamic program for the station problem, present an integer L-shaped method as a more general alternative, and use nested decomposition for the reconfiguration problem. We strengthen both decomposition methods with structure-based cuts. 
    
Our contributions to the dynamic matrix-layout and partially adaptive MSP literatures can be summarized as follows.
\begin{itemize}
    \item \textbf{Problem definition.}
    We introduce the novel \emph{layout revision timing problem}, extending dynamic matrix assembly layout design to jointly optimize layout plans and revision timing under demand uncertainty. This moves beyond the fixed revision calendars of existing dynamic models and exposes the resulting flexibility-commitment trade-off.

    \item \textbf{Modeling framework.}
    We broaden partially adaptive MSP along two complementary dimensions. To our knowledge, our cost-based variant is the first to capture explicitly the economic incentive structure of advance commitment, whereas our scenario-adaptive variant generalizes stage-wise revision timing to policies conditioned on observed scenario histories. The resulting framework provides an economically grounded, information-responsive alternative to prescribing the degree of adaptivity.

    \item \textbf{Solution methodology.}
    We characterize the commitment structure induced by revision policies and exploit station-problem separability to derive an exact dynamic program for all four partially adaptive model variants. The dynamic program uses polynomially many station subproblems whose optimal solutions can be reused across variants, budgets, and discount functions. For the more general integer L-shaped method, we establish monotonicity in revision authorizations for budgeted partially adaptive MSPs, yielding, to our knowledge, the first monotonicity-based optimality cuts for this class, a generic enhancement beyond the matrix-layout application. The same separability yields enhanced cuts for all four station problem variants. For the reconfiguration problem, we extend asymmetric Lagrangian switch cuts of~\citep{bayraktar2} to commitment-dependent costs within nested decomposition.
    
    \item \textbf{Computational and managerial insights.}
    Our computational study demonstrates both algorithmic impact and managerial value. In the cut-comparison experiments, the enhanced cuts reduce average integer L-shaped iteration counts by roughly two to three orders of magnitude.
    For each scenario tree, the dynamic program solves the shared station subproblems and all four station-problem variants at their tested parameter settings in under 41 minutes, even for trees with nearly 10,000 nodes.
    The model comparisons quantify the value of economically pricing commitments: at the base-case discount, cost-based scenario-adaptive planning lowers expected station-use cost by up to 4.99\% relative to full adaptivity and 28.65\% relative to full commitment across the tested datasets.    
    More broadly, the policy experiments reveal how demand-path heterogeneity and commitment incentives jointly govern the value of adaptivity: scenario adaptivity is often most valuable when demand paths diverge and commitment incentives are intermediate.
\end{itemize}

The remainder of the paper is organized as follows. Section~\ref{sec:literature} positions our work within the matrix assembly layout and partially adaptive MSP literatures. Section~\ref{sec:model} develops the proposed modeling framework and its matrix-layout specialization. Section~\ref{sec:solution} presents the solution methodology. Section~\ref{sec:experiments} reports the computational study and resulting managerial insights. Section~\ref{sec:conclusions} provides concluding remarks.

% ---------------------------------- %
\section{Literature Review}
\label{sec:literature}

Our contributions draw on and advance two literature streams, which we review in turn: matrix assembly layout design and partially adaptive multistage stochastic programming.
  
\subsection{Matrix Assembly Layout Design}
The optimization literature on matrix assembly layouts originated in \emph{static, deterministic} design. \citet{hottenrott2019flexible} formulated the first model, jointly determining station locations, task allocations, and anticipated product routes while lexicographically minimizing the number of opened stations and transportation effort generated by those routes. Within this static setting, \citet{schumacher2021economic} incorporated heterogeneous production resources and replaced the lexicographic criteria with a monetary objective covering opened stations, deployed resources, and product transportation.

The \emph{dynamic} stream began with \citet{bayraktar}, who considered known, time-varying demand over a sequence of one-year planning periods. Their model selects a layout for each year and lexicographically minimizes the total number of active stations over the horizon, the number of initial task assignments and subsequent layout reconfigurations, and the total distance traveled by products along their AGV routes. \citet{bayraktar2} extended this setting to uncertain demand through a multistage stochastic program: each year's layout is selected before demand for that year is revealed, whereas operational production quantities are determined after demand is observed.
  
These dynamic models share two \emph{restrictions} central to our study. Both operate on an annual planning grid and represent demand through yearly aggregate product volumes. At that resolution, they neither represent within-year demand patterns nor certify production feasibility at a finer calendar resolution. Moreover, although they optimize reconfigurations between consecutive layouts, they do not optimize when the future layout plan should be reconsidered. The deterministic model fixes the multi-year layout plan at the outset, whereas the stochastic model permits reconsideration at every annual decision point; neither selects which revision opportunities should be used or how that selection should adapt to observed demand. 
At the broader field level, the review of 125 publications by \citet{li2026matrix} identifies limited attention to dynamic reconfiguration during operation and to the timing of reconfiguration triggers. Although that observation concerns a broader operational setting, it reinforces the unresolved tactical planning question studied in this paper: how should a manufacturer time revisions to its layout plan while balancing demand responsiveness against the coordination and economic benefits of committing to layouts in advance?

\subsection{Partially Adaptive Multistage Stochastic Programming}

\emph{Two-stage and multistage stochastic programming} provide benchmark information structures for multi-period planning decisions: a two-stage model fixes their full-horizon sequence before uncertainty is observed, whereas an MSP model allows their successive values to depend on the uncertainty history available at each stage~\citep{birge2011introduction}. 
Motivated by the computational difficulty of multistage formulations \citep{SHAPIRO20061}, \citet{zou2018partially} introduced a hybrid approximation for generation expansion planning. Capacity-expansion decisions remain adaptive through a pre-specified switching period; at the switching point, the expansion plan for the remainder of the horizon is determined based on the uncertainty history observed up to that point and kept fixed thereafter.
  
A different motivation underlies the subsequent \emph{revision-timing} literature: contractual, lead-time, and coordination requirements can make stage-by-stage revision infeasible or undesirable, creating a practical need to balance flexibility with commitment. To strike this balance, \citet{basciftci2024adaptive} introduced adaptive two-stage stochastic programming, which optimizes one revision point for each decision component requiring limited flexibility before uncertainty unfolds. Their capacity-expansion experiments showed that a policy even with one optimally timed revision per strategic decision component can yield objective values close to the MSP optimum. \citet{sezenece} extended this framework to partially adaptive multistage stochastic programming by allowing multiple revision points for each state-variable type, with a given revision budget capping the number of such points for each type.  They treat the entire revision schedule as a here-and-now decision, optimizing it before uncertainty unfolds and applying the same revision points across all scenarios. This budgeted, stage-wise framework provides the principal benchmark for our work. We extend it in two complementary directions: (i) our cost-based variant replaces the pre-specified revision budget with commitment-dependent costs that explicitly represent the economic incentives for advance commitment; and (ii) our scenario-adaptive variant allows revision timing to adapt to the observed uncertainty history.
  
In work developed in parallel with ours, \citet{wang2026balancing} introduce $K$-revision MSP, a related framework closely connected to our scenario-adaptive extension (ii). Their framework begins with a plan for \emph{all} future strategic decisions and limits its replacement along each scenario path to at most $K$ times. $K$-revision uses same-stage (\emph{lag-0}) timing: after a stage's uncertainty is revealed, replacing the active plan may alter the current strategic decision and provisionally replan the remaining horizon. Our extension (ii) instead uses one-period-ahead (\emph{lag-1}) timing: the revision policy selects each revision opportunity one stage before the layout plan may be updated, providing time to coordinate implementation; the resulting commitment then covers only the layout sequence until the next selected revision point. The central additional distinction is economic. $K$-revision, like budgeted stage-wise framework of \citet{sezenece}, limits revisions through a pre-specified cap; our cost-based extension (i) instead lets commitment-dependent costs determine when and how often to revise. To our knowledge, it is the first partially adaptive MSP model to represent explicitly the economic value of advance commitment.

Methodologically, \citet{sezenece} develop a decomposition method for their budgeted, stage-wise partially adaptive MSP formulation. They place the here-and-now revision schedule in the master problem and solve each fixed-schedule MSP in extensive form, generating Benders cuts from its LP relaxation and selectively adding integer L-shaped cuts from its mixed-integer formulation; exact reductions of nonanticipativity constraints and heuristic preprocessing further accelerate the method. For our four partially adaptive variants, we develop a decomposition approach that uses integer optimality cuts as a common baseline and exploits two complementary value-function properties. Monotonicity yields, to our knowledge, the first monotonicity-based optimality cuts for the budgeted partially adaptive MSP class, while application-specific separability yields an interval decomposition that localizes each value-function component to its relevant commitment interval. These results underpin enhanced optimality cuts for all four variants. Moreover, the interval decomposition yields an exact dynamic program for all four station-problem variants, with station subproblem solutions reusable across variants and parameter settings. For the reconfiguration problem, we extend the asymmetric Lagrangian switch cuts of~\citet{bayraktar2} to commitment-dependent costs within nested decomposition.
% ---------------------------------- %
% ---------------------------------- %
\section{Modeling Framework and Matrix-Layout Specialization}
\label{sec:model}

This section develops the proposed revision-policy framework and its matrix-layout specialization. We first describe the underlying sequential decision process (Section~\ref{sec:sequential-process}), then formulate the budgeted, stage-wise partially adaptive baseline (Section~\ref{sec:budgeted-stagewise-pa}), its cost-based and scenario-adaptive extensions (Sections~\ref{sec:cost-based-pa} and~\ref{sec:scenario-adaptive}, respectively), and the resulting four-model family (Section~\ref{sec:combined-framework}). Throughout these subsections, we present the framework in general terms, with brief matrix-layout interpretations of its decisions and uncertainties. We finally provide the full matrix-layout specialization by defining its decisions, hierarchical objectives, and feasibility structure (Section~\ref{sec:objfunc-mald}).

\subsection{Sequential Decision-Making Process}\label{sec:sequential-process}

We consider a finite planning horizon consisting of $T$ decision-making stages, indexed by $t\in\mathcal{T}:=[1,T]$, where $[a,b]:=\{a,a+1,\ldots,b\}$ for integers $a\leq b$. Each stage corresponds to a production period in our application (e.g., a month, quarter, or year) and comprises two sequential substages. In the \emph{strategic substage}, a strategic decision establishes the system configuration that links successive stages; in the subsequent \emph{operational substage}, an operational decision governs within-stage execution under that configuration. Uncertainty is revealed across these substages: strategic uncertainty describes how the planning environment evolves across stages and is observed before the strategic decision, whereas operational uncertainty is revealed after that decision and before the operational response.

We represent strategic uncertainty by the stochastic process $\{\boldsymbol{\xi}_t\}_{t\in\mathcal{T}}$ (market conditions in our application). At the strategic substage of stage $t$, the history $\boldsymbol{\xi}^t := (\boldsymbol{\xi}_1,\ldots,\boldsymbol{\xi}_t)$ is observed, after which the strategic decision vector $\boldsymbol{s}_t$ (the layout configuration for the corresponding production period) is selected. We represent operational uncertainty by the uncertain vector $\boldsymbol{D}_t$ (demand volumes by product), whose admissible realizations depend on this history. At the operational substage, a realization $\boldsymbol{d}_t$ of $\boldsymbol{D}_t$ is observed, after which the operational decision vector $\boldsymbol{o}_t$ (task-level production allocations across stations) is selected subject to $\boldsymbol{s}_t$ and $\boldsymbol{d}_t$. Under full adaptivity, $\boldsymbol{s}_t$ may depend on all strategic information available at stage $t$; the revision policies introduced later restrict this dependence. The resulting integrated planning problem under full adaptivity has the following nested representation:
\begin{align*}
    &\min_{\boldsymbol{s}_1\in \mathcal{S}_1(\boldsymbol{s}_0)} \Bigg\{ f_1^{\text{str}}(\boldsymbol{s}_0, \boldsymbol{s}_1) + \sup_{\boldsymbol{d}_1\in\operatorname{supp}(\boldsymbol{D}_1\mid\boldsymbol{\xi}^1)} \big\{\min_{\boldsymbol{o}_1\in\mathcal{O}_1(\boldsymbol{s}_1,\boldsymbol{d}_1)}f^{\text{opr}}_1(\boldsymbol{o}_1)\big\}  \\ 
    &\quad + \mathbb{E}_{\boldsymbol{\xi}_2\mid\boldsymbol{\xi}^1}\bigg[
    \min_{\boldsymbol{s}_2 \in \mathcal{S}_2(\boldsymbol{s}_1)} \Big\{  f_2^{\text{str}}(\boldsymbol{s}_1, \boldsymbol{s}_2) + \sup_{\boldsymbol{d}_2\in\operatorname{supp}(\boldsymbol{D}_2\mid\boldsymbol{\xi}^2)}\big\{\min_{\boldsymbol{o}_2\in \mathcal{O}_2(\boldsymbol{s}_2,\boldsymbol{d}_2)}f^{\text{opr}}_2(\boldsymbol{o}_2)\big\} + \dots \\
    &\quad \ \  + \mathbb{E}_{\boldsymbol{\xi}_T\mid\boldsymbol{\xi}^{T-1}}\bigg[
    \min_{\boldsymbol{s}_T \in \mathcal{S}_T(\boldsymbol{s}_{T-1})} \{  f_T^{\text{str}}(\boldsymbol{s}_{T-1}, \boldsymbol{s}_T) + \sup_{\boldsymbol{d}_T\in\operatorname{supp}(\boldsymbol{D}_T\mid\boldsymbol{\xi}^T)}\{\min_{\boldsymbol{o}_T\in \mathcal{O}_T(\boldsymbol{s}_T,\boldsymbol{d}_T)}f^{\text{opr}}_T(\boldsymbol{o}_T)\} \} \bigg] \dots \Big\} \bigg] \Bigg\},
\end{align*}
where $\boldsymbol{s}_0$ denotes the initial system state; $\operatorname{supp}(\boldsymbol{D}_t\mid\boldsymbol{\xi}^t)$ denotes the conditional support of $\boldsymbol{D}_t$ given $\boldsymbol{\xi}^t$; $\mathcal{S}_t$ and $\mathcal{O}_t$ denote the feasible sets, including variable domains, for strategic and operational decisions, respectively; and $f_t^{\mathrm{str}}$ and $f_t^{\mathrm{opr}}$ are the corresponding cost functions.

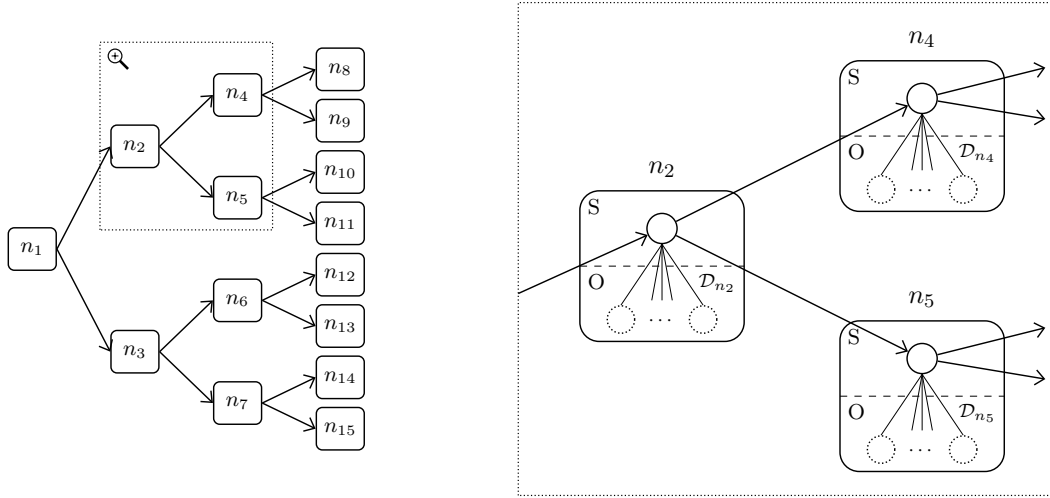
\begin{figure}[htbp!]
    \centering
    \resizebox{0.85\linewidth}{!}{\input{./Figures/tree_overview}}
    \caption{Composite scenario tree with strategic (S) and operational (O) substages.}
    \label{fig:tree-overview}
\end{figure}

To obtain a computationally tractable representation of this two-substage structure, we employ a scenario tree with composite nodes. Let $\mathcal{N}_t$ denote the set of nodes at stage $t$. Each node $n\in\mathcal{N}_t$ combines one possible history $\boldsymbol{\xi}^t$ of strategic uncertainty with a finite set $\mathcal{D}_n$ that approximates the conditional support of $\boldsymbol{D}_t$ given that history. We define $\mathcal{N}:=\bigcup_{t\in\mathcal{T}}\mathcal{N}_t$. For each $t\in\mathcal{T}$, let $\Npairs_t:=\{\{n_1,n_2\}\subseteq\mathcal{N}_t:n_1\neq n_2\}$ denote the set of unordered pairs of distinct stage-$t$ nodes. The links between nodes describe the evolution of strategic uncertainty. Figure~\ref{fig:tree-overview} illustrates this structure, enlarging nodes $n_2$, $n_4$, and $n_5$ to distinguish their strategic information states (solid circles) from the associated conditional demand realizations (dotted circles).

For each $n\in\mathcal{N}$, let $t(n)$ denote its stage, $a(n)$ its immediate predecessor, and $\mathcal{C}(n)$ its set of children. For two nodes $n_1,n_2\in\mathcal{N}$, let $a(n_1,n_2)$ denote their most recent common ancestor, i.e., their common ancestor with the largest stage index. Moreover, $q_{nn'}$ denotes the conditional transition probability from $n$ to $n'\in\mathcal{C}(n)$, while $p_n>0$ denotes the probability of reaching the strategic history associated with $n$. Accordingly, $p_{n'}=p_nq_{nn'}$ for every $n'\in\mathcal{C}(n)$. The root node is unique, with $\mathcal{N}_1=\{1\}$ and $p_1=1$. We introduce a virtual predecessor $a(1)$ and set $\boldsymbol{s}_{a(1)}:=\boldsymbol{s}_0$.

Under this finite representation, each conditional supremum in the nested model becomes a maximum over $\mathcal{D}_n$. The fully adaptive node-based formulation, which we use as the MSP benchmark, is
\begin{equation*}
\begin{aligned}
(\msp{}):
\min_{\{\boldsymbol{s}_n\}_{n\in\mathcal{N}}} \ \ 
& \sum_{n\in\mathcal{N}}p_n
\left[
f_n^{\mathrm{str}}(\boldsymbol{s}_{a(n)},\boldsymbol{s}_n)
+
\max_{\boldsymbol{d}_n\in\mathcal{D}_n}
\left\{
\min_{\boldsymbol{o}_n\in\mathcal{O}_n(\boldsymbol{s}_n,\boldsymbol{d}_n)}
f_n^{\mathrm{opr}}(\boldsymbol{o}_n)
\right\}
\right] \\
\text{s.t.} \ \ 
&(\boldsymbol{s}_{a(n)},\boldsymbol{s}_n)\in\mathcal{S}_n
\quad \forall n\in\mathcal{N},
\end{aligned}
\end{equation*}
where $\mathcal{S}_n$, $\mathcal{O}_n$, $f_n^{\mathrm{str}}$, and $f_n^{\mathrm{opr}}$ are the node-based counterparts of the stage-indexed feasible sets and cost functions defined above. More broadly, \msp{} shares the strategic-operational nesting of multi-horizon stochastic programming \citep{kaut2014multi}. In our formulation, the operational block is evaluated robustly, through the maximum over $\mathcal{D}_n$, rather than through the expected-recourse evaluation commonly used in multi-horizon applications. More fundamentally, our uncertainty structure departs from the canonical multi-horizon construction: the current operational realization may influence the strategic uncertainty revealed at the next stage, as clarified in the following remark.

\begin{remark}[Operational-to-strategic uncertainty dependence]
For presentational simplicity, both the nested and node-based representations above leave implicit the possible dependence of the next strategic state on the current operational realization. Our framework accommodates this dependence by allowing $\boldsymbol{\xi}_{t+1}$ to follow a conditional distribution given $(\boldsymbol{\xi}^t,\boldsymbol{d}_t)$. In the finite scenario tree, child states and transition probabilities can depend on the current demand realization; this dependence is encoded in the tree construction, while \msp{} retains its stated form. In the standard datasets of Section~\ref{sec:exp-design}, each realized demand vector at a non-leaf node determines the baseline for its corresponding child. The regime-driven datasets instead use the latent-state construction described in Section~\ref{sec:value-pa}. Neither construction alters the requirement that $\boldsymbol{s}_n$ admit a feasible operational response for every $\boldsymbol{d}_n\in\mathcal{D}_n$.
\end{remark}

The revision policies developed next restrict the strategic adaptivity of \msp{} while preserving demand-responsive operational decisions within each committed layout configuration.

\subsection{Budgeted, Stage-Wise Partial Adaptivity}\label{sec:budgeted-stagewise-pa}

Following \citet{sezenece}, we first impose a budgeted, stage-wise revision schedule on the strategic decisions in \msp{}. Under this framework, the complete schedule of revision stages is selected before uncertainty unfolds and applies uniformly across scenario histories. Stage~1 begins the initial commitment interval, and each designated revision stage begins a new one. The strategic decisions within an interval may depend on the information available when that interval begins but cannot respond to information revealed later within the same interval.

Let $\mathcal{I}$ index the components of the strategic decision vector $\boldsymbol{s}_n=(s_{ni})_{i\in\mathcal{I}}$. To represent the coordinated revision of a layout plan, we assign a single revision schedule to the entire vector, rather than a separate schedule to each state-variable type, as permitted by the underlying partially adaptive framework. In our application, $\boldsymbol{s}_n$ comprises two binary decision families, namely, station activation and task-to-station assignment, both of which follow this common schedule. For each $i\in\mathcal{I}$, let $M_i\geq 0$ be a finite upper bound on $|s_{n_1,i}-s_{n_2,i}|$ for any two nodes at the same stage. If $L_i\leq s_{ni}\leq U_i$ at every node, we may take $M_i=U_i-L_i$; in particular, $M_i=1$ suffices for binary components.

Let $\boldsymbol{r}:=(r_t)_{t\in[2,T]}$ denote the here-and-now revision schedule, where $r_t=1$ if stage $t$ is a revision stage and $r_t=0$ otherwise. The parameter $\mu\in[0,T-1]$ bounds the number of revision stages in $[2,T]$; the initial planning point at stage~1 does not count against this budget. The resulting budgeted, stage-wise partially adaptive MSP is
\begin{subequations}
\begin{align}
\hspace*{-0.23cm} (\bstgmsp{}): 
\min_{\substack{
{\{\boldsymbol{s}_n\}}_{n\in\mathcal{N}},\\
{\{r_t\}}_{t\in[2,T]}
}} \ \  
& \mathrlap{\sum_{n\in\mathcal{N}}p_n
\Big[
f_n^{\mathrm{str}}(\boldsymbol{s}_{a(n)},\boldsymbol{s}_n)
+
\max_{\boldsymbol{d}_n\in\mathcal{D}_n}
\Big\{
\min_{\boldsymbol{o}_n\in\mathcal{O}_n(\boldsymbol{s}_n,\boldsymbol{d}_n)}
f_n^{\mathrm{opr}}(\boldsymbol{o}_n)
\Big\}
\Big]} \label{form:pamsp-objective}\\
\text{s.t.} \ \ 
&(\boldsymbol{s}_{a(n)},\boldsymbol{s}_n)\in\mathcal{S}_n
&& \forall n\in\mathcal{N}, \label{form:pamsp-drmal}\\*[0.1cm]
&\sum_{t \in [2,T]} r_t\leq\mu \label{form:pamsp-cons3} \\*[-0.3cm]
&s_{n_1,i}-s_{n_2,i}
\leq M_i \cdot \hspace*{-0.8cm} \sum_{t'=t(a(n_1,n_2))+1}^{t} \hspace*{-0.7cm} r_{t'}
&& \forall t\in[2,T],\ \{n_1,n_2\}\in\Npairs_t,\ i\in\mathcal{I},
\label{form:pamsp-cons1}\\
&s_{n_2,i}-s_{n_1,i}
\leq M_i \cdot \hspace*{-0.8cm} \sum_{t'=t(a(n_1,n_2))+1}^{t} \hspace*{-0.7cm} r_{t'} \ \ 
&& \forall t\in[2,T],\ \{n_1,n_2\}\in\Npairs_t,\ i\in\mathcal{I},
\label{form:pamsp-cons2}\\
&r_t\in\{0,1\}
&& \forall t\in[2,T]. \label{form:pamsp-cons4}
\end{align}
\end{subequations}

Constraint~\eqref{form:pamsp-drmal} imposes strategic feasibility; constraint~\eqref{form:pamsp-cons3} limits the number of revision stages. Given two distinct nodes $n_1$ and $n_2$ at stage $t$, the sum in constraints~\eqref{form:pamsp-cons1}-\eqref{form:pamsp-cons2} counts the revision stages after their most recent common ancestor and through stage $t$. A zero count makes the paired inequalities enforce $s_{n_1,i}=s_{n_2,i}$ for every $i\in\mathcal{I}$, whereas a positive count makes them nonbinding under the bounds $M_i$. Hence, these constraints allow stage-$t$ decisions to distinguish the two histories only if a revision is scheduled after their most recent common ancestor and no later than stage $t$.

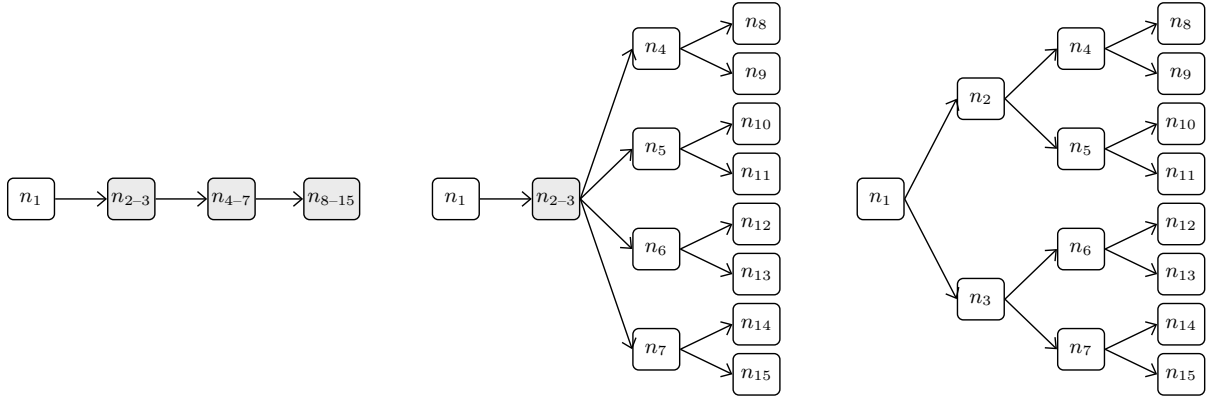
\begin{figure*}[htbp!]
    \centering
    \begin{subfigure}[c]{0.32\textwidth}
        \resizebox{0.95\textwidth}{!}{\input{Figures/tree_mu0}}
        \caption{$\mu=0$: Full commitment (\twosp{})}
    \end{subfigure}~\hspace{0.01\textwidth}
    \begin{subfigure}[c]{0.32\textwidth} 
        \resizebox{0.95\textwidth}{!}{\input{Figures/tree_stage_mu2}}
        \caption{$\mu=2$: $r_2=0$, $r_3=r_4=1$}
        \label{fig:tree-limitedcommit}
    \end{subfigure}~\hspace{0.01\textwidth}
    \begin{subfigure}[c]{0.32\textwidth} 
        \resizebox{0.95\textwidth}{!}{\input{Figures/tree_msp}}
        \caption{$\mu=T-1$: Full adaptivity (\msp{})}
        \label{fig:tree-nocommit}
    \end{subfigure}
    \caption{Strategic-decision structures induced by different stage-wise revision schedules. Merged nodes represent strategic decisions constrained to be identical; the underlying uncertainty tree remains unchanged.}
    \label{fig:partialtree}
\end{figure*}

The extreme values of the revision budget $\mu$ show how \bstgmsp{} interpolates between the full-commitment and fully multistage benchmarks. Figure~\ref{fig:partialtree} illustrates this spectrum through both endpoints and an intermediate revision schedule. At one endpoint, with $\mu=0$, no stage after stage~1 begins a new commitment interval; the full stage-indexed strategic plan is therefore fixed before uncertainty unfolds and cannot vary across scenario histories. We denote this full-commitment benchmark by \twosp{}: the entire strategic plan is selected here and now, while demand-responsive operational decisions constitute the recourse block. At the other endpoint, with $\mu=T-1$, setting $r_t=1$ for all $t\in[2,T]$ makes every later stage a revision stage, removes the additional cross-history linking restrictions, and recovers \msp{}.

\subsection{Cost-Based Partial Adaptivity} \label{sec:cost-based-pa}

In practice, committing to a configuration is not merely an organizational constraint but also an economic decision. By locking in equipment orders, workforce contracts, or supplier agreements in advance, manufacturers can secure preferential pricing, bulk discounts, or priority delivery; longer commitments can yield greater benefits. This creates a trade-off: frequent revisions preserve adaptivity but relinquish these benefits, whereas longer commitments secure them at the cost of flexibility. Existing budget-based formulations of partial adaptivity do not capture this trade-off: they limit revisions without representing the value created by advance commitment. We therefore develop a cost-based model that embeds this value directly in the objective.

To formalize this economic mechanism, we first develop it for the stage-wise revision structure introduced above; Section~\ref{sec:combined-framework} later combines it with scenario-adaptive revision timing. We remove the exogenous revision budget and instead link each commitment-eligible strategic cost to how far in advance the corresponding decision is fixed. Under a stage-wise schedule, the most recent revision stage determines this lead time. Treating stage~1 as the initial planning point, we set the \emph{commitment lag} $\delta_1(\boldsymbol{r}):=0$ and, for each $t\in[2,T]$, define
\begin{equation}\label{eq:stg-lag}
\delta_t(\boldsymbol{r})
:=t-1-\max_{t'\in[2,t]} \big\{(t'-1)r_{t'} \big\}.
\end{equation}
The lag is zero at the initial planning point and at every revision stage, and it increases by one at each subsequent stage of a commitment interval. We translate this lead time into a \emph{cost multiplier} by setting $w_n(\boldsymbol{r}):=h(\delta_{t(n)}(\boldsymbol{r}))$ for every node $n$, where $h:[0,T-1]\rightarrow(0,1]$ is nonincreasing and normalized by $h(0)=1$. Thus, all nodes at the same stage receive the same weight, interpreted as the fraction of the undiscounted cost that remains. A smaller weight therefore represents a greater benefit from advance commitment.

Writing $f_n^{\mathrm{str}}$ for the commitment-eligible strategic costs and leaving $f_n^{\mathrm{opr}}$ unweighted yields our cost-based, stage-wise partially adaptive MSP
\begin{equation*}
\begin{aligned}
(\cstgmsp{}): \ \ 
\min_{\substack{
{\{\boldsymbol{s}_n\}}_{n\in\mathcal{N}},\\
{\{r_t\}}_{t\in[2,T]}
}} \ \ 
&\sum_{n\in\mathcal{N}}p_n
\Big[
w_n(\boldsymbol{r})f_n^{\mathrm{str}}(\boldsymbol{s}_{a(n)},\boldsymbol{s}_n)
+
\max_{\boldsymbol{d}_n\in\mathcal{D}_n}
\Big\{
\min_{\boldsymbol{o}_n\in\mathcal{O}_n(\boldsymbol{s}_n,\boldsymbol{d}_n)}
f_n^{\mathrm{opr}}(\boldsymbol{o}_n)
\Big\}
\Big] \\
\text{s.t.} \ \ 
&\eqref{form:pamsp-drmal},\quad
\eqref{form:pamsp-cons1}-\eqref{form:pamsp-cons4}.
\end{aligned}
\end{equation*}

In \cstgmsp{}, the complete revision schedule remains a here-and-now decision. However, unlike \bstgmsp{}, the model determines when and how often to revise by balancing flexibility against commitment benefits, without imposing an exogenous cap on the number of revisions. Each revision relaxes the cross-history linking restrictions but resets the commitment lag to zero, sacrificing the advance-commitment discount at that stage. Extending a commitment interval can yield greater discounts but prevents later strategic decisions within it from responding to information revealed since the interval began. As detailed in Section \ref{sec:objfunc-mald}, in our matrix-layout specialization, the multiplier applies to the station-use and reconfiguration terms in their respective hierarchical problems; other applications may apply it selectively to any cost components that benefit from advance commitment. Because $w_n(\boldsymbol{r})$ is a decision-dependent multiplier of the strategic cost, the objective is nonlinear; Section~\ref{sec:solution} develops our solution approach.

The relationship to the benchmark models clarifies the role of $h$. If $h\equiv1$, commitment has no economic benefit; setting $r_t=1$ for all $t\in[2,T]$ then recovers \msp{}. Conversely, sufficiently strong commitment benefits may make a schedule with no post-initial revisions optimal, recovering the strategic-decision structure of \twosp{}. In our computational experiments (Section~\ref{sec:experiments}), we evaluate linear, convex, and concave specifications of $h$, representing constant, diminishing, and increasing marginal benefits from longer commitments, respectively.

\subsection{Scenario-Adaptive Revision Timing}\label{sec:scenario-adaptive}

Budgeted, stage-wise partial adaptivity fixes a common revision schedule before
uncertainty unfolds and applies it to every scenario history. Once histories
diverge, however, the value of a future revision may differ across them. We
therefore extend \bstgmsp{} to the budgeted \emph{scenario-adaptive} formulation
\bscenmsp{}, in which revision timing may depend on the observed strategic history.

For each non-leaf node $n\in\mathcal{N}\setminus\mathcal{N}_T$, we introduce a
binary indicator $\rho_n$, selected using the strategic information available at
$n$, before operational uncertainty is realized. These indicators form the
revision policy $\boldsymbol{\rho}$. If $\rho_n=1$, a revision window opens for the
children of $n$: each child $n'\in\mathcal{C}(n)$ begins a new commitment interval
and may select its strategic decision $\boldsymbol{s}_{n'}$ using the information
available at that child. If $\rho_n=0$, the previously committed strategic sequence
remains in force across the children. Thus, $\rho_n$ first affects strategic
decisions one stage ahead and leaves the strategic decision $\boldsymbol{s}_n$
already selected at $n$ unchanged. We call $n$ a \emph{revision-enabling node} when
$\rho_n=1$ and a \emph{commitment-preserving node} otherwise. An enabled revision
permits the strategic plan to change; it does not require different decisions to
be selected.

The children of a revision-enabling node may select a strategic sequence and
decide whether and when to revise it next, using the information available at
those children. The children of a commitment-preserving node must follow the
previously committed sequence and revision timing. We impose
\emph{strategic-decision linking restrictions} on the strategic decisions and
\emph{revision-policy linking restrictions} on the revision indicators to enforce
these rules across histories.

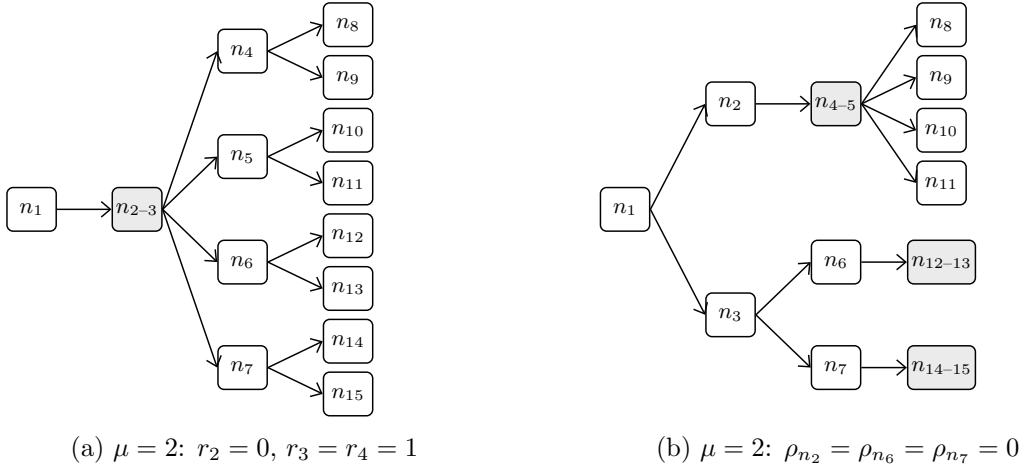
\begin{figure*}[htbp!]
    \centering
    \begin{subfigure}[c]{0.4\textwidth} 
        \resizebox{0.8\textwidth}{!}{\input{Figures/tree_stage_mu2}}
        \caption{$\mu=2$: $r_2=0$, $r_3=r_4=1$}\label{fig:sca-tree-limitedcommit}
    \end{subfigure}~\hspace{0.06\textwidth}
    \begin{subfigure}[c]{0.4\textwidth} 
        \resizebox{0.8\textwidth}{!}{\input{Figures/tree_scen_mu2}}
        \caption{$\mu=2$: $\rho_{n_2}=\rho_{n_6}=\rho_{n_7}=0$}
    \end{subfigure}
    \caption{Strategic-decision structures under (a) stage-wise and (b) scenario-adaptive revision timing with $\mu=2$. Merged boxes denote strategic decisions constrained to coincide. In panel~(b), all node indicators not specified equal one.}\label{fig:scenario-adaptive-tree}
\end{figure*}

Figure~\ref{fig:scenario-adaptive-tree} contrasts revision timing under
\bstgmsp{} and \bscenmsp{} for $\mu=2$. The stage-wise schedule uses the same two
revision stages along every scenario path, whereas the scenario-adaptive policy
allocates the two revision
opportunities according to the observed history. In panel~(b), $n_1$ is a
revision-enabling node and $n_2$ a commitment-preserving node, so the strategic
decisions for $n_4$ and $n_5$ are committed at $n_2$ using stage-2 information and
implemented identically at stage~3. Their revision indicators must also coincide;
here, $\rho_{n_4}=\rho_{n_5}=1$ opens revision opportunities at stage~4 along both
histories.

To express the linking restrictions for arbitrary pairs of nodes at the same stage,
let $\mathcal{P}(u,v)$ denote the nodes on the unique path from an ancestor $u$ to
a descendant $v$, including $u$ and excluding $v$. For two distinct stage-$t$
nodes $n_1$ and $n_2$, we define
\begin{equation}\label{eq:sca-gamma}
\Gamma_{n_1,n_2}(\boldsymbol{\rho})
:=
\sum_{n'\in
\mathcal{P}(a(n_1,n_2),n_1)\,\cup\,
\mathcal{P}(a(n_1,n_2),n_2)}\rho_{n'}.
\end{equation}
$\Gamma_{n_1,n_2}(\boldsymbol{\rho})$ is positive exactly when a revision window
opens after the most recent common ancestor of $n_1$ and $n_2$ and by stage $t$
along the path to at least one of these nodes. If it is zero, both nodes remain
bound to the same earlier commitment, so their strategic decisions and, at
non-leaf nodes, revision indicators must coincide. Otherwise, at least one node
may select these choices using information revealed after their common ancestor.
Accordingly, we model the strategic-decision linking restrictions as
\begin{subequations}\label{eq:sca-nac}
\begin{align}
&s_{n_1,i}-s_{n_2,i}\leq M_i \cdot \Gamma_{n_1,n_2}(\boldsymbol{\rho})
&& \forall t\in[2,T],\ \{n_1,n_2\}\in\Npairs_t,\ i\in\mathcal{I},
\label{eq:sca-nac-1}\\
&s_{n_2,i}-s_{n_1,i}\leq M_i \cdot \Gamma_{n_1,n_2}(\boldsymbol{\rho})
&& \forall t\in[2,T],\ \{n_1,n_2\}\in\Npairs_t,\ i\in\mathcal{I},
\label{eq:sca-nac-2}
\end{align}
\end{subequations}
and the revision-policy linking restrictions as
\begin{subequations}\label{eq:sca-rev-link}
  \begin{align}
      &\rho_{n_1}-\rho_{n_2}\leq\Gamma_{n_1,n_2}(\boldsymbol{\rho})
      && \forall t\in[2,T-1],\ \{n_1,n_2\}\in\Npairs_t,
      \label{eq:sca-rev-link-1}
      \\
      &\rho_{n_2}-\rho_{n_1}\leq\Gamma_{n_1,n_2}(\boldsymbol{\rho})
      && \forall t\in[2,T-1],\ \{n_1,n_2\}\in\Npairs_t.
      \label{eq:sca-rev-link-2}
  \end{align}
\end{subequations}
With the strategic bounds $M_i$ and binary revision indicators, each pair of linking inequalities enforces equality when $\Gamma_{n_1,n_2}(\boldsymbol{\rho})=0$ and is nonbinding when $\Gamma_{n_1,n_2}(\boldsymbol{\rho})>0$. In the latter case, revision indicators may differ even when the selected strategic decisions coincide.

Finally, \bscenmsp{} limits each scenario path to at most $\mu$ revision opportunities. For a leaf $n\in\mathcal{N}_T$, the sum of the indicators on $\mathcal{P}(1,n)$ counts the authorized revision opportunities at stages~2 through~$T$. This sum includes $\rho_1$, which governs the stage-2 opportunity; the initial stage-1 commitment consumes no budget. We impose
\begin{align}
    &\sum_{n' \in \mathcal{P}(1,n)} \rho_{n'} \;\le\; \mu
    && \forall n\in\mathcal{N}_T,
    \label{eq:path-budget}\\
    &\rho_n\in\{0,1\}
    && \forall n\in\mathcal{N}\setminus\mathcal{N}_T.
    \label{eq:sca-rho-domain}
\end{align}

\bscenmsp{} retains the objective and strategic feasibility of \bstgmsp{}
and replaces its stage-wise adaptivity constraints with the scenario-adaptive
constraints defined above. The resulting formulation is
\begin{equation*}
\begin{aligned}
(\bscenmsp{}):\quad
\min_{\substack{
{\{\boldsymbol{s}_n\}}_{n\in\mathcal{N}},\\
{\{\rho_n\}}_{n\in\mathcal{N}\setminus\mathcal{N}_T}
}}\quad
&\sum_{n\in\mathcal{N}}p_n
\Big[
f_n^{\mathrm{str}}(\boldsymbol{s}_{a(n)},\boldsymbol{s}_n)
+
\max_{\boldsymbol{d}_n\in\mathcal{D}_n}
\Big\{
\min_{\boldsymbol{o}_n\in\mathcal{O}_n(\boldsymbol{s}_n,\boldsymbol{d}_n)}
f_n^{\mathrm{opr}}(\boldsymbol{o}_n)
\Big\}
\Big] \\
\text{s.t.} \ 
&\eqref{form:pamsp-drmal},\quad\eqref{eq:sca-nac-1}-\eqref{eq:sca-nac-2},\quad\eqref{eq:sca-rev-link-1}-\eqref{eq:sca-rev-link-2},\quad\eqref{eq:path-budget}-\eqref{eq:sca-rho-domain}.
\end{aligned}
\end{equation*}

\subsection{Cost-Based, Scenario-Adaptive Framework}\label{sec:combined-framework}

The preceding subsections develop two complementary dimensions of partial
adaptivity: the revision mechanism (budgeted or cost-based) and revision timing
(stage-wise or scenario-adaptive).
We combine the commitment-cost mechanism of \cstgmsp{} with the
revision policy of \bscenmsp{} to obtain the cost-based, scenario-adaptive
formulation \cscenmsp{}.
Table~\ref{tab:cpasca} summarizes the resulting four formulations.

\begin{table}[htbp!]
    \centering
    \scalebox{0.9}{
    \begin{tabular}{lcc}
        \toprule
        \textbf{Revision timing} & \textbf{Budgeted} & \textbf{Cost-based} \\
        \midrule
        Stage-wise & \bstgmsp{} & \cstgmsp{} \\
        Scenario-adaptive & \bscenmsp{} & \cscenmsp{} \\
        \bottomrule
    \end{tabular}
    }
    \caption{The four partially adaptive MSP variants.}\label{tab:cpasca}
\end{table}

\cscenmsp{} retains the strategic-decision linking restrictions and the
revision-policy linking restrictions of \bscenmsp{} and omits its pathwise
revision budget. To value advance commitment as in \cstgmsp{}, we extend the
stage-wise commitment lag to individual nodes because revision timing may
differ across scenario paths. An authorization $\rho_u=1$ begins a new
commitment interval at the children of
$u$, at stage $t(u)+1$. We set $\delta_1(\boldsymbol{\rho}):=0$ and, for each
$n\in\mathcal{N}\setminus\{1\}$, define
\begin{equation}\label{eq:sca-lag}
\delta_n(\boldsymbol{\rho})
:=t(n)-1-\max_{u\in\mathcal{P}(1,n)}\{t(u)\rho_u\}.
\end{equation}
The lag is zero
at the root and at every child of a revision-enabling node, and increases by
one at each subsequent stage within a commitment interval. We set
$w_n(\boldsymbol{\rho}):=h(\delta_n(\boldsymbol{\rho}))$, using the same
nonincreasing function $h$ as in \cstgmsp{}, with $h(0)=1$.
Because commitment intervals may begin at different stages along different
scenario paths, nodes in the same stage may have different lags and cost
multipliers.

As in \cstgmsp{}, the multiplier applies only to commitment-eligible
strategic costs; the operational cost term remains unweighted.
The resulting formulation is
\begin{equation*}
\begin{aligned}
(\cscenmsp{}):\quad
\min_{\substack{
{\{\boldsymbol{s}_n\}}_{n\in\mathcal{N}},\\
{\{\rho_n\}}_{n\in\mathcal{N}\setminus\mathcal{N}_T}
}}\quad
&\sum_{n\in\mathcal{N}}p_n\cdot
\Big[
    w_n(\boldsymbol{\rho})f_n^{\mathrm{str}}(\boldsymbol{s}_{a(n)},\boldsymbol{s}_n)
    + 
    \max_{\boldsymbol{d}_n\in\mathcal{D}_n} 
    \Big\{
        \min_{\boldsymbol{o}_n\in\mathcal{O}_n(\boldsymbol{s}_n,\boldsymbol{d}_n)} 
        f_n^{\mathrm{opr}}(\boldsymbol{o}_n)
    \Big\}
\Big] \\
\text{s.t.} \ 
&\eqref{form:pamsp-drmal},\quad\eqref{eq:sca-nac-1}-\eqref{eq:sca-nac-2},\quad\eqref{eq:sca-rev-link-1}-\eqref{eq:sca-rev-link-2},\quad\eqref{eq:sca-rho-domain}.
\end{aligned}
\end{equation*}

Taking $h\equiv1$ and imposing the pathwise revision
budget~\eqref{eq:path-budget} yields \bscenmsp{}. For either revision
mechanism, restricting scenario-adaptive timing to a common here-and-now
schedule recovers the stage-wise formulation. The following proposition
compares the optimal values of all four formulations and the fully adaptive
benchmark.

\begin{restatable}[Optimal-value comparisons]{proposition}{optimalvaluecomparisons}\label{prop:scamsp-generalizes-pamsp}
Let $z_{\msp{}}$, $z_{\bstgmsp{}}(\mu)$, $z_{\bscenmsp{}}(\mu)$,
$z_{\cstgmsp{}}(h)$, and $z_{\cscenmsp{}}(h)$ denote the finite optimal
values of the corresponding formulations defined above. Then
\[
z_{\cscenmsp{}}(h)
\le z_{\cstgmsp{}}(h)
\le z_{\msp{}}
\le z_{\bscenmsp{}}(\mu)
\le z_{\bstgmsp{}}(\mu).
\]
Moreover,
\[
\begin{aligned}
h\equiv1 &\quad\Longrightarrow\quad
z_{\cscenmsp{}}(h)=z_{\cstgmsp{}}(h)=z_{\msp{}},\\
\mu=T-1 &\quad\Longrightarrow\quad
z_{\bscenmsp{}}(\mu)=z_{\bstgmsp{}}(\mu)=z_{\msp{}}.
\end{aligned}
\]
\end{restatable}

\subsection{Matrix Assembly Layout Specialization}\label{sec:objfunc-mald}

We now specialize the four partially adaptive MSP formulations to matrix assembly layout design, jointly optimizing layout plans and revision timing under demand uncertainty. In this section, we define the application-specific decisions and hierarchical objectives within the proposed revision-policy framework; we provide the detailed physical feasibility constraints and complete formulations in Appendix~\ref{app:formulation}.

Each stage corresponds to a production period. Let $\mathcal{M}$ denote the set of product models (e.g., variants of different car models), 
%in the automotive industry), 
$\mathcal{K}$ the set of assembly tasks, and $\mathcal{L}$ the candidate workstations at fixed locations within the matrix arrangement illustrated in Figure~\ref{fig:matrixillustration}. Assembling one unit of model $m$ requires the tasks in
$\mathcal{K}_m\subseteq\mathcal{K}$, with model-specific processing times and precedence relations.

The layout configuration at each scenario-tree node $n\in\mathcal{N}$ is represented by two binary decision families: $x_{n\ell}$ equals 1 if station $\ell$ is active, and $y_{nk\ell}$ equals 1 if task $k$ is assigned to station $\ell$; they are zero otherwise. Together, these variables constitute the strategic decision vector $\boldsymbol{s}_n=(\boldsymbol{x}_n,\boldsymbol{y}_n)$. 
We require at least one active station in each column of the matrix layout and permit task assignments only to active stations. Under forward movement across columns, precedence requires successor tasks to be assigned to the same or later columns than their predecessors. A station may perform several assigned tasks, one at a time. Assigning a task to several stations, termed task duplication, allows its workload to be distributed among them. Both components of $\boldsymbol{s}_n$ follow the strategic-decision linking restrictions of the chosen partially adaptive variant.

For demand realization $\boldsymbol{d}_n={(d_{nm})}_{m\in\mathcal{M}}\in\mathcal{D}_n$, the continuous variable $o_{nmk\ell}(\boldsymbol{d}_n)$ represents the production volume of model $m$ for which task $k\in\mathcal{K}_m$  is performed at station $\ell$. These allocations must process the realized demand for every required task while respecting the task assignments and the processing time available at each station in the period. We use these aggregate allocations to certify production feasibility in the layout design model, setting $f_n^{\mathrm{opr}}\equiv0$. That is, each committed configuration must admit a feasible production allocation for every $\boldsymbol{d}_n\in\mathcal{D}_n$, while the allocations may differ across demand realizations.

We evaluate strategic layout plans through station use and reconfiguration. Keeping a station active entails recurring costs for staffing and operating equipment and infrastructure in each production period. Task reassignments also incur implementation costs, such as equipment movement and worker training, but these are typically much smaller than the costs of maintaining additional active stations~\citep{bayraktar}. This cost structure motivates a \emph{hierarchical optimization} approach, with station use as the primary criterion and reconfiguration as the secondary criterion. We refer to the corresponding problems as the station problem (\sta{}) and reconfiguration problem (\rec{}). For the same reason, we determine the revision policy in \sta{} and fix it in \rec{}, as detailed next.

For a selected partially adaptive framework, let $\boldsymbol{\eta}$ denote the stage-wise schedule $\boldsymbol{r}$ or scenario-adaptive policy $\boldsymbol{\rho}$. Let $\mathcal{A}$ be the set of pairs $(\boldsymbol{s},\boldsymbol{\eta})$ satisfying the framework's strategic-decision linking restrictions, revision-variable domains, and applicable revision-policy linking and budget constraints. We write the objectives in a common form across all four variants, using $w_n(\boldsymbol{\eta})$ to weight the contribution of node $n$. In the budgeted variants, $w_n(\boldsymbol{\eta})=1$ for all $n\in\mathcal{N}$; in the cost-based variants, it is the commitment-dependent cost multiplier defined earlier. 

The station model is
\begin{equation}\label{form:sta-generic}
\begin{aligned}
    (\text{\sta{}}):\quad
    \min_{\boldsymbol{s},\boldsymbol{o},\boldsymbol{\eta}}\quad &
    \sum_{n\in\mathcal{N}}p_n\cdot w_n(\boldsymbol{\eta}) \cdot
    \sum_{\ell\in\mathcal{L}}x_{n\ell}\\
    \text{s.t.}\quad &
    (\boldsymbol{s}_{a(n)},\boldsymbol{s}_n)
    \in\mathcal{S}_n
    && \forall n\in\mathcal{N},\\
    &\boldsymbol{o}_n(\boldsymbol{d}_n)
    \in\mathcal{O}_n(\boldsymbol{s}_n,\boldsymbol{d}_n)
    && \forall n\in\mathcal{N},\
       \boldsymbol{d}_n\in\mathcal{D}_n,\\
    &(\boldsymbol{s},\boldsymbol{\eta})\in\mathcal{A}.
\end{aligned}
\end{equation}
In the objective, the inner sum counts the active stations at node $n$. With unit weights, the objective thus measures expected station use summed across production periods, counting a station in every period in which it is active. With commitment-dependent weights, it represents discounted station-use cost. Note that the operational variables enter the model only through the feasibility constraints, which must hold for every node and demand realization. If $\mathcal{O}_n(\boldsymbol{s}_n,\boldsymbol{d}_n)=\emptyset$ for any such pair, the corresponding layout is infeasible.

From a selected optimal \sta{} solution, let $\boldsymbol{\eta}^*$ denote its revision decisions and set $\zsta{n}:=\sum_{\ell\in\mathcal{L}}x_{n\ell}^*$ for all $n\in\mathcal{N}$. The reconfiguration model is
\begin{equation}\label{form:rec-generic}
\begin{aligned}
    (\text{\rec{}}):\quad
    \min_{\boldsymbol{s},\boldsymbol{o}}\quad &
    \sum_{n\in\mathcal{N}}p_n\cdot w_n(\boldsymbol{\eta}^*) \cdot
    \sum_{k\in\mathcal{K}}\sum_{\ell\in\mathcal{L}}
    {[y_{nk\ell}-y_{a(n),k\ell}]}^+\\
    \text{s.t.}\quad &
    (\boldsymbol{s}_{a(n)},\boldsymbol{s}_n)
    \in\mathcal{S}_n
    && \forall n\in\mathcal{N},\\
    &\boldsymbol{o}_n(\boldsymbol{d}_n)
    \in\mathcal{O}_n(\boldsymbol{s}_n,\boldsymbol{d}_n)
    && \forall n\in\mathcal{N},\
       \boldsymbol{d}_n\in\mathcal{D}_n,\\
    &(\boldsymbol{s},\boldsymbol{\eta}^*)\in\mathcal{A}\\
    &\sum_{\ell\in\mathcal{L}}x_{n\ell}\leq\zsta{n}
    && \forall n\in\mathcal{N},
\end{aligned}
\end{equation}
where ${[v]}^+:=\max\{v,0\}$. In the objective, ${[y_{nk\ell}-y_{a(n),k\ell}]}^+$ equals one when task $k$ is assigned to station $\ell$ at node $n$ but not at its parent. Thus, each new task-station assignment counts as one reconfiguration, including additions prescribed within committed layout sequences. In our numerical study, we assume a greenfield start with no active stations or task assignments before the first period, so we set $y_{a(1),k\ell}=0$ and count initial assignments in the objective as well. With unit weights, the objective measures the expected number of reconfigurations; with
commitment-dependent weights, it represents discounted reconfiguration cost.

In \rec{}, the revision decisions $\boldsymbol{\eta}^*$ and associated weights $w_n(\boldsymbol{\eta}^*)$ are fixed, and the last constraints bound station use at each node by $\zsta{n}$. The following result justifies these bounds.

\begin{restatable}[Preservation of primary optima]{proposition}{primaryoptima}\label{prop:rec-primary-optimality}
For each of the four partially adaptive layout variants, a pair $(\boldsymbol{s},\boldsymbol{o})$ is feasible for~\eqref{form:rec-generic} if and only if $(\boldsymbol{s},\boldsymbol{o},\boldsymbol{\eta}^*)$ is optimal for~\eqref{form:sta-generic}. Every such pair satisfies $\sum_{\ell\in\mathcal{L}}x_{n\ell}=\zsta{n}$ for all $n\in\mathcal{N}$. Consequently, an optimal \rec{} solution minimizes reconfiguration cost over all \sta{}-optimal layouts with revision decisions $\boldsymbol{\eta}^*$.
\end{restatable}

In what follows, we name the variants by replacing the final \texttt{MSP} tag in Table~\ref{tab:cpasca} with \texttt{STA} or \texttt{REC}. For example, \cscensta{} and \cscenrec{} denote the cost-based, scenario-adaptive station and reconfiguration problems, respectively.

% ---------------------------------- %
\section{Solution Approach} \label{sec:solution}

\begin{figure}[t]
    \centering
    \scalebox{0.8}{\input{Figures/SolutionApproachOverview.tex}}
     %\scalebox{0.85}{\input{Figures/SolutionApproachOverviewCompact.tex}}
    \caption{Solution approach for the four partially adaptive layout variants.}
    \label{fig:solappover}
\end{figure}
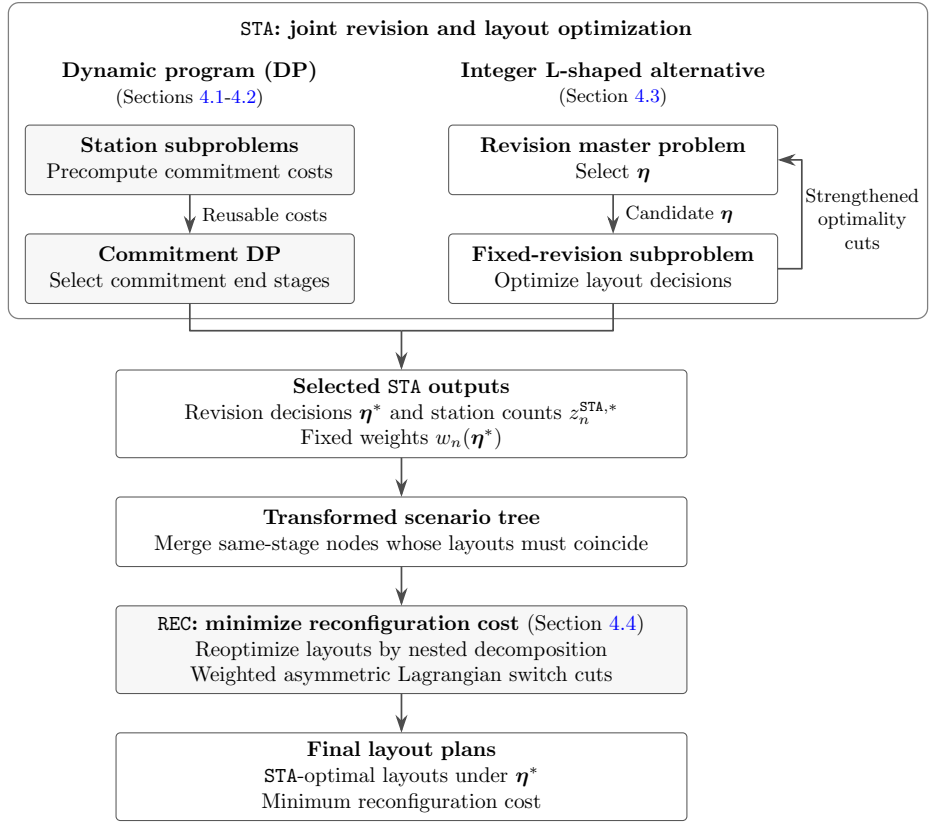

This section presents our solution approach for the four partially adaptive
layout variants, which is summarized in Figure~\ref{fig:solappover}. We first
solve \sta{}, jointly optimizing the revision policy and layouts. For a fixed
revision policy, \sta{} decomposes into independent station subproblems whose
optimal costs can be reused across revision policies (Section~\ref{sec:structure}).
Our main method precomputes these costs and combines them in an exact commitment
dynamic program to determine an optimal revision policy (Section~\ref{sec:dp}).
We present an integer L-shaped alternative that also accommodates costs or
constraints coupling consecutive stages. It evaluates the fixed-revision
subproblem under candidate revision policies to generate cuts for the revision
master problem. We strengthen this method with monotonicity-based cuts for the
budgeted variants and enhanced cuts exploiting \sta{}'s separability for all
four variants (Section~\ref{sec:lshaped}). We then solve \rec{}
with the selected revision decisions fixed and the station-count bounds inherited
from \sta{}, reoptimizing layouts to minimize reconfiguration cost. As these
costs couple consecutive layouts, we use nested decomposition on a transformed
scenario tree. To strengthen this procedure, we develop weighted asymmetric
Lagrangian switch cuts for commitment-dependent reconfiguration costs
(Section~\ref{sec:lexicographic}). All proofs are provided in
Appendix~\ref{app:proofs}.

\subsection{Commitment Structure and Cost Decomposition} \label{sec:structure}

We first characterize the commitment structure induced by a revision policy
and derive a cost decomposition, providing the basis for the dynamic program
developed in Section~\ref{sec:dp}. For a fixed revision policy, we identify the groups of same-stage nodes
whose strategic decisions must coincide and exploit the separability of
\sta{} to express its optimal cost as the sum of the costs obtained by
optimizing each group independently.

To formalize these results, we use scenario-adaptive
notation, with $\boldsymbol{\rho}$ describing the revision policy. The stage-wise
variants correspond to policies whose revision indicators are identical across
nodes at the same stage, under the mapping $\rho_n:=r_{t(n)+1}$ for each
non-leaf node, as in the proof of Proposition~\ref{prop:scamsp-generalizes-pamsp}.
Thus, the results below also apply to the stage-wise variants.

For a node $n$ and a stage $t\in[t(n),T]$, let $\mathcal{C}_t(n)$ denote its
stage-$t$ descendants. Thus, $\mathcal{C}_{t(n)}(n)=\{n\}$ and
$\mathcal{C}_{t(n)+1}(n)=\mathcal{C}(n)$; we also set
$\mathcal{C}_{T+1}(n):=\emptyset$.
For a given $\boldsymbol{\rho}$, we define the set of \emph{commitment-start nodes} as
\[
\mathcal{B}(\boldsymbol{\rho}):=\{1\} \hspace*{0.2cm} \cup \hspace*{-0.4cm}
\bigcup_{n\in\mathcal{N}\setminus\mathcal{N}_T:\,\rho_n=1} \hspace*{-0.55cm} \mathcal{C}(n).
\]
For each node $n$, let $\sigma_n(\boldsymbol{\rho})$ denote the most recent
commitment-start node on the root-to-$n$ path, including $n$ itself. Thus,
$\sigma_n(\boldsymbol{\rho})$ identifies the beginning of the commitment
interval containing $n$.
For each $u\in\mathcal{B}(\boldsymbol{\rho})$, we also define the \emph{end stage}
of the commitment beginning at $u$ as
$e_u(\boldsymbol{\rho}):=\max\{t(n):n\in\mathcal{N},
\sigma_n(\boldsymbol{\rho})=u\}$, which lies in $[t(u),T]$.
If $\boldsymbol{\rho}$ satisfies~\eqref{eq:sca-rev-link}, then for each
$u\in\mathcal{B}(\boldsymbol{\rho})$ and
$t\in[t(u),e_u(\boldsymbol{\rho})]$, we call $\mathcal{C}_t(u)$ the
\emph{stage-$t$ commitment group} starting at $u$. Part~(c) of
Lemma~\ref{lem:commitment-structure} shows that its nodes have $u$ as
their most recent commitment start.
For example, in Figure~\ref{fig:scenario-adaptive-tree}(b),
$\mathcal{B}(\boldsymbol{\rho})=\{n_1,n_2,n_3,n_6,\ldots,n_{11}\}$,
$\sigma_{n_2}(\boldsymbol{\rho})=\sigma_{n_4}(\boldsymbol{\rho})
=\sigma_{n_5}(\boldsymbol{\rho})=n_2$ with $e_{n_2}(\boldsymbol{\rho})=3$, and
$\sigma_{n_6}(\boldsymbol{\rho})=\sigma_{n_{12}}(\boldsymbol{\rho})
=\sigma_{n_{13}}(\boldsymbol{\rho})=n_6$ with $e_{n_6}(\boldsymbol{\rho})=4$.
Thus, $\mathcal{C}_2(n_2)=\{n_2\}$ and
$\mathcal{C}_3(n_2)=\{n_4,n_5\}$ are distinct commitment groups.

\begin{restatable}[Commitment structure]{lemma}{commitmentstructure} \label{lem:commitment-structure}
Let $\boldsymbol{\rho}$ be a binary vector indexed by the non-leaf nodes.
\begin{enumerate}[label=(\alph*)]
    \item For any $t\in[2,T]$ and $\{n_1,n_2\}\in\Npairs_t$, $\Gamma_{n_1,n_2}(\boldsymbol{\rho})=0$ if and only if $\sigma_{n_1}(\boldsymbol{\rho})=\sigma_{n_2}(\boldsymbol{\rho})$.
    \item $\delta_n(\boldsymbol{\rho})=t(n)-t(\sigma_n(\boldsymbol{\rho}))$ for every $n\in\mathcal{N}$.
    \item If $\boldsymbol{\rho}$ satisfies~\eqref{eq:sca-rev-link}, then for every $u\in\mathcal{B}(\boldsymbol{\rho})$, the nodes in $\mathcal{C}_t(u)$ are commitment-preserving for $t(u)\le t<e_u(\boldsymbol{\rho})$, the nodes in $\mathcal{C}_{e_u(\boldsymbol{\rho})}(u)$ are revision-enabling if $e_u(\boldsymbol{\rho})<T$, and
    \[
        \{n\in\mathcal{N}:\sigma_n(\boldsymbol{\rho})=u\}=\bigcup_{t=t(u)}^{e_u(\boldsymbol{\rho})}\mathcal{C}_t(u).
    \]
\end{enumerate}
\end{restatable}

For any $\boldsymbol{\rho}$ satisfying~\eqref{eq:sca-rev-link}, part~(c) of
Lemma~\ref{lem:commitment-structure} partitions the nodes into
stage-$t$ commitment groups $\mathcal{C}_t(u)$.
Strategic decisions are linked within each group; the following property
of \sta{} allows us to optimize the groups independently.

\begin{property}[Node-local strategic costs and feasibility] \label{property:node-local}
In all four \sta{} variants, changing the preceding layout while holding the
current layout fixed leaves the strategic cost and physical feasibility
unchanged. Specifically, for every $n\in\mathcal{N}$, candidate current
strategic vector $\bar{\boldsymbol{s}}\in\mathbb{R}^{|\mathcal{I}|}$, and
candidate preceding strategic vectors
$\bar{\boldsymbol{s}}^{\mathrm{prev}},\widetilde{\boldsymbol{s}}^{\mathrm{prev}}\in\mathbb{R}^{|\mathcal{I}|}$,
\[
\begin{aligned}
    f_n^{\mathrm{str}}(\bar{\boldsymbol{s}}^{\mathrm{prev}},\bar{\boldsymbol{s}})
    &=f_n^{\mathrm{str}}(\widetilde{\boldsymbol{s}}^{\mathrm{prev}},\bar{\boldsymbol{s}}),\\
    (\bar{\boldsymbol{s}}^{\mathrm{prev}},\bar{\boldsymbol{s}})\in\mathcal{S}_n
    &\quad\Longleftrightarrow\quad
    (\widetilde{\boldsymbol{s}}^{\mathrm{prev}},\bar{\boldsymbol{s}})\in\mathcal{S}_n.
\end{aligned}
\]
\end{property}

This property follows directly from the \sta{} formulation:
the strategic cost at node $n$ is $\sum_{\ell\in\mathcal{L}}x_{n\ell}$,
and the constraints defining $\mathcal{S}_n$,
namely~\eqref{app:sta-col}-\eqref{app:sta-prec} and~\eqref{app:sta-bin},
involve only $\boldsymbol{x}_n$ and $\boldsymbol{y}_n$.

For formulations satisfying Property~\ref{property:node-local}, we use
$\bar{\boldsymbol{s}}$ for a candidate strategic vector at a single node
and write $f_n^{\mathrm{str}}(\bar{\boldsymbol{s}})$ and
$\bar{\boldsymbol{s}}\in\mathcal{S}_n$ for
the strategic cost and feasibility condition, respectively. Accounting for
operational feasibility and cost, we define the feasible strategic set and
local cost at node $n$, for a lag $\delta$, as
\[
\begin{aligned}
    \mathcal{X}_n&:=\big\{\bar{\boldsymbol{s}}:\bar{\boldsymbol{s}}\in\mathcal{S}_n,\ \mathcal{O}_n(\bar{\boldsymbol{s}},\boldsymbol{d}_n)\neq\emptyset\ \ \forall\boldsymbol{d}_n\in\mathcal{D}_n\big\},\\
    \varphi_n(\bar{\boldsymbol{s}},\delta)&:=h(\delta)f_n^{\mathrm{str}}(\bar{\boldsymbol{s}})+\max_{\boldsymbol{d}_n\in\mathcal{D}_n}\min_{\boldsymbol{o}_n\in\mathcal{O}_n(\bar{\boldsymbol{s}},\boldsymbol{d}_n)}f_n^{\mathrm{opr}}(\boldsymbol{o}_n).
\end{aligned}
\]
For the budgeted variants, we take $h\equiv1$.
For a commitment interval beginning
at $u$ and reaching stage $t\in[t(u),T]$, its stage-$t$ commitment group
$\mathcal{C}_t(u)$ must share a strategic vector feasible at every node.
Its nodes also share the commitment lag $t-t(u)$. We refer to this stage
contribution as a \emph{unit}, indexed by $(u,t)$. The associated
\emph{unit problem} minimizes
the group's probability-weighted cost over the common feasible strategic
vector, with optimal value
\begin{equation}
    c(u,t):=\min\Big\{\sum_{n\in\mathcal{C}_t(u)}p_n\,\varphi_n\big(\bar{\boldsymbol{s}},t-t(u)\big):\ \bar{\boldsymbol{s}} \ \in \hspace*{-0.2cm} \bigcap_{n\in\mathcal{C}_t(u)} \hspace*{-0.2cm} \mathcal{X}_n\Big\}.
    \label{eq:unit-cost}
\end{equation}
Throughout the paper, we assign value $+\infty$ to infeasible minimization
problems and assume that every feasible operational, unit, or fixed-revision
subproblem admits an optimal solution with finite objective value.
For \sta{}, this assumption holds automatically:
each feasible unit or fixed-revision subproblem minimizes over finitely many binary layouts,
and every feasible operational response has zero cost.
For fixed model data and $h$, the unit cost depends only on
$u$ and $t$; it is independent of the selected revision policy and the end
stage of the commitment interval, provided the interval reaches $t$.
The following proposition shows that, under a fixed revision policy, the unit
problems can be optimized independently; 
Section~\ref{sec:dp} uses this decomposition to optimize the revision policy.

\begin{restatable}[Interval decomposition]{proposition}{intervaldecomposition} \label{prop:interval-decomp}
For a formulation satisfying Property~\ref{property:node-local}, let $\boldsymbol{\rho}$ be a binary vector satisfying the revision-policy linking restrictions~\eqref{eq:sca-rev-link} and the applicable revision-budget constraints. Let $v(\boldsymbol{\rho})$ denote the optimal value of the formulation when $\boldsymbol{\rho}$ is fixed. Then
\begin{equation}
    v(\boldsymbol{\rho})=\sum_{u\in\mathcal{B}(\boldsymbol{\rho})}\ \sum_{t=t(u)}^{e_u(\boldsymbol{\rho})}c(u,t).
    \label{eq:interval-decomp}
\end{equation}
\end{restatable}

We conclude this subsection by showing how to compute the unit costs for
\sta{} from station problems for node groups. For any nonempty node group
$\mathcal{N}_g\subseteq\mathcal{N}$, let $\mathrm{SP}^*(\mathcal{N}_g)$ denote the minimum number
of active stations in a common layout
$\bar{\boldsymbol{s}}=(\bar{\boldsymbol{x}},\bar{\boldsymbol{y}})$ feasible
at every node in $\mathcal{N}_g$:
\begin{equation}
    \mathrm{SP}^*(\mathcal{N}_g):=\min\Big\{\sum_{\ell\in\mathcal{L}}\bar{x}_\ell:\ \bar{\boldsymbol{s}}\in\mathcal{S}_n\ \forall n\in\mathcal{N}_g,\ \mathcal{O}_n(\bar{\boldsymbol{s}},\boldsymbol{d}_n)\neq\emptyset\ \forall n\in\mathcal{N}_g,\ \boldsymbol{d}_n\in\mathcal{D}_n\Big\}.
    \label{eq:station-problem}
\end{equation}
We set $\mathrm{SP}^*(\mathcal{N}_g)=+\infty$ when no such layout exists. Using
$f_n^{\mathrm{opr}}\equiv0$ and
$\sum_{n\in\mathcal{C}_t(u)}p_n=p_u$ for $t\in[t(u),T]$,
the unit cost reduces to
\begin{equation}
    c(u,t)=p_u\,h\big(t-t(u)\big)\,\mathrm{SP}^*\big(\mathcal{C}_t(u)\big).
    \label{eq:unit-cost-sta}
\end{equation}
%\begin{remark}[Unit precomputation] \label{rem:dp-implementation} The station problems~\eqref{eq:station-problem} for all unit indices can be obtained by evaluating the units of $T$ stage-wise schedules separately. Full commitment covers units $(1,t)$ for $t\in[1,T]$; a schedule with a single revision at stage $t\in[2,T]$ covers units $(u,t')$ for $u\in\mathcal{N}_t$ and $t'\in[t,T]$. \end{remark}
A \emph{key computational advantage} is that the optimal station
counts can be computed once and reused across all four variants, budgets,
and functions $h$, since $h$ enters only through the positive scalar multiplier
in~\eqref{eq:unit-cost-sta}.

\subsection{Commitment Dynamic Program} \label{sec:dp}

We develop a dynamic program (DP) to optimize revision policies using the
precomputed unit costs of Section~\ref{sec:structure}. For scenario-adaptive
policies, it optimizes the duration of a commitment beginning at a given
node, using optimal policies computed recursively from the nodes where
subsequent commitments begin. We first present these recursions, then
adapt them to the stage-wise variants.

For \cscenmsp{}, let $V(u)$ denote the cost-to-go value at $u\in\mathcal{N}$:
the minimum contribution of its subtree to the full objective, given that a
commitment begins at $u$. Choosing $e\in[t(u),T]$ as its end stage incurs the interval
cost $\sum_{t=t(u)}^e c(u,t)$. If $e<T$, a new commitment begins at each
$u'\in\mathcal{C}_{e+1}(u)$, and these successor subtrees can be optimized
independently. Thus,
\begin{equation}
    V(u)=\min_{e\in[t(u),T]}\Big\{\sum_{t=t(u)}^{e}c(u,t)  \ + \hspace*{-0.2cm} \sum_{u'\in\mathcal{C}_{e+1}(u)} \hspace*{-0.25cm}  V(u')\Big\}, \label{eq:dp-cscen}
\end{equation}
where $\mathcal{C}_{T+1}(u)=\emptyset$ makes the cost-to-go sum zero when
$e=T$. In particular, $V(u)=c(u,T)$ for every leaf $u\in\mathcal{N}_T$.
The unit costs already include the node probabilities, so successor
cost-to-go values are summed without additional probability weights.

For \bscenmsp{}, the state also records the remaining pathwise revision
budget. Let $V(u,b)$ denote the corresponding cost-to-go value with $b\in[0,\mu]$
revision opportunities remaining on each path through $u$. Ending the current
interval at $e<T$ uses one opportunity on every such path, leaving $b-1$
for each successor commitment. Hence, for every $u\in\mathcal{N}$,
\begin{equation}
    V(u,b)=\min_{e\in[t(u),T]}\Big\{\sum_{t=t(u)}^{e}c(u,t) \ + \hspace*{-0.2cm} \sum_{u'\in\mathcal{C}_{e+1}(u)} \hspace*{-0.25cm} V(u',b-1)\Big\}, \quad b\in[1,\mu]. \label{eq:dp-bscen}
\end{equation}
For $b=0$, commitment lasts through stage $T$, so
$V(u,0):=\sum_{t=t(u)}^T c(u,t)$.

In the stage-wise variants, all commitments beginning at stage $t\in[1,T]$
must share an end stage. Let $W(t)$ denote the optimal cost-to-go value for stages
$t,\ldots,T$ in \cstgmsp{} when a commitment begins at every node in
$\mathcal{N}_t$, and let $W(t,b)$ denote its \bstgmsp{} counterpart with
$b\in[0,\mu]$ revisions remaining. The interval costs now sum over
$u\in\mathcal{N}_t$, and choosing end stage $e<T$ starts new commitments at
every node in $\mathcal{N}_{e+1}$. The recursions are
\begin{align}
    W(t)&=\min_{e\in[t,T]}\Big\{\sum_{u\in\mathcal{N}_t}\sum_{t'=t}^{e}c(u,t')+W(e+1)\Big\}, \label{eq:dp-cstg}\\
    W(t,b)&=\min_{e\in[t,T]}\Big\{\sum_{u\in\mathcal{N}_t}\sum_{t'=t}^{e}c(u,t')+W(e+1,b-1)\Big\}, \quad b\in[1,\mu], \label{eq:dp-bstg}
\end{align}
with $W(t,0):=\sum_{u\in\mathcal{N}_t}\sum_{t'=t}^T c(u,t')$ and terminal values
$W(T+1)=W(T+1,b)=0$ for all $b\in[0,\mu]$.

Starting from the boundary values above, we evaluate all states from
stage $T$ backward to stage~1. For each state, we compare the admissible
end stages using the interval costs and the already computed
cost-to-go values of successor commitments, retaining a minimizing
end stage for policy recovery. The following theorem establishes
the exactness of this procedure.

\begin{restatable}[Exactness of the commitment DP]{theorem}{commitmentdpexactness} \label{thm:dp}
Suppose Property~\ref{property:node-local} holds and each operational and unit
problem is either infeasible or admits an optimal solution with finite
objective value. Then the optimal values of the four partially adaptive
formulations satisfy
\[
    z_{\cscenmsp{}}(h)=V(1), \ \ 
    z_{\bscenmsp{}}(\mu)=V(1,\mu),\ \
    z_{\cstgmsp{}}(h)=W(1),\ \ 
    z_{\bstgmsp{}}(\mu)=W(1,\mu).
\]
For each formulation, if the corresponding DP value above is finite, an optimal
revision policy and associated strategic decisions can be recovered by following
minimizing end stages from stage~1 and combining optimal solutions of the
corresponding unit problems~\eqref{eq:unit-cost}.
A DP value of $+\infty$ certifies infeasibility.
\end{restatable}

The \emph{computational advantage} of the DP is that it optimizes the revision policy
using polynomially many operations once a polynomial number of unit
problems has been solved. There is one unit problem for each $u\in\mathcal{N}$
and $t\in[t(u),T]$, giving
$|\mathcal{N}|+\sum_{u\in\mathcal{N}}(T-t(u))$ unit problems, at most
$T|\mathcal{N}|$ in total. Once their costs are available, the
scenario-adaptive recursions require $O(T|\mathcal{N}|)$ arithmetic
operations in \cscenmsp{} and $O((\mu+1)T|\mathcal{N}|)$ in \bscenmsp{}.
For comparison, \cscenmsp{} with $T\ge2$ has at least
$2^{|\mathcal{N}_{T-1}|}$ binary revision policies satisfying the
revision-policy linking restrictions~\eqref{eq:sca-rev-link}.

\begin{remark}[Extensions] \label{rem:dp-extensions}
The DP also accommodates more general cost weights and restrictions on revision
stages. The strategic-cost multiplier in the unit problems~\eqref{eq:unit-cost}
may be any given positive function $w(t(u),t)$ of the commitment's start and
current stages without changing the form of the recursions. Also, we can
restrict revisions to specified stages by allowing only commitment end stages
immediately before a permitted revision stage, or $T$ for commitment through
the horizon.
In our numerical experiments in Section~\ref{sec:granularity},
we use these extensions to compare monthly, quarterly, and yearly revision
calendars on scenario trees with monthly stages, with cost weights based on
commitment lag measured in quarters.
\end{remark}

\subsection{Integer L-Shaped Decomposition and Optimality Cuts} \label{sec:lshaped}

As an alternative to the DP, we present an integer L-shaped
method~\citep{LAPORTE1993133} that optimizes revision policies through
repeated fixed-policy evaluations. A revision master problem selects a
candidate policy, and a fixed-revision subproblem optimizes the strategic
and operational decisions under that policy. Its optimal cost yields cuts
that guide subsequent policy choices. This approach also accommodates
interstage coupling without requiring Property~\ref{property:node-local}.
We first present the general procedure, then strengthen its cuts using
monotonicity in the general budgeted partially adaptive formulations and, for \sta{},
the commitment structure shared by all four variants.

We continue with scenario-adaptive notation; the corresponding results
also apply to the stage-wise variants under the mapping introduced earlier. 
We call a revision policy
$\boldsymbol{\rho}$ \emph{admissible} if it satisfies the revision-policy
linking restrictions~\eqref{eq:sca-rev-link} and binary
domains~\eqref{eq:sca-rho-domain}, together with the pathwise
budget~\eqref{eq:path-budget} in the budgeted variant.
For an admissible policy, the fixed-revision subproblem optimizes the
remaining strategic and operational decisions under the associated
strategic-decision linking restrictions and cost weights.
As in Section~\ref{sec:structure}, we denote its optimal value by
$v(\boldsymbol{\rho})$.
For \sta{}, the evaluation reduces to independent station
problems~\eqref{eq:station-problem}, whose weighted values are combined
as in~\eqref{eq:interval-decomp}.

The revision master problem introduces an epigraph variable $\theta$
for a lower approximation of $v(\boldsymbol{\rho})$. Given a finite bound
$L\le v(\boldsymbol{\rho})$ for every admissible policy,
it solves, at iteration $j$,
\begin{equation}
    \min\big\{\theta:\ \boldsymbol{\rho}\text{ admissible},\ \theta\ge L,\
    \text{cuts generated in iterations }1,\ldots,j-1\big\}.
    \label{eq:master}
\end{equation}
For \sta{}, nonnegative station costs allow $L=0$.
The optimal revision master problem value bounds the overall optimum
from below. Evaluating its candidate policy
$\bar{\boldsymbol{\rho}}^j$ gives
$v^j:=v(\bar{\boldsymbol{\rho}}^j)$; the best finite value found gives
an upper bound. If the fixed-revision subproblem is infeasible, we can exclude
the candidate by adding the no-good cut
\[
    \sum_{\substack{n\in\mathcal{N}\setminus\mathcal{N}_T \, : \, \bar{\rho}_n^j=1}} \hspace*{-0.45cm}  (1-\rho_n)
    \ + \hspace*{-0.45cm}  \sum_{\substack{n\in\mathcal{N}\setminus\mathcal{N}_T \, : \, \bar{\rho}_n^j=0}}\hspace*{-0.45cm}  \rho_n\ge1,
\]
otherwise, we can add the integer optimality cut (IOC)~\citep{LAPORTE1993133}
\begin{equation}
    \theta\ge L+(v^j-L)\Big(1 \ -
    \hspace*{-0.45cm} \sum_{\substack{n\in\mathcal{N}\setminus\mathcal{N}_T  \, : \, \bar{\rho}_n^j=1}}  \hspace*{-0.45cm} (1-\rho_n)
    \ -  \hspace*{-0.45cm} \sum_{\substack{n\in\mathcal{N}\setminus\mathcal{N}_T \, : \, \bar{\rho}_n^j=0}}  \hspace*{-0.45cm} \rho_n\Big).
    \label{eq:standard-cut}
\end{equation}
The right-hand side equals $v^j$ at the evaluated policy and is at most $L$
at every other binary policy. Thus, the IOC is valid and tight at
$\bar{\boldsymbol{\rho}}^j$.
For the stage-wise variants, the no-good cut and IOC use
$r_t$ and $\bar{r}_t^j$, with one term per revision stage $t\in[2,T]$.
We re-solve the revision master problem with all generated cuts until
the upper and lower bounds differ by at most a prescribed tolerance
$\varepsilon\ge0$ or the revision master problem is infeasible.
Since each iteration evaluates a new policy or terminates,
the finite policy set guarantees finite termination.

We next strengthen the IOC for the budgeted partially adaptive MSP class
by exploiting monotonicity in revision authorizations. Removing these
authorizations restricts strategic decisions
without changing their costs, so the evaluated value remains a lower bound
for admissible policies obtained in this way.
\begin{restatable}[Monotonicity]{proposition}{revisionmonotonicity} \label{prop:monotonic}
In \bscenmsp{}, let $\boldsymbol{\rho}$ and $\boldsymbol{\rho}'$ be
admissible revision policies satisfying
$\boldsymbol{\rho}\le\boldsymbol{\rho}'$ componentwise. Then
$v(\boldsymbol{\rho}')\le v(\boldsymbol{\rho})$.
The result also holds for \bstgmsp{}.
\end{restatable}

For a finite evaluation $v^j$, Proposition~\ref{prop:monotonic} gives
$v(\boldsymbol{\rho})\ge v^j$ for every admissible policy with
$\boldsymbol{\rho}\le\bar{\boldsymbol{\rho}}^j$. Following~\citet{Hooker2024},
we enforce this bound through the \emph{monotone optimality cut} (MOC)
\begin{equation}
    \theta\ge L+(v^j-L)\Big(1 \ -
    \hspace*{-0.45cm} \sum_{\substack{n\in\mathcal{N}\setminus\mathcal{N}_T \, : \, \bar{\rho}_n^j=0}} \hspace*{-0.45cm} \rho_n\Big).
    \label{eq:monotonic-cut:sn}
\end{equation}
% or, using stage-wise indicators in \bstgmsp{},
% \begin{equation}
%     \theta\ge L+(v^j-L)\Big(1 \ - \hspace*{-0.45cm} \sum_{t\in[2,T] \, : \, \bar{r}_t^j=0} \hspace*{-0.45cm} r_t\Big).
%     \label{eq:monotonic-cut:st}
% \end{equation}
The right-hand side equals $v^j$ when no new authorization is added to
the evaluated policy, and is at most $L$ otherwise. The MOC therefore
strengthens the IOC by retaining the bound when some evaluated
authorizations are removed.
We note that this monotonicity need not hold in the cost-based variants
since an additional authorization can shorten commitment lags and increase
cost multipliers.

For all four \sta{} variants, the decomposition in
Section~\ref{sec:structure} permits further strengthening by bounding
each node's contribution separately. This contribution depends only on
the node's commitment start, so its bound can remain informative when
the policy changes elsewhere in the tree.
For an admissible policy $\boldsymbol{\rho}$ and $n\in\mathcal{N}$, let
\[
    v_n(\boldsymbol{\rho}):=w_n(\boldsymbol{\rho})\,
    \mathrm{SP}^*\big(\mathcal{C}_{t(n)}(\sigma_n(\boldsymbol{\rho}))\big),
\]
%with $w_n\equiv1$ in the budgeted variants. 
be the optimal weighted
station count at $n$. Then, \eqref{eq:interval-decomp}
and~\eqref{eq:unit-cost-sta} give
$v(\boldsymbol{\rho})=\sum_{n\in\mathcal{N}}p_n v_n(\boldsymbol{\rho})$.
We introduce a lower approximation $\theta_n\ge0$ for each
$v_n(\boldsymbol{\rho})$ and replace $\theta$ throughout the revision
master problem~\eqref{eq:master} by
$\sum_{n\in\mathcal{N}}p_n\theta_n$, using $L=0$.
For an evaluated policy with $v^j<\infty$, let
$v_n^j:=v_n(\bar{\boldsymbol{\rho}}^j)$ and
$\sigma_n^j:=\sigma_n(\bar{\boldsymbol{\rho}}^j)$, the commitment-start node
for $n$ under the evaluated policy.

If the indicators on the path $\mathcal{P}(\sigma_n^j,n)$ remain zero, the commitment
start for $n$ stays at $\sigma_n^j$ or moves to an ancestor.
Its commitment group therefore stays the same or expands, preserving
the lower bound $v_n^j$ in a budgeted variant.
In a cost-based variant, we also retain the authorization at the parent 
$a(\sigma_n^j)$ when $\sigma_n^j\neq1$, thereby keeping the commitment
start and its cost multiplier unchanged.
These observations yield the \emph{enhanced optimality cuts}
(EOCs), for every $n\in\mathcal{N}$,
\begin{align}
    \theta_n&\ge v_n^j\Big(1 \ - \hspace*{-0.45cm} \sum_{u\in\mathcal{P}(\sigma_n^j,n)} \hspace*{-0.45cm} \rho_u\Big)
    &&\text{in \bscensta{}}, \label{eq:enhanced-cut:bpa-sca}\\
    \theta_n&\ge v_n^j\Big(1 \ - (1-\rho_{a(\sigma_n^j)}) \ -
    \hspace*{-0.45cm} \sum_{u\in\mathcal{P}(\sigma_n^j,n)} \hspace*{-0.45cm} \rho_u\Big)
    &&\text{in \cscensta{}}, \label{eq:enhanced-cut:cpa-sca}
\end{align}
where $1-\rho_{a(\sigma_n^j)}$ is omitted when $\sigma_n^j=1$.
The following result establishes the validity of the EOCs and their
dominance over the earlier cut families.

\begin{restatable}[Validity and dominance of EOCs]{proposition}{eocvalidity} \label{prop:valid-separable-eoc}
In \bscensta{} or \cscensta{}, let $\bar{\boldsymbol{\rho}}^j$ be an
admissible revision policy with $v^j<\infty$. Then:
\begin{enumerate}[label=(\alph*)]
    \item For every admissible policy $\boldsymbol{\rho}$ with
    $v(\boldsymbol{\rho})<\infty$ and every $n\in\mathcal{N}$,
    the corresponding EOC holds at $\theta_n=v_n(\boldsymbol{\rho})$,
    with equality when $\boldsymbol{\rho}=\bar{\boldsymbol{\rho}}^j$.
    \item Under $\theta=\sum_{n\in\mathcal{N}}p_n\theta_n$, the EOCs
    imply the MOC~\eqref{eq:monotonic-cut:sn}
    in \bscensta{} and the IOC~\eqref{eq:standard-cut} in \cscensta{},
    generated at $\bar{\boldsymbol{\rho}}^j$ with $L=0$.
\end{enumerate}
The same validity and dominance properties hold for \bstgsta{} and \cstgsta{}
under the stage-wise mapping.
\end{restatable}

%Table~\ref{tab:character} summarizes the general cuts and the additional \sta{} cuts available for each revision mechanism.
%\begin{table}[htbp!]    
%\centering
%    \begin{tabular}{lccc}
 %       \toprule
 %       \textbf{Variants} & \textbf{\shortstack{Monotonicity\\guaranteed}} &
 %       \textbf{\shortstack{General\\cut}} & \textbf{Additional \sta{} cut} \\
 %       \midrule
 %       Budgeted & Yes & MOC & EOC~\eqref{eq:enhanced-cut:bpa-sca} \\
 %       Cost-based & No & IOC & EOC~\eqref{eq:enhanced-cut:cpa-sca} \\
 %       \bottomrule
 %   \end{tabular}
 %   \caption{Optimality cuts for stage-wise and scenario-adaptive variants.}
 %   \label{tab:character}
%\end{table}

\begin{remark}[Relation to the DP] \label{rem:lshaped-dp}
For \sta{}, the DP precomputes all unit costs, whereas the integer
L-shaped method solves station problems~\eqref{eq:station-problem}
during candidate-policy evaluations and can reuse previously computed
values. The latter may therefore solve fewer station problems but
requires repeated revision master problem solves. We compare these approaches in Section~\ref{sec:perf-analysis}.
\end{remark}

\subsection{Reconfiguration Optimization under a Fixed Revision Policy} \label{sec:lexicographic}

We now solve \rec{} under the revision policy and station-count bounds
selected by \sta{}. Reconfiguration costs couple consecutive layouts
through their task assignments. We encode the fixed strategic-decision
linking restrictions on a transformed scenario tree and optimize its
layouts by nested decomposition, strengthened by weighted asymmetric
Lagrangian switch cuts.

Using the node representation $\boldsymbol{\rho}^*$ of the selected
policy $\boldsymbol{\eta}^*$, we merge each stage-$t$ commitment group
$\mathcal{C}_t(u)$, for
$u\in\mathcal{B}(\boldsymbol{\rho}^*)$ and
$t\in[t(u),e_u(\boldsymbol{\rho}^*)]$, into one node of a transformed
tree with node set $\Nrec$.
Lemma~\ref{lem:commitment-structure} ensures that one layout per merged
node enforces exactly the fixed strategic-decision linking restrictions.
We retain the original tree's parent-child connections after merging
each group. For the policy depicted in panel~(b) of
Figure~\ref{fig:scenario-adaptive-tree}, the transformed node for
$\{n_4,n_5\}$ has parent $n_2$ and children $n_8,\ldots,n_{11}$.
Each merged node
inherits the group's common station bound $\zsta{n}$ and weight
$w_n(\boldsymbol{\eta}^*)$, the union of its demand sets, and the sum of
its probabilities.
Throughout this subsection, node indices, parents, children, and
probabilities refer to the transformed tree, with conditional child
probabilities $q_{nn'}:=p_{n'}/p_n$ for $n'\in\mathcal{C}(n)$.

Physical feasibility and the inherited station bounds are independent of
the preceding layout. Thus, only the task-assignment vector is needed as
the state passed between successive nodes. For
$n\in\Nrec$ and an incoming binary vector $\boldsymbol{y}_{a(n)}$, let
$\mathcal{Q}_n(\boldsymbol{y}_{a(n)})$ denote the minimum expected
reconfiguration cost from $n$ onward, conditional on reaching $n$.
The cost-to-go value satisfies
\begin{equation}\label{eq:rec-bellman}
\begin{aligned}
    \mathcal{Q}_n(\boldsymbol{y}_{a(n)})
    =\min_{\substack{\boldsymbol{s}_n\in\mathcal{X}_n\\
                    \sum_{\ell\in\mathcal{L}}x_{n\ell}\le\zsta{n}}}
    \Big\{&w_n(\boldsymbol{\eta}^*)
    \sum_{k\in\mathcal{K}}\sum_{\ell\in\mathcal{L}}
    (1-y_{a(n),k\ell})y_{nk\ell}
    +\sum_{n'\in\mathcal{C}(n)}q_{nn'}\,
    \mathcal{Q}_{n'}(\boldsymbol{y}_n)\Big\}.
\end{aligned}
\end{equation}
Here, $\mathcal{X}_n$ is the feasible strategic set defined in
Section~\ref{sec:structure}, evaluated using the merged demand set.
The first term is the immediate reconfiguration cost, and the second is
the expected child cost-to-go. The selected \sta{} layouts satisfy these
nodal constraints, so the minimization is feasible for every incoming
assignment. The child sum is zero at leaves, and the
\rec{} optimum is $\mathcal{Q}_1(\boldsymbol{0})$ under the greenfield
initial condition.

Nested decomposition~\citep{bayraktar2,zou2019sddip} approximates these
values by replacing the expected child cost-to-go at each non-leaf node
$n$ with $\psi_n$, constrained by lower-bounding cuts and initially by
$\psi_n\ge0$. A forward pass solves the nodal approximations to obtain
candidate assignments, while a backward pass evaluates child problems at these
assignments to generate cuts. The root approximation gives a lower bound.
Evaluating the selected \sta{} layouts on the transformed tree gives an
initial upper bound, which can be improved by subsequent feasible plans. We
repeat these passes until the lower and upper bounds differ by at most a
prescribed tolerance.

We strengthen the nodal approximations by extending the asymmetric
Lagrangian switch cut of~\citep{bayraktar2} to commitment-dependent weights.
Changing one incoming task-station assignment from zero to one can reduce
a child's cost-to-go by at most that child's weight, whereas changing it from one
to zero cannot reduce the cost-to-go. For a non-leaf node $n\in\Nrec$,
the child weights coincide. If the commitment continues, $n$ has one merged
child; otherwise, every child starts a new commitment with lag zero.
We denote this common weight by $w_{\mathcal{C}(n)}$, so
$w_{n'}(\boldsymbol{\eta}^*)=w_{\mathcal{C}(n)}$ for every
$n'\in\mathcal{C}(n)$. In the budgeted variants, $w_{\mathcal{C}(n)}=1$.

\begin{restatable}[Weighted asymmetric Lagrangian switch cut]{proposition}{weightedswitchcut} \label{prop:extended-lsc}
Fix a non-leaf node $n\in\Nrec$ and a binary candidate assignment
$\boldsymbol{y}_n^i$. Let
$v_{n'}^i:=\mathcal{Q}_{n'}(\boldsymbol{y}_n^i)$ for every
$n'\in\mathcal{C}(n)$. Then the cut
\begin{equation}
    \psi_n\ge
    \sum_{n'\in\mathcal{C}(n)}q_{nn'}\,v_{n'}^i
    -w_{\mathcal{C}(n)}\sum_{k\in\mathcal{K}}\sum_{\ell\in\mathcal{L}}
    (1-y_{nk\ell}^i)y_{nk\ell}
    \label{eq:extended-lsc}
\end{equation}
has a right-hand side at most
$\sum_{n'\in\mathcal{C}(n)}q_{nn'}\mathcal{Q}_{n'}(\boldsymbol{y}_n)$
for every binary assignment $\boldsymbol{y}_n$, with equality when
$\boldsymbol{y}_n=\boldsymbol{y}_n^i$.
\end{restatable}

The cut also has a direct Lagrangian interpretation. For a non-leaf node
$n\in\Nrec$ and a child $n'\in\mathcal{C}(n)$, suppose we introduce a copy
$\bar{\boldsymbol{y}}_n$ of the incoming assignment in the child problem
and relax the equality
$\bar{\boldsymbol{y}}_n=\boldsymbol{y}_n^i$ with the term
$\boldsymbol{\pi}^{\top}(\boldsymbol{y}_n^i-\bar{\boldsymbol{y}}_n)$.
With $\pi_{k\ell}=w_{n'}(\boldsymbol{\eta}^*)\big(y_{nk\ell}^i-1\big)$,
the coefficient of $\bar{y}_{nk\ell}$ is
$w_{n'}(\boldsymbol{\eta}^*)\big(1-y_{nk\ell}^i-y_{n'k\ell}\big)$
for every $(k,\ell)\in\mathcal{K}\times\mathcal{L}$. This coefficient is
nonnegative when $y_{nk\ell}^i=0$ and nonpositive when $y_{nk\ell}^i=1$.
Thus, minimizing over
$\bar{\boldsymbol{y}}_n\in[0,1]^{|\mathcal{K}||\mathcal{L}|}$ attains its
minimum at $\bar{\boldsymbol{y}}_n=\boldsymbol{y}_n^i$ for every child layout.
The relaxation attains $v_{n'}^i$, and the multipliers yield the
coefficients in~\eqref{eq:extended-lsc} after probability-weighted aggregation.

In the nested procedure, we replace each exact $v_{n'}^i$ with the optimal
value of the current child approximation. The approximated future cost
remains independent of the incoming assignment, so the same argument
establishes validity for these approximations and equality at the candidate
assignment. Since the approximations lower-bound the true child values,
the cut remains valid for the true expected child cost-to-go, but tightness
for that value requires exact child values.

% ---------------------------------- %
\section{Computational Experiments} \label{sec:experiments}

\subsection{Experimental Design} \label{sec:exp-design}

We evaluate the solution methods and revision policies on matrix assembly layout instances motivated by automotive original equipment manufacturers (OEMs) navigating vehicle electrification. Our standard instances vary along two dimensions: (i) \emph{time granularity}, where we consider monthly (M) and quarterly (Q) planning periods that capture shorter demand cycles relevant to OEM practice, and (ii) \emph{demand heterogeneity}, where we distinguish homogeneous (HM) environments from heterogeneous (HT) ones subject to exogenous shocks that generate localized demand volatility.

The product portfolio comprises four vehicle types, treated as product models: battery electric (BEV), classic hybrid (CH), plug-in hybrid (PH), and combustion engine (CE). The raw task list and precedence relationships for 11 task types are obtained from the benchmark dataset of \citet{scholl1993data}. Following \citet{hottenrott2019flexible}, we generate model-specific task times within $25\%$ of each task's demand-weighted mean. Following \citet{bayraktar}, we normalize the task times, set annual station capacity to 288,000 minutes (25 working days per month and 16 working hours per day), and use compact grid configurations. We divide this capacity by $12$ and $4$ for monthly and quarterly instances, respectively. For each instance, all station problems, i.e., those of the individual nodes and of the commitment groups, use one common compact grid: starting with three candidate workstation locations, we expand the grid until every station problem of the instance is feasible.

Demand data are drawn from vehicle registration records of the German Federal Motor Transport Authority \citep{kbadata}. We remove seasonality with an STL decomposition, fit a vector autoregressive model with a one-period lag, VAR(1), to capture inter-model demand dependencies, and add the seasonal component back to the forecasts; the period length is one month for M instances and one quarter for Q instances. Combining the two granularity levels with the two heterogeneity levels yields four demand datasets, which we refer to as the standard datasets. For HT datasets, each non-root node is shocked independently with probability $0.1$, and a shocked node draws its demand scenarios with forecast errors scaled up by a factor of $5$, producing localized volatility around the VAR(1) forecast; HM datasets use the unscaled VAR(1) errors.

The scenario trees have branching factor $B$ and are generated as follows. Each node $n$ carries a set $\mathcal{D}_n$ of $B$ demand realizations, obtained by adding $B$ forecast errors, sampled from the VAR(1) error distribution by Latin hypercube sampling, to the forecast demand baseline of the node. The baseline of the root node is the one-step VAR(1) forecast from the last observation. At each non-leaf node $n$, every realization $\boldsymbol{d}_n\in\mathcal{D}_n$ leads to one child with conditional probability $1/B$, and the forecast baseline of this child is the one-step VAR(1) forecast from $\boldsymbol{d}_n$.

For each of the four datasets, we generate 5 instances per stage-branching configuration using different random seeds, each corresponding to a distinct scenario tree. We label datasets by granularity and demand type (e.g., M-HM) and stage-branching configurations (e.g., S4-B15). To control tree size, both granularities use S4-B15, S6-B5, and S12-B2; Q also uses S8-B3, and M uses S9-B3. These configurations yield 80 standard scenario trees. The regime-driven datasets introduced in Section~\ref{sec:value-pa} add 40 trees with the same configurations. For the performance analysis in Section~\ref{sec:perf-analysis}, we additionally generate the S5-B2 configuration for the four datasets with five seeds each (20 trees); these trees are not used in the managerial analyses.

In the cost-based variants, we use three forms of the commitment multiplier $h$:
  \begin{enumerate}
      \item Linear: $h(\delta) = \max\{0.01, 1 - c\delta\}$
      \item Convex: $h(\delta) = \max\{0.01, (1-c)^{\delta}\}$
      \item Concave: $h(\delta) = \max\{0.01, 1 - c\delta - c^2\delta(\delta-1)/2\}$
  \end{enumerate}
where $c \in \{0.01, 0.02, \ldots, 0.10\}$ is the discount coefficient per planning period of lag. We use $c=0.03$ as the base case. The $0.01$ floor keeps weights positive; the shape labels refer to the uncapped functions. The three forms and ten coefficients yield 30 combinations per cost-based variant. For budgeted variants, we use $\mu\in\{0,1,\ldots,T-1\}$ and report $\mu\leq3$ in Section~\ref{sec:cpa-vs-bpa}.

Unless stated otherwise, reported values are averages over the scenario trees of a configuration. For cost-based results, we first average the reported measure for each tree over the three discount shapes and, when no coefficient is specified, over all discount coefficients. Shaded bands show $\pm 1$ standard error across trees. Sections~\ref{sec:val-pamsp}--\ref{sec:value-pa} report the 4- to 9-stage trees; the 12-stage trees are reported separately in Figure~\ref{fig:vsa-sensitivity} and used in Section~\ref{sec:granularity}. A model matches a benchmark when their optimal values agree within a relative tolerance of $10^{-6}$; shares of matching cases are computed over all combinations of tree and discount shape.

Experiments ran on an Intel(R) Xeon(R) Platinum 8280L workstation (2.70~GHz, 16~GB RAM) using Python 3.12 and Gurobi 13.0. For each tree, the DP solves the station problems~\eqref{eq:station-problem} once and reuses their optimal station counts across all four \sta{} variants and parameter settings; these problems use Gurobi's default tolerances. The L-shaped \sta{} runs in Section~\ref{sec:perf-analysis} have a 2-hour limit, and their revision master problems use Gurobi's default relative MIP gap of $10^{-4}$. 

Our computational analysis centers on \sta{}, which determines the revision policies and station-use profiles; \rec{} refines the layouts under the selected policy while preserving station use and is not part of the reported experiments. We use \cscensta{} as the primary model. Section~\ref{sec:perf-analysis} evaluates the solution methods, Sections~\ref{sec:val-pamsp}--\ref{sec:value-pa} compare revision policies, and Section~\ref{sec:granularity} studies revision frequency.

\subsection{Performance Analysis} \label{sec:perf-analysis}

We first compare the L-shaped cut families on small trees, where IOC remains tractable, and then compare the DP with the L-shaped method using EOC on the scenario trees of the main study.

\subsubsection{Optimality Cuts of the L-Shaped Method} \label{exp:cuts}

We compare the integer optimality cuts (IOC), the monotone optimality cuts (MOC), and the enhanced optimality cuts (EOC) of Section~\ref{sec:lshaped}.
We use 20 scenario trees with the S5-B2 configuration (31 nodes each).
We fix the revision budget to $\mu = 2$ in \bscensta{} (20 runs per cut family) and use the base case $c = 0.03$ with all three discount shapes in \cscensta{} (60 runs per cut family).
Of the two variants tested here, only \bscensta{} admits MOC, since the required monotonicity need not hold in \cscensta{}.

Table~\ref{tab:performance-cuts} reports the overall search time, time spent solving the revision master problem, and iteration counts.
Search time runs from the start of the L-shaped iterations until optimality is proven; it excludes node-problem initialization, which takes about one second for all cut families.

\begin{table}[htbp!]
      \centering
      \setlength{\tabcolsep}{4pt}
      \scalebox{0.9}{
      \begin{tabular}{llrrrcrcrrr}
           & & \multicolumn{3}{c}{Search time (sec.)} & & Master (sec.) & & \multicolumn{3}{c}{Iterations} \\
           \cmidrule(lr){3-5} \cmidrule(lr){7-7} \cmidrule(lr){9-11}
           & Cut & Avg. & Std. Dev. & Max & & Avg. & & Avg. & Std. Dev. & Max \\ \hline
           \multirow{3}{*}{\bscensta{}} & IOC & 10.03 & 1.49 & 12.08 & & 7.03 & & 355.00 & 0.00 & 355 \\
                                        & MOC & 2.20 & 0.39 & 3.16 & & 0.38 & & 65.20 & 9.35 & 89 \\
                                        & EOC & 1.81 & 0.37 & 2.69 & & 0.05 & & 5.00 & 0.00 & 5 \\ \hline
           \multirow{2}{*}{\cscensta{}} & IOC & 1,157.46 & 145.20 & 1,645.71 & & 1,140.02 & & 4,470.00 & 0.00 & 4,470 \\
                                        & EOC & 1.76 & 0.34 & 2.85 & & 0.05 & & 5.00 & 0.00 & 5 \\ \hline
      \end{tabular}
      }
      \caption{Performance of the L-shaped cut families on S5-B2 trees of the HM and HT datasets.}
      \label{tab:performance-cuts}
\end{table}

In \bscensta{}, IOC evaluates all 355 admissible revision policies in every run before proving optimality.
Exploiting monotonicity through MOC reduces the number of iterations to $65.2$ on average, an $81.6\%$ reduction, and the search time from $10.03$ to $2.20$ seconds.
EOC further exploits the interval decomposition to bound node contributions separately, reducing the iteration count to 5 in every run and the average search time to $1.81$ seconds.

In \cscensta{}, the gain from EOC is larger.
IOC evaluates all 4,470 admissible revision policies in every run, with an average search time of $1{,}157$ seconds, of which $98.5\%$ is spent solving the revision master problem.
EOC requires 5 iterations in every run and averages $1.76$ seconds.
With EOC, time spent solving the revision master problem averages only $0.05$ seconds in both variants.
We therefore use EOC in all subsequent L-shaped runs.

\subsubsection{Dynamic Program vs.\ L-Shaped Method} \label{exp:scalability}

We compare the two exact methods on the 80 standard scenario trees.
For each tree, the DP reuses one set of station-problem solutions to evaluate all $60+2T$ settings: 30 combinations of discount shape and coefficient for each cost-based variant, and all budgets $\mu\in[0,T-1]$ for each budgeted variant.
The L-shaped method solves one setting per run.
We evaluate it for \cscensta{} and \cstgsta{} with three shapes and $c\in\{0.03,0.05,0.10\}$ and for \bscensta{} with $\mu\le2$ (S4-B15, S6-B5) or $\mu\le4$ (otherwise).

Table~\ref{tab-scale} reports DP time \emph{per tree} for all settings together and L-shaped time \emph{per run}, in seconds.
Each row combines the HM and HT instances for the indicated time granularity and stage-branching configuration.
Units counts station problems; Node reports time for individual-node problems.
Lim. and Gap give the number of time-limited runs and their mean optimality gap (\%); Master gives the revision master problem's share of computation time (\%).

  \begin{table}[htbp!]
    \centering
    \footnotesize
    \setlength{\tabcolsep}{4pt}
    \begin{tabular}{ll rr rr c rrrr c rr c r}
      & & & & \multicolumn{2}{c}{DP} & & \multicolumn{9}{c}{L-shaped with EOC} \\
      \cmidrule(lr){5-6} \cmidrule(lr){8-16}
      & & & & & & & \multicolumn{4}{c}{\cscensta{}} & & \multicolumn{2}{c}{\bscensta{}} & & \cstgsta{} \\
      \cmidrule(lr){8-11} \cmidrule(lr){13-14} \cmidrule(lr){16-16}
                    & Config. & Nodes & Units & Node & Total & & Avg. & Lim. & Gap & Master & & Avg. & Lim. & & Avg. \\ \hline
\multirow{4}{*}{Q}  & S4-B15 & 3,616 & 3,874 & 519 & 700 & & 633 & 0 & -- & 0.4 & & 320 & 0 & & 347 \\
                    & S6-B5 & 3,906 & 4,881 & 289 & 581 & & 566 & 0 & -- & 3.1 & & 252 & 0 & & 295 \\
                    & S8-B3 & 3,280 & 4,916 & 186 & 527 & & 1,543 & 7 & 0.83 & 55.7 & & 257 & 0 & & 279 \\
                    & S12-B2 & 4,095 & 8,178 & 186 & 780 & & 3,123 & 26 & 1.30 & 76.0 & & 857 & 1 & & 388 \\ \hline
\multirow{4}{*}{M}  & S4-B15 & 3,616 & 3,874 & 449 & 605 & & 567 & 0 & -- & 0.5 & & 289 & 0 & & 309 \\
                    & S6-B5 & 3,906 & 4,881 & 317 & 588 & & 581 & 0 & -- & 3.9 & & 258 & 0 & & 299 \\
                    & S9-B3 & 9,841 & 14,757 & 685 & 1,755 & & 5,059 & 48 & 1.71 & 59.4 & & 1,590 & 3 & & 912 \\
                    & S12-B2 & 4,095 & 8,178 & 199 & 748 & & 3,543 & 34 & 1.15 & 80.6 & & 1,937 & 9 & & 354 \\ \hline
    \end{tabular}
    \caption{Computational performance of the DP and L-shaped method on the HM and HT datasets. L-shaped averages are over 90 runs per row for \cscensta{} and \cstgsta{} (three shapes, $c\in\{0.03,0.05,0.10\}$, ten trees) and over 30 or 50 runs for \bscensta{}.}
    \label{tab-scale}
  \end{table}

The DP takes between 527 and 1,755 seconds per tree on average.
The maximum over all 120 trees, including the 40 regime-driven trees introduced in Section~\ref{sec:value-pa}, is 2,408 seconds.
Station-problem solves dominate computation: individual-node problems account for $24\%$-$74\%$ of the average time, whereas the recursions take at most 0.36 seconds per setting.
The interval decomposition thus makes comparisons across variants, budgets, and commitment multipliers inexpensive once the station problems have been solved.
Precomputation can still be demanding, since the largest station problem pools the demand sets of all leaf nodes.

The L-shaped method solves all \cscensta{} runs on four- and six-stage trees to optimality. An average run takes about as long as the DP for all settings; the revision master problem accounts for at most $3.9\%$ of the time.
With at least eight stages, an average run takes $2.9$ to $4.7$ times as long as the DP for all settings, and $7$ to $48$ of the $90$ runs per row reach the time limit.
With EOC, each run takes at most 18 iterations, but solving the revision master problem consumes $56\%$-$81\%$ of the time on these deeper trees.
For the stage-wise counterpart \cstgsta{}, the L-shaped method solves all runs to optimality and spends less than $3.2\%$ of its time in the revision master problem.
Its revision master problem has $T-1$ revision variables, whereas that of \cscensta{} has one per non-leaf node.
The L-shaped method reaches the time limit in 13 \bscensta{} runs.

For the time-limited L-shaped runs, the reported optimality gaps average $1.27\%$ (at most $4.80\%$), and most of their computation time is spent in the revision master problem. The main remaining difficulty is therefore improving the revision master problem's lower bound to certify solution quality.
The DP offers a clear computational advantage for the full \sta{} parameter study.
The integer L-shaped framework nevertheless remains applicable when interstage coupling prevents the interval decomposition required by the DP.

\subsection{Value of Partial Adaptivity} \label{sec:val-pamsp}

\begin{figure}[t]
    \centering
    \includegraphics[width=0.77\linewidth]{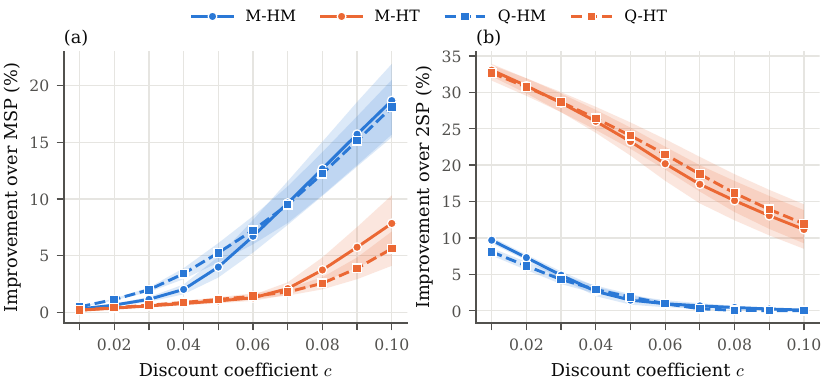}
    \caption{Improvement (\%) of \cscensta{} over (a) \msp{} and (b) \twosp{}, by discount coefficient $c$.}
    \label{fig:val-pamsp}
  \end{figure}

We assess the value of partial adaptivity by comparing \cscensta{} with the benchmarks in Section~\ref{sec:budgeted-stagewise-pa}: \msp{}, which retains full adaptivity without commitment discounts, and \twosp{}, which commits the entire layout plan at the first stage and receives the largest discounts.
Since both benchmark policies are feasible in \cscensta{}, its optimal cost cannot exceed either benchmark's cost.
Figure~\ref{fig:val-pamsp} reports percentage reductions in expected discounted station-use cost for the 60 HM and HT scenario trees with four to nine stages.
Ranges below span the monthly and quarterly datasets.

At the base case $c = 0.03$, partial adaptivity yields its largest gains over full commitment under heterogeneous demand.
The improvement over \twosp{} is $4.29$-$4.88\%$ under homogeneous demand and $28.55$-$28.65\%$ under heterogeneous demand, with at least $18.7\%$ in every heterogeneous instance.
Improvements over \msp{} are smaller, at $1.14$-$1.98\%$ and $0.56$-$0.62\%$, respectively.
The corresponding objective values are reported in Table~\ref{tab:obj-levels} in Appendix~\ref{app:obj-levels}.
These results show that selective commitment can yield a net benefit from discounts while retaining enough adaptability to reduce cost relative to full commitment, especially under heterogeneous demand.

As $c$ grows, stronger commitment incentives increase the gains over \msp{} and shift the selected revision policies toward full commitment.
At $c = 0.10$, the improvement over \msp{} reaches $18.09$-$18.67\%$ under homogeneous demand and $5.61$-$7.84\%$ under heterogeneous demand.
Full commitment becomes nearly optimal on average under homogeneous demand: the gain over \twosp{} falls below $0.5\%$ from $c = 0.08$ onward.
Under heterogeneous demand, the gain remains $11.16$-$11.94\%$ even at $c = 0.10$, indicating that selective revision remains valuable throughout the tested discount range.

\begin{figure}[htbp!]
    \centering
    \includegraphics[width=0.77\linewidth]{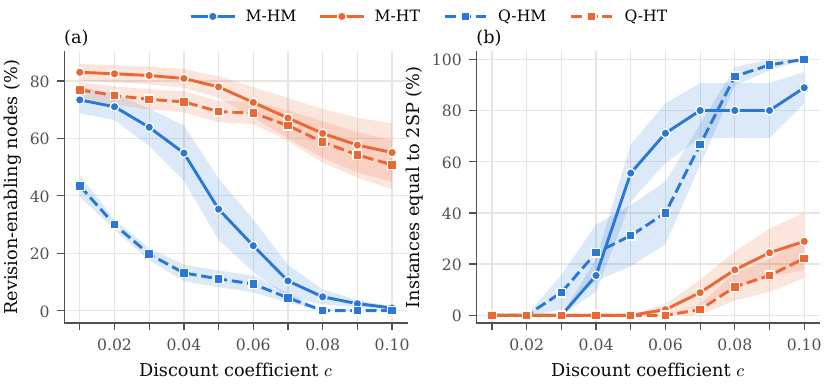}
    \caption{(a) Revision-enabling nodes (\% of non-leaf nodes) in selected \cscensta{} policies and (b) agreement with \twosp{} in optimal cost, by discount coefficient $c$.}
    \label{fig:val-pamsp-rev}
\end{figure}

Figure~\ref{fig:val-pamsp-rev} shows that revision opportunities decline as $c$ grows, but more sharply under homogeneous demand.
Panel~(a) reports the share of non-leaf nodes that enable a revision, i.e., the share of non-root nodes with a revision opportunity, and panel~(b) reports the frequency of matching \twosp{} in optimal cost. Under homogeneous demand, the revision rate falls from $43.47$-$73.43\%$ at $c = 0.01$ to $19.68$-$63.89\%$ at the base case and to less than $1\%$ at $c = 0.10$. Under heterogeneous demand, it falls only from $76.93$-$83.05\%$ at $c = 0.01$ to $50.84$-$55.08\%$ at $c = 0.10$. At $c = 0.10$, the frequency of matching \twosp{} is $89$-$100\%$ under homogeneous demand and at most $29\%$ under heterogeneous demand. Weighted by node probability, the base-case policies enable a revision at $38.2$-$60.8\%$ (homogeneous) and $76.5$-$78.3\%$ (heterogeneous) of the stages along a scenario path. The model thus commits selectively under homogeneous demand and retains most revision opportunities under heterogeneous demand, while still collecting part of the commitment discount.

\begin{figure}[htbp!]
    \centering
    \includegraphics[width=0.77\linewidth]{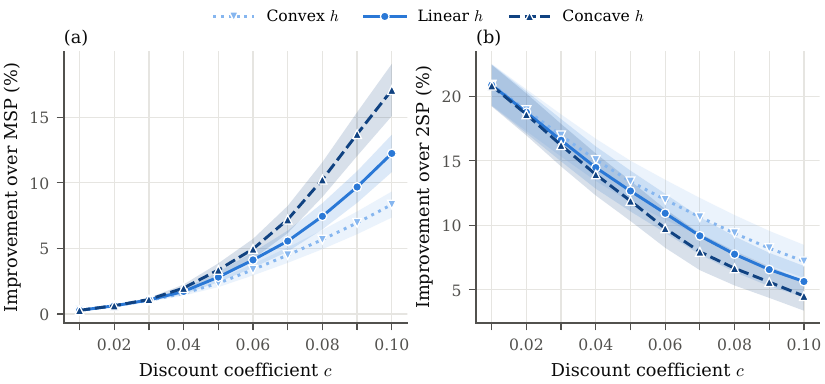}
    \caption{Improvement (\%) of \cscensta{} over (a) \msp{} and (b) \twosp{}, by discount coefficient $c$ and shape of $h$.}
    \label{fig:val-pamsp-shape}
  \end{figure}

Finally, we examine how the shape of $h$ affects these gains, pooling the four datasets (Figure~\ref{fig:val-pamsp-shape}).
Averaged over the coefficients, the improvement over \msp{} is largest under the concave shape and smallest under the convex shape, while the improvement over \twosp{} follows the reverse order.
This pattern reflects the stronger discounts offered by the concave specification at long lags.
Differences are small at low coefficients: for $c \leq 0.03$, mean improvements differ by at most $0.07$ percentage points over \msp{} and $0.75$ over \twosp{}.
At $c = 0.10$, however, the improvement over \msp{} is $8.35\%$ under the convex, $12.24\%$ under the linear, and $17.07\%$ under the concave shape.

The results suggest a practical guideline for planners.
Under homogeneous demand and strong commitment incentives, full commitment can be an effective simplification: it is nearly optimal on average in our experiments.
Under heterogeneous demand, selective revision remains valuable throughout the tested discount range.
By jointly optimizing layout plans and revision timing, \cscensta{} determines how much adaptability to retain without requiring the planner to choose between full adaptivity and full commitment in advance.

\subsection{Effect of the Revision Budget} \label{sec:cpa-vs-bpa}

\begin{figure}[htbp!]
  \centering
  \includegraphics[width=0.77\linewidth]{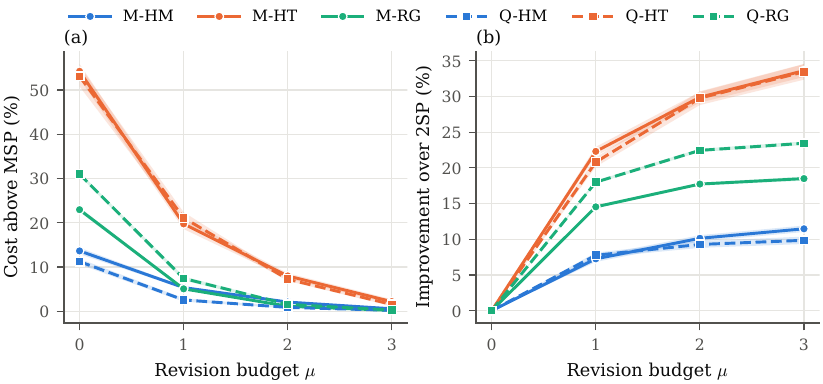}
  \caption{(a) Cost (\%) of \bscensta{} above \msp{} and (b) improvement (\%) of \bscensta{} over \twosp{}, by revision budget $\mu$.}
  \label{fig:bpa}
\end{figure}

Having examined the benefits of commitment discounts, we now evaluate the cost of limiting revision opportunities through an exogenous budget.
Figure~\ref{fig:bpa} compares \bscensta{} with \msp{} and \twosp{} as $\mu$ grows from $0$ to $3$, on the four- to nine-stage trees.
Both benchmarks use undiscounted station-use costs here: \twosp{} is the zero-budget case, and \msp{} provides a lower bound.
In addition to HM and HT, the figure includes regime-driven demand (RG, introduced in Section~\ref{sec:value-pa}).

A single revision opportunity per scenario path closes most of the cost gap between full commitment and full adaptivity in the dataset averages: $61$-$77\%$ under homogeneous demand, $60$-$64\%$ under heterogeneous demand, and $76$-$78\%$ under regime-driven demand, across the monthly and quarterly datasets.
Additional opportunities yield smaller marginal cost reductions in the reported averages.
With $\mu = 3$, \bscensta{} is within $0.23$-$0.58\%$ (HM) and $0.21$-$0.28\%$ (RG) of \msp{}, whereas heterogeneous demand still leaves a gap of $1.55$-$2.20\%$ (HT).
The remaining gaps arise from the deeper trees, since $\mu=3$ permits full adaptivity on the four-stage trees.

For planners, these results suggest that a few well-chosen revision opportunities can retain much of the benefit of full adaptivity, although heterogeneous demand incurs a larger residual cost at the same budget.

\subsection{Value of Scenario Adaptivity} \label{sec:value-pa}
We next compare the cost-based variants \cscensta{} and \cstgsta{} to assess the value of allowing revision timing to vary across scenario paths.
Proposition~\ref{prop:scamsp-generalizes-pamsp} guarantees that scenario adaptivity cannot increase the optimal cost; here, we investigate when its benefit is substantial.
We first examine the results on the standard datasets, then analyze when stage-wise commitment can fall short.
This analysis guides an experiment on a regime-driven (RG) dataset motivated by the electrification transition, with persistent shifts in product demand.

Table~\ref{tab:vsca} reports averages over 15 trees per row. On the standard datasets, \cscensta{} improves on \cstgsta{} by $0.40$-$0.71\%$ at the base case.
In what follows, we examine commitment incentives to help explain why these gains are modest and when they can be larger.
Standard errors across trees are at most $1.38$ percentage points for all reported improvements and $0.23$ percentage points for the improvement over \cstgsta{}.

\begin{table}[htbp!]
  \centering
  \scalebox{0.9}{
  \begin{tabular}{l|cccc}
           & \msp{}   & \twosp{}  & \cstgsta{} & DI (\%)    \\ \hline
      M-HM & 1.14  & 4.88  & 0.71      & 0.65  \\
      M-HT & 0.56  & 28.55 & 0.40      & 0.46  \\
      M-RG & 4.99  & 15.31 & 2.45      & 0.98  \\
      Q-HM & 1.98  & 4.29  & 0.65      & 0.78  \\
      Q-HT & 0.62  & 28.65 & 0.59      & 0.70  \\
      Q-RG & 4.43  & 20.45 & 2.58      & 0.99
  \end{tabular}
  }
  \caption{\cscensta{} improvement (\%) over alternative policies, and dispersion index, at $c=0.03$.}
  \label{tab:vsca}
\end{table}

\subsubsection{Commitment Incentives}

For fixed $t\in[2,T]$, we compare policies that are fully adaptive except possibly at nodes in $\mathcal{N}_{t-1}$. For $u\in\mathcal{N}_{t-1}$, we define the \emph{commitment gain} from closing only its revision window as the cost reduction relative to full adaptivity, using the station-problem values~\eqref{eq:station-problem}:
\[
\gamma_u := \sum_{n\in\mathcal{C}(u)}p_n\,\mathrm{SP}^*(\{n\})
-h(1)p_u\,\mathrm{SP}^*(\mathcal{C}(u)).
\]
The first term permits separate layouts at the children, while the second requires a shared layout with a lag-one discount. The sign of $\gamma_u$ indicates whether the discount savings exceed or fall short of the added station-use cost of sharing a layout.

\begin{restatable}[Single-stage commitment gains]{proposition}{commitmentgain} \label{prop:commitment-gain}
Suppose \msp{} is feasible and the children of every non-leaf node admit a common feasible layout. For any $t\in[2,T]$, consider revision policies that are fully adaptive except possibly at nodes in $\mathcal{N}_{t-1}$. Relative to \msp{}, the largest cost reductions in this class are $\big[\sum_{u\in\mathcal{N}_{t-1}}\gamma_u\big]^+$ under stage-wise timing and $\sum_{u\in\mathcal{N}_{t-1}}[\gamma_u]^+$ under scenario-adaptive timing.
If the optimal \cstgsta{} cost equals the \msp{} cost, then $\sum_{u\in\mathcal{N}_{t-1}}\gamma_u\leq0$ for every $t\in[2,T]$. In that case, \cscensta{} improves on \cstgsta{} by at least $\max_{t\in[2,T]}\sum_{u\in\mathcal{N}_{t-1}}[\gamma_u]^+$.
\end{restatable}

At a fixed stage, \cscensta{} can close just the windows with $\gamma_u>0$, giving the sum of positive gains in Proposition~\ref{prop:commitment-gain}. By contrast, \cstgsta{} must close every window at that stage or none, giving the positive part of their total gain. We define the \emph{dispersion index} by summing these differences across stages and expressing the result as a percentage of the \msp{} cost:
\[
\mathrm{DI}:=\frac{100}{z_{\msp{}}}
\sum_{t=2}^{T}\Big(\sum_{u\in\mathcal{N}_{t-1}}[\gamma_u]^+
-\big[\sum_{u\in\mathcal{N}_{t-1}}\gamma_u\big]^+\Big).
\]
DI need not equal or bound the optimal cost difference between \cstgsta{} and \cscensta{}, since commitments at different stages may interact. It is independent of discount shape because $h(1)=1-c$ for all three shapes. A shared layout cannot use fewer stations than each child separately, giving $0\le\mathrm{DI}\le100c$. This bound may be loose, but its $3\%$ value in our base case provides a scale for interpreting the DI values.
For the standard datasets, Table~\ref{tab:vsca} reports average base-case DI values of $0.46$-$0.78\%$, amounting to up to roughly one-quarter of the $3\%$ upper bound.

\subsubsection{Regime-Driven Demand}

Proposition~\ref{prop:commitment-gain} shows that scenario-adaptive timing has a single-stage advantage when commitment pays off at some nodes but not others.
The standard datasets use the same VAR(1) dynamics throughout the tree, with localized forecast-error shocks in HT.
These generators do not explicitly represent persistent branch-specific demand regimes, although such regimes may arise in practice.
To examine how demand regimes affect commitment incentives, we construct the RG dataset with persistent shifts toward electric or combustion vehicles and small variation within each regime.

The root starts in the baseline regime, and a node is \emph{stable} if its path has not yet entered an EV- or CE-dominant regime.
Each child of a stable node has conditional probability $1/B$. For $B\ge3$, $\lfloor B/3\rfloor$ children remain stable, $\lfloor (B-\lfloor B/3\rfloor)/2\rfloor$ enter an EV-dominant regime, and the remaining children enter a CE-dominant regime, i.e., 5/5/5 children for $B=15$, 1/2/2 for $B=5$, and 1/1/1 for $B=3$.
In the EV-dominant regime, BEV and PH demand is scaled up and CH and CE demand down by a factor of $3$; the CE-dominant regime reverses this scaling.
Regimes are absorbing, i.e., all descendants of a regime node remain in its regime.
Each node carries a latent state that follows the VAR(1) model with forecast errors scaled to $5\%$ of their estimated size, independently of the regime, and its demand is this state scaled by the regime factors. Unlike in the standard datasets, children branch from the parent's latent state rather than from its demand realizations. Demand scenarios within a node are identical on regime paths and tightly concentrated on stable paths.
For S12-B2, each stable node has one stable and one CE-dominant child.

Within a persistent regime, similar demands at the children of a regime node can limit the extra station use required by a shared layout, favoring commitment.
Children of a stable node, however, may require different layouts as they enter different regimes.
Scenario adaptivity can then retain revision opportunities while a regime change remains possible and allow commitment after the regime is established.

The RG rows of Table~\ref{tab:vsca} show that, at the base case, the improvement of \cscensta{} over \cstgsta{} rises to $2.45\%$ for M-RG and $2.58\%$ for Q-RG.
\cscensta{} also improves on \msp{} by $4.43$-$4.99\%$ and on \twosp{} by $15.31$-$20.45\%$.
Even \cstgsta{} improves on \msp{} by $1.90$-$2.63\%$, so a common revision calendar remains valuable.
DI also rises to $0.98$-$0.99\%$, roughly one-third of the $3\%$ upper bound.

\begin{figure}[!htb]
    \centering
    \includegraphics[width=0.77\linewidth]{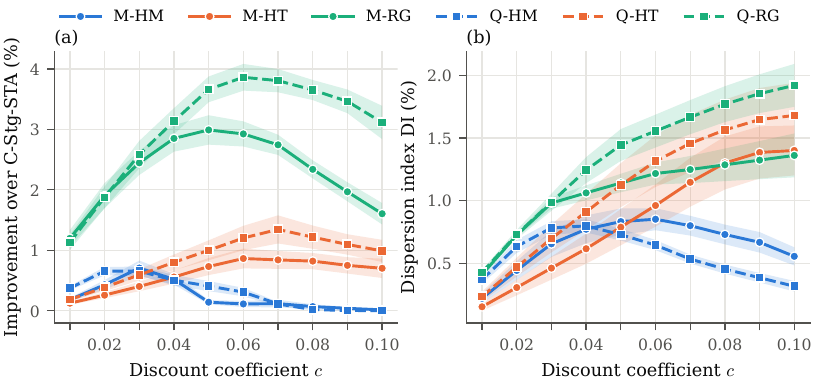}
    \caption{(a) Improvement (\%) of \cscensta{} over \cstgsta{} and (b) DI (\%), by discount coefficient $c$.}
    \label{fig:vsca-c}
  \end{figure}

Figure~\ref{fig:vsca-c}(a) shows that the improvement over \cstgsta{} first rises and then falls as $c$ increases.
Under homogeneous demand, it is greatest near the base case and nearly disappears at strong discounts; under heterogeneous demand, it peaks later and is larger on average across the tested coefficients, reversing the base-case ordering.
On RG, the improvement remains above $1\%$ throughout and reaches $2.99$-$3.86\%$ at intermediate discounts.
Panel~(b) illustrates the local scope of DI, which continues to rise on HT and RG after the realized improvement begins to decline.
DI uses only the lag-one discount around full adaptivity, whereas the optimized policies can exploit longer commitment intervals.

\begin{figure}[!htb]
    \centering
    \includegraphics[width=0.77\linewidth]{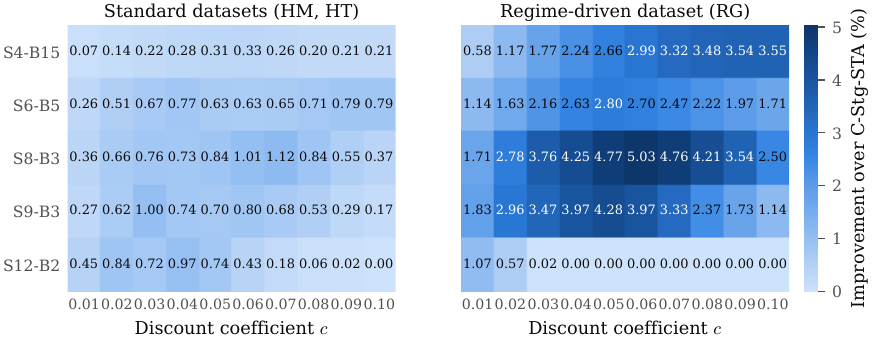}
    \caption{Improvement (\%) of \cscensta{} over \cstgsta{}, by tree configuration and discount coefficient $c$.}
    \label{fig:vsa-sensitivity}
\end{figure}

Figure~\ref{fig:vsa-sensitivity} separates the results by tree configuration, including the 12-stage trees, with HM and HT pooled.
In these configuration averages, the improvement never exceeds $1.12\%$ on the standard datasets. On RG, it reaches $5.03\%$ on S8-B3 and remains at least $1.14\%$ at $c=0.10$ across all four- to nine-stage configurations.
The six- to nine-stage RG gains rise and then fall with $c$, whereas the S4-B15 gain continues to increase through the largest tested coefficient.

On the 12-stage RG trees, by contrast, the improvement is nearly zero from $c=0.03$; we observe that both models match the \twosp{} cost in $87\%$ of cases at that coefficient and in all cases from $c=0.05$.
Commitment lags can reach 11 periods, while the probability of remaining on the only stable path halves each stage, favoring full commitment.
On the standard 12-stage trees, both models match the \twosp{} cost in all cases at $c=0.10$.

The results suggest a practical decision rule. A common revision calendar may suffice when demand paths call for similar commitment timing or when discounts make full commitment optimal. Scenario-adaptive timing is worth considering when persistent demand shifts make commitment attractive on some paths while revision opportunities remain valuable on others, a trade-off that can be especially consequential at intermediate discounts. For automotive planners facing electrification-driven demand shifts, the key question is whether different market paths settle into regimes at different times; where they do, path-specific revision timing can capture value that a common calendar leaves unrealized.

\subsection{Revision Frequency: Monthly vs. Quarterly vs. Yearly} \label{sec:granularity}

Matrix layouts can be revised more often than conventional lines, but each revision requires planning effort, and a firm may prefer a coarser revision calendar.
We therefore assess the value of monthly revision opportunities in \cscensta{} relative to quarterly and yearly calendars.

To isolate revision frequency, we compare three calendars on the same M-S12-B2 scenario trees of the M-HM, M-HT, and M-RG datasets, keeping demand, station capacity, and monthly decision stages unchanged.
The monthly calendar allows revision opportunities at every stage after the first, as in \cscensta{}.
The quarterly calendar allows them only at the start of months $4$, $7$, and $10$, with $\rho_n=0$ for every non-leaf node $n$ whose children do not start a quarter.
The yearly calendar allows no revision within the horizon and thus coincides with \twosp{}.
Revision opportunities retain the one-stage-ahead authorization rule.
These calendars restrict adaptive layout revisions; committed plans may still specify different layouts each month, and operational decisions remain adaptive.

For all three calendars, we measure commitment lag in quarters by setting $w_n(\boldsymbol{\rho})=h(\delta_n(\boldsymbol{\rho})/3)$.
Thus, $c$ applies per quarter of lag here, rather than per planning period as in the preceding experiments.
The feasible policy sets are nested, so the monthly optimal cost cannot exceed the quarterly cost, which cannot exceed the yearly cost.
The DP extensions in Remark~\ref{rem:dp-extensions} handle calendar restrictions and time scaling, reusing precomputed station counts without further station-problem solves.
Table~\ref{tab:granularity} reports the base case $c=0.03$, averaged over five trees per dataset and three discount shapes, with improvements relative to the more restricted calendar's cost.

\begin{table}[t]
      \centering
      \scalebox{0.9}{
      \begin{tabular}{l ccc cc}
       & \multicolumn{3}{c}{Objective} & \multicolumn{2}{c}{Improvement (\%)} \\
      \cmidrule(lr){2-4}\cmidrule(lr){5-6}
      Dataset & Monthly & Quarterly & Yearly (\twosp{}) & M vs.\ Q & Q vs.\ Y \\ \hline
      M-HM & 57.42 & 59.16  & 62.17  & 2.94 & 4.81  \\
      M-HT & 58.97 & 62.48  & 81.53  & 5.61 & 23.35 \\
      M-RG & 98.10 & 100.64 & 100.64 & 2.52 & 0.00  \\ \hline
      \end{tabular}
      }
      \caption{Monthly, quarterly, and yearly revision calendars on the same monthly trees at $c=0.03$.}
      \label{tab:granularity}
\end{table}

Monthly revision opportunities reduce discounted station-use cost relative to quarterly opportunities by $2.94\%$ under homogeneous demand and $5.61\%$ under heterogeneous demand, on average. Across trees, the corresponding gains are $2.43$-$3.59\%$ and $4.66$-$6.53\%$, respectively.
Quarterly opportunities in turn improve on \twosp{} by $4.81\%$ and $23.35\%$, respectively.
The larger monthly gain under heterogeneous demand reflects the value of adapting layout plans between quarter starts, beyond the operational flexibility retained under either calendar.
At the base case, the selected quarterly policies authorize revision windows at $88\%$ (HM) and $99\%$ (HT) of eligible nodes, while the monthly cost is within $0.6\%$ of \msp{}.

The regime-driven dataset shows a different pattern: the quarterly variant matches the discounted \twosp{} cost in all instances, while monthly opportunities reduce it by $2.52\%$.
In these binary RG trees, half of the probability enters the regime in the second month, and the stable path has probability $2^{-(t-1)}$ in month $t$.
The selected monthly policies authorize revision at every non-leaf stable node, so that new commitment intervals begin when the regime is revealed, and mostly remain committed after regime entry; they authorize further revision windows at only $9$-$24\%$ of the non-leaf regime nodes, mostly in month $11$.
The first quarterly opportunity is month $4$, when the probability of remaining stable is only $1/8$.
Because no earlier window is permitted, all month-3 nodes make the same choice about opening the month-4 window. Opening it resets commitment lags on all branches, including those already in the persistent regime.
At the base case, the resulting layout savings do not outweigh the lost discounts.
The value of a finer calendar therefore depends on whether its additional opportunities align with regime entry, not just on their frequency.

\begin{figure}[htbp!]
\centering
\includegraphics[width=0.77\linewidth]{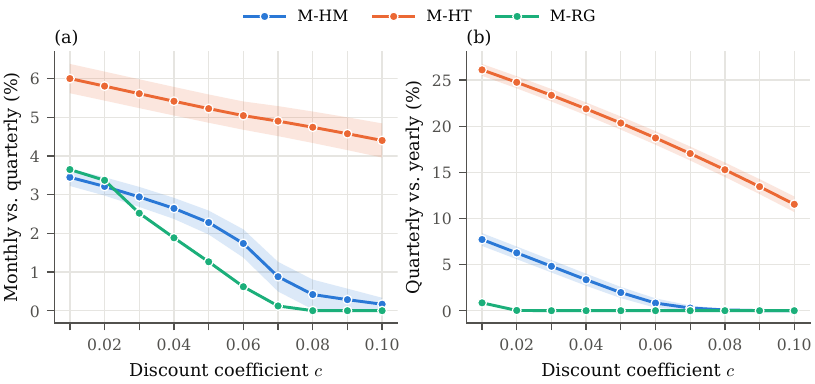}
\caption{(a) Improvement (\%) of monthly over quarterly revision calendars and (b) of quarterly over yearly calendars, by discount coefficient $c$.}
\label{fig:granularity-c}
\end{figure}

Figure~\ref{fig:granularity-c} shows how the gains from finer revision calendars vary with the discount coefficient.
Under homogeneous demand, both gains decline as discounts strengthen; at $c=0.10$, quarterly costs match \twosp{} in all instances.
Under heterogeneous demand, both gains remain substantial throughout the tested range.
On RG, quarterly opportunities offer little improvement over yearly commitment, while monthly opportunities help at low discounts. From $c=0.08$, full commitment is optimal under all three calendars.

\begin{figure}[htbp!]
  \centering
  \includegraphics[width=0.77\linewidth]{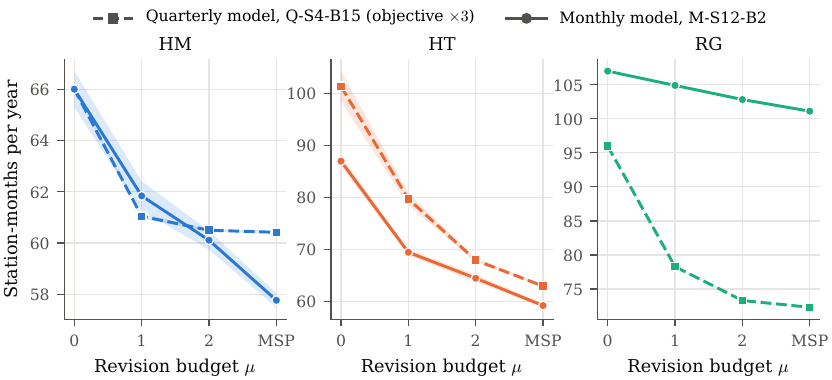}
  \caption{Expected station-months of \bscensta{} on quarterly and monthly instances, by revision budget $\mu$, with \msp{} as reference.}
  \label{fig:granularity-exp-bpa}
\end{figure}

As a complementary comparison, Figure~\ref{fig:granularity-exp-bpa} reports \bscensta{} results on separately generated Q-S4-B15 and M-S12-B2 instances at budgets $\mu\in\{0,1,2\}$. We multiply quarterly objectives by $3$ to express both in station-months. Even with full adaptivity, the monthly \msp{} cost is $4.4\%$ lower for HM and $6.0\%$ lower for HT, but $39.8\%$ higher for RG, relative to the quarterly instances. At the tested budgets, monthly costs are nearly equal to quarterly costs under HM, lower under HT, and higher under RG. The cost differences here also reflect changes in demand models, scenario trees, and period lengths; the same-tree comparison above isolates revision frequency.

For planners, these results suggest matching the revision calendar to demand dynamics.
Quarterly opportunities can be sufficient under homogeneous demand and strong commitment discounts, whereas monthly opportunities retain value under heterogeneous demand throughout the tested range.
The RG results further show that a calendar can miss valuable opportunities near regime entry even when it permits several revisions later.
The measured savings provide a basis for weighing revision frequency against planning effort.

% ---------------------------------- %
\section{Conclusion}\label{sec:conclusions}

Matrix assembly layouts can be revised more frequently than the annual schedules represented in existing models. The planning challenge is to use this flexibility selectively: committing to layouts in advance can support implementation and earn discounts, while retaining revision opportunities allows the plan to respond to evolving demand. We formalized this choice as the layout revision timing problem. Our cost-based and scenario-adaptive extensions to partially adaptive multistage stochastic programming jointly determine the revision policy and committed layout plans. For matrix layouts, a hierarchical procedure subsequently minimizes reconfiguration while preserving the optimal station-use level.

The structure induced by a revision policy also makes the problem more tractable. For station use, commitment costs decompose, yielding an exact dynamic program whose station subproblem solutions can be reused across all four partially adaptive models and parameter settings. The integer L-shaped method applies more generally, including when decisions are coupled across stages; its monotonicity-based cuts strengthen budgeted partially adaptive models beyond the layout application. For reconfiguration, we develop commitment-dependent cuts within nested decomposition.

The experiments show that selective commitment can reduce expected station-use cost relative to both full adaptivity and full commitment. In the budgeted experiments, a single revision opportunity per scenario path recovers most of the cost advantage of full adaptivity in the dataset averages. A common revision calendar performs well when demand paths call for similar commitment timing. Scenario-adaptive timing becomes more valuable when persistent shifts make commitment attractive on some paths while revisions remain valuable on others, as in the electrification-driven instances. The monthly-versus-quarterly comparison further shows that additional revision opportunities matter when they occur near demand changes that call for different layouts.

Future work could estimate commitment benefits and persistent demand transitions from industrial data, testing whether these policy patterns extend beyond the stylized discounts and constructed regimes considered here. Applications in other multistage planning settings would test the broader framework and motivate stronger bounds for the integer L-shaped method when interstage coupling prevents the interval decomposition underlying our dynamic program.

%%%%%%%%%%%%%%%%%%%%%%%%%%%%%%%%%%%%%%%%%%%%%%%%
%%%%%%%%%%%%%%%%%%%%%%%%%%%%%%%%%%%%%%%%%%%%%%%%

% ---------------------------------- %
%\bibliographystyle{plain}
\bibliographystyle{apalike}
\bibliography{References}
% ---------------------------------- %

\newpage
\appendix
\section{Detailed Matrix Assembly Layout Formulations}\label{app:formulation}

We expand the physical feasibility sets and give the extensive \sta{} and \rec{} formulations, using the notation of Section~\ref{sec:objfunc-mald}.

\emph{Additional notation.} Let $\mathcal{M}_k:=\{m\in\mathcal{M}:k\in\mathcal{K}_m\}$ be the set of product models requiring task $k$, and let $\mathcal{K}^{\mathrm{succ}}_{mk}$ be the set of successors of task $k$ for model $m$. The processing time of task $k$ for model $m$ is $R_{mk}$, and each station has $C$ units of processing time available per period. Let $E^{\mathrm{col}}_\ell$ denote the given matrix column of station $\ell$, and let $\mathcal{E}$ be the set of columns. We define $\mathcal{L}_e:=\{\ell\in\mathcal{L}:E^{\mathrm{col}}_\ell=e\}$ for each $e\in\mathcal{E}$ and $\mathcal{L}^{<}_\ell:=\{\ell'\in\mathcal{L}:E^{\mathrm{col}}_{\ell'}<E^{\mathrm{col}}_\ell\}$ for each station $\ell$.

\subsection{Station problem}

The fully adaptive station formulation, before adding partial adaptivity constraints, is
\begin{subequations}\label{app:sta}
  \small
  \begin{align}
      (\text{\sta{}}): \ \min_{\boldsymbol{s},\boldsymbol{o}}\ & \sum_{n \in \mathcal{N}} p_n \sum_{\ell \in \mathcal{L}} x_{n\ell} && \label{app:sta-obj} \\
      \text{s.t.} \ & \sum_{\ell \in \mathcal{L}_e} x_{n\ell} \geq 1 && \forall e \in \mathcal{E},\ n \in \mathcal{N}, \label{app:sta-col} \\
      & x_{n\ell} \geq y_{nk\ell} && \forall k \in \mathcal{K},\ \ell \in \mathcal{L},\ n \in \mathcal{N}, \label{app:sta-xy} \\
      & |\mathcal{L}^{<}_\ell| \cdot (1 - y_{nk\ell}) \geq \sum_{\ell' \in \mathcal{L}^{<}_\ell} y_{nk'\ell'} && \forall m \in \mathcal{M},\ k \in \mathcal{K}_m,\ k' \in \mathcal{K}^{\mathrm{succ}}_{mk},\ \ell \in \mathcal{L},\ n \in \mathcal{N}, \label{app:sta-prec} \\
      & \sum_{m \in \mathcal{M}_k} o_{nmk\ell}(\boldsymbol{d}_n) \leq y_{nk\ell} \cdot \hspace*{-0.15cm} \sum_{m \in \mathcal{M}_k} d_{nm} && \forall k \in \mathcal{K},\ \ell \in \mathcal{L},\ n \in \mathcal{N},\ \boldsymbol{d}_n \in \mathcal{D}_n, \label{app:sta-oy} \\
      & \sum_{\ell \in \mathcal{L}} o_{nmk\ell}(\boldsymbol{d}_n) = d_{nm} && \forall m \in \mathcal{M},\ k \in \mathcal{K}_m,\ n \in \mathcal{N},\ \boldsymbol{d}_n \in \mathcal{D}_n, \label{app:sta-dem} \\
      & \sum_{m \in \mathcal{M}} \sum_{k \in \mathcal{K}_m} R_{mk} \cdot o_{nmk\ell}(\boldsymbol{d}_n) \leq C && \forall \ell \in \mathcal{L},\ n \in \mathcal{N},\ \boldsymbol{d}_n \in \mathcal{D}_n, \label{app:sta-cap} \\
      & x_{n\ell},\ y_{nk\ell} \in \{0,1\} && \forall k \in \mathcal{K},\ \ell \in \mathcal{L},\ n \in \mathcal{N}, \label{app:sta-bin} \\
      & o_{nmk\ell}(\boldsymbol{d}_n) \geq 0 && \forall m \in \mathcal{M},\ k \in \mathcal{K}_m,\ \ell \in \mathcal{L},\ n \in \mathcal{N},\ \boldsymbol{d}_n \in \mathcal{D}_n. \label{app:sta-nonneg}
  \end{align}
  \end{subequations}

The physical feasibility constraints~\eqref{app:sta-col}-\eqref{app:sta-cap} are adapted from~\citet{bayraktar}. Constraints~\eqref{app:sta-col} require at least one active station in each matrix column, and constraints~\eqref{app:sta-xy} permit task assignments only to active stations. In constraints~\eqref{app:sta-prec}, assigning task $k$ to station $\ell$ forces all assignments of its successors to earlier columns to zero. When $y_{nk\ell}=0$, the coefficient $|\mathcal{L}^{<}_\ell|$ makes the constraint nonbinding because it bounds the number of stations counted on the right-hand side. Constraints~\eqref{app:sta-oy} permit production of a task at a station only if it is assigned there. Constraints~\eqref{app:sta-dem} require the demand of each model to be processed for every required task, and constraints~\eqref{app:sta-cap} limit station workloads. Constraints~\eqref{app:sta-col}-\eqref{app:sta-prec} and \eqref{app:sta-bin} define $\mathcal{S}_n$, while constraints~\eqref{app:sta-oy}-\eqref{app:sta-cap} and \eqref{app:sta-nonneg} define $\mathcal{O}_n(\boldsymbol{s}_n,\boldsymbol{d}_n)$.

Each partially adaptive variant also optimizes the corresponding revision decisions $\boldsymbol{\eta}$ and adds the partial adaptivity constraints listed in Table~\ref{tab:app-sta-variants}. With $\boldsymbol{s}_n=(\boldsymbol{x}_n,\boldsymbol{y}_n)$, the strategic-decision linking restrictions apply to both station activation variables $x_{n\ell}$ and task assignment variables $y_{nk\ell}$. The budgeted variants retain objective~\eqref{app:sta-obj}; the cost-based variants use
\begin{align}
  \min_{\boldsymbol{s},\boldsymbol{o},\boldsymbol{\eta}} \sum_{n \in \mathcal{N}} p_n \cdot w_n(\boldsymbol{\eta}) \sum_{\ell \in \mathcal{L}} x_{n\ell}, \label{app:sta-obj-cpa}
\end{align}
with the commitment-dependent weights defined in Section~\ref{sec:objfunc-mald}.

\begin{table}[htbp!]
  \centering
  \small
  \renewcommand{\arraystretch}{1.15}
  \begin{tabular}{lll}
      \toprule
      \textbf{Variant} & \textbf{Objective} & \textbf{Partial adaptivity constraints} \\
      \midrule
      \bstgsta{} & \eqref{app:sta-obj} & \eqref{form:pamsp-cons3}-\eqref{form:pamsp-cons4} \\
      \cstgsta{} & \eqref{app:sta-obj-cpa} & \eqref{form:pamsp-cons1}-\eqref{form:pamsp-cons4} \\
      \bscensta{} & \eqref{app:sta-obj} & \eqref{eq:sca-nac}, \eqref{eq:sca-rev-link}, \eqref{eq:path-budget}-\eqref{eq:sca-rho-domain} \\
      \cscensta{} & \eqref{app:sta-obj-cpa} & \eqref{eq:sca-nac}, \eqref{eq:sca-rev-link}, \eqref{eq:sca-rho-domain} \\
      \bottomrule
  \end{tabular}
  \caption{Objective and partial adaptivity constraints added to the common physical constraints~\eqref{app:sta-col}-\eqref{app:sta-nonneg}.}
  \label{tab:app-sta-variants}
\end{table}

\subsection{Reconfiguration problem}

We express~\eqref{form:rec-generic} on the transformed tree $\Nrec$ described in Section~\ref{sec:lexicographic}. Under the selected revision decisions $\boldsymbol{\eta}^*$, each merged node represents a group of same-stage nodes whose strategic decisions must coincide. We retain $a(n)$ for its parent, $p_n$ for the sum of the merged nodes' probabilities, and $\mathcal{D}_n$ for the union of their demand sets. The selected counts $\zsta{n}$ and fixed weights $w_n(\boldsymbol{\eta}^*)$ are common within each group and are inherited by its merged node.

To linearize the positive-part terms in the \rec{} objective, we introduce nonnegative variables $\zeta_{nk\ell}$. The resulting formulation is
\begin{subequations}\label{app:rec}
  \small
  \begin{align}
      (\text{\rec{}}): \  \min_{\boldsymbol{s},\boldsymbol{o},\boldsymbol{\zeta}}\ & \sum_{n \in \Nrec} p_n \cdot w_n(\boldsymbol{\eta}^*) \cdot \sum_{k \in \mathcal{K}} \sum_{\ell \in \mathcal{L}} \zeta_{nk\ell} && \label{app:rec-obj} \\
      \text{s.t.} \ & \eqref{app:sta-col}\text{-}\eqref{app:sta-nonneg} && \text{for all } n \in \Nrec, \nonumber \\
      & \sum_{\ell \in \mathcal{L}} x_{n\ell} \leq \zsta{n} && \forall n \in \Nrec, \label{app:rec-z} \\
      & \zeta_{nk\ell} \geq y_{nk\ell} - y_{a(n),k\ell} && \forall k \in \mathcal{K},\ \ell \in \mathcal{L},\ n \in \Nrec, \label{app:rec-lin} \\
      & \zeta_{nk\ell} \geq 0 && \forall k \in \mathcal{K},\ \ell \in \mathcal{L},\ n \in \Nrec, \label{app:rec-nonneg} \\
            & y_{a(1),k\ell} = 0 && \forall k \in \mathcal{K},\ \ell \in \mathcal{L}. \label{app:rec-green} 
  \end{align}
\end{subequations}

Constraints~\eqref{app:rec-z} impose the station-count limits inherited from \sta{} in the hierarchical optimization described in Section~\ref{sec:objfunc-mald}. Constraints~\eqref{app:rec-lin} and \eqref{app:rec-nonneg} linearize the assignment additions. Since $p_n w_n(\boldsymbol{\eta}^*)>0$, every optimal solution satisfies $\zeta_{nk\ell}={[y_{nk\ell}-y_{a(n),k\ell}]}^+$, so objective~\eqref{app:rec-obj} equals that in~\eqref{form:rec-generic}. Constraint~\eqref{app:rec-green} imposes the greenfield initial condition used in the numerical study.

%%%%%%%%%%%%%%%%%%%%%%%%%%%%%%%%%%%%%%%%
%% Proofs
\section{Proofs}\label{app:proofs}

\optimalvaluecomparisons*
\begin{proof}
For the budgeted pair, we keep the strategic decisions of a feasible
\bstgmsp{} solution and set $\rho_n:=r_{t(n)+1}$ at each non-leaf node.
The revision-policy linking restrictions hold because these indicators are
equal within each stage. For two distinct nodes at stage $t$, by~\eqref{eq:sca-gamma},
$\Gamma_{n_1,n_2}(\boldsymbol{\rho})=0$ exactly when no revision opportunity
is scheduled after the stage of their common ancestor and through stage $t$.
This is also when \bstgmsp{} requires their strategic decisions to coincide.
For the other node pairs, the strategic linking restrictions in both models
are nonbinding under the bounds $M_i$. Every scenario path has
$\sum_{t\in[2,T]}r_t$ revision opportunities, preserving the budget.
Thus, the mapped solution is feasible for \bscenmsp{} with the same objective
value, giving $z_{\bscenmsp{}}(\mu)\le z_{\bstgmsp{}}(\mu)$.

We apply the same mapping to a feasible \cstgmsp{} solution. The linking
restrictions are satisfied by the preceding argument. Along every scenario
path, the mapped policy authorizes revision opportunities at exactly the
stages selected by $\boldsymbol{r}$. Both formulations measure lag from the
latest such stage reached, using stage~1 before the first revision opportunity.
Hence, by~\eqref{eq:stg-lag} and~\eqref{eq:sca-lag},
$\delta_n(\boldsymbol{\rho})=\delta_{t(n)}(\boldsymbol{r})$ at every node.
The cost multipliers and objective value are therefore unchanged, giving
$z_{\cscenmsp{}}(h)\le z_{\cstgmsp{}}(h)$.

Every feasible budgeted solution has MSP-feasible strategic decisions and
the MSP objective, so $z_{\msp{}}\le z_{\bscenmsp{}}(\mu)$.
Conversely, setting every $r_t=1$ in \cstgmsp{} makes the strategic linking
restrictions nonbinding and all lags zero. Any MSP-feasible strategic decisions
therefore remain feasible with the same objective because $h(0)=1$, proving
$z_{\cstgmsp{}}(h)\le z_{\msp{}}$.

When $h\equiv1$, every feasible cost-based solution has MSP-feasible strategic
decisions and the MSP objective. Full revision gives the reverse comparison,
proving both cost-based equalities. When $\mu=T-1$, full revision is also
budget-feasible in both timing formulations, proving both budgeted equalities.
\end{proof}

\primaryoptima*
\begin{proof}
Let $(\boldsymbol{s},\boldsymbol{o})$ be feasible for~\eqref{form:rec-generic}.
Then $(\boldsymbol{s},\boldsymbol{o},\boldsymbol{\eta}^*)$ is feasible
for~\eqref{form:sta-generic}, and the station-count bounds ensure that its
\sta{} objective cannot exceed that of the selected optimum. Optimality
therefore forces equal objective values. Since
$p_n w_n(\boldsymbol{\eta}^*)>0$ for every node, every station-count bound
holds with equality.

Conversely, let $(\boldsymbol{s},\boldsymbol{o},\boldsymbol{\eta}^*)$ be
\sta{}-optimal. With $\boldsymbol{\eta}^*$ fixed, the strategic-decision linking restrictions require common layouts within groups of same-stage nodes, while physical feasibility is independent across these groups. If any group were to use more stations than in the selected solution, replacing its layout and node-specific operational allocations by those of the selected solution would preserve feasibility and strictly decrease the \sta{} objective, contradicting the optimality of $(\boldsymbol{s},\boldsymbol{o},\boldsymbol{\eta}^*)$. Thus all station-count bounds hold, and the remaining \rec{} constraints already hold. Hence $(\boldsymbol{s},\boldsymbol{o})$ is \rec{}-feasible. The equivalence also establishes the final claim. 
\end{proof}

\commitmentstructure*
\begin{proof}
(a) Fix $t\in[2,T]$ and $\{n_1,n_2\}\in\Npairs_t$, and let
$a:=a(n_1,n_2)$. By~\eqref{eq:sca-gamma},
$\Gamma_{n_1,n_2}(\boldsymbol{\rho})$ is a sum of binary
indicators and hence is zero exactly when all indicators on
$\mathcal{P}(a,n_1)\cup\mathcal{P}(a,n_2)$ are zero. For $i\in\{1,2\}$,
each child reached along the path from $a$ to $n_i$ then retains its parent's
commitment start. Thus,
$\sigma_{n_1}(\boldsymbol{\rho})=\sigma_a(\boldsymbol{\rho})
=\sigma_{n_2}(\boldsymbol{\rho})$.

Conversely, suppose
$u:=\sigma_{n_1}(\boldsymbol{\rho})=\sigma_{n_2}(\boldsymbol{\rho})$.
The common start $u$ lies on both paths from the root to $n_1$ and $n_2$,
and hence on the path from the root to $a$.
For $i\in\{1,2\}$, any revision-enabling node in $\mathcal{P}(a,n_i)$
would start a new commitment at its child on that path, strictly after $a$.
This would contradict $\sigma_{n_i}(\boldsymbol{\rho})=u$, since $u$ is
the most recent commitment start on the path to $n_i$.
Thus, all indicators on $\mathcal{P}(a,n_1)\cup\mathcal{P}(a,n_2)$ are zero, giving
$\Gamma_{n_1,n_2}(\boldsymbol{\rho})=0$.

(b) At the root, $\delta_1(\boldsymbol{\rho})=0$ and
$\sigma_1(\boldsymbol{\rho})=1$. Now fix $n\in\mathcal{N}\setminus\{1\}$ and
let $u:=\sigma_n(\boldsymbol{\rho})$. If $u\neq1$, its parent $a(u)$ is
revision-enabling and all nodes in $\mathcal{P}(u,n)$ are
commitment-preserving. The maximum in~\eqref{eq:sca-lag} therefore equals
$t(a(u))=t(u)-1$. If $u=1$, all indicators on the nonempty path
$\mathcal{P}(1,n)$ are zero, so that maximum is $0=t(1)-1$.
Hence $\delta_n(\boldsymbol{\rho})=t(n)-t(u)$.

(c) Suppose $\boldsymbol{\rho}$ satisfies~\eqref{eq:sca-rev-link} and fix
$u\in\mathcal{B}(\boldsymbol{\rho})$. For $t\in[t(u),T]$, define
$\mathcal{H}_t:=\{n\in\mathcal{N}_t:\sigma_n(\boldsymbol{\rho})=u\}$.
Initially, $\mathcal{H}_{t(u)}=\{u\}=\mathcal{C}_{t(u)}(u)$.
At any stage $t$ with $t(u)\le t<T$ and
$\mathcal{H}_t=\mathcal{C}_t(u)$, part~(a) and~\eqref{eq:sca-rev-link}
give all nodes in this group the same indicator, trivially so for a singleton.
If this indicator is zero, their children retain starting node $u$, giving
$\mathcal{H}_{t+1}=\mathcal{C}_{t+1}(u)$. If it is one, their children
start new intervals and no later descendant has starting node $u$;
thus, $e_u(\boldsymbol{\rho})=t$. If no such stage occurs, then
$\mathcal{H}_T=\mathcal{C}_T(u)$ and $e_u(\boldsymbol{\rho})=T$.
This proves the stated partition and indicator pattern, including $t(u)=T$.
\end{proof}

\intervaldecomposition*
\begin{proof}
Fix a binary vector $\boldsymbol{\rho}$ satisfying~\eqref{eq:sca-rev-link}
and the applicable revision-budget constraints.
By Lemma~\ref{lem:commitment-structure}(c), the stage-$t$ commitment groups $\mathcal{C}_t(u)$,
with $u\in\mathcal{B}(\boldsymbol{\rho})$ and
$t\in[t(u),e_u(\boldsymbol{\rho})]$, partition $\mathcal{N}$.
For any such $u$ and $t$, part~(a) and the bounds $M_i$, $i\in\mathcal{I}$,
imply that~\eqref{eq:sca-nac} imposes a common strategic vector within
$\mathcal{C}_t(u)$ and adds no restriction between distinct groups.
Lemma~\ref{lem:commitment-structure}\textup{(b)} gives the common lag $t-t(u)$
within this group.
Property~\ref{property:node-local} makes strategic cost and feasibility
node-local, while operational cost and feasibility are node-local in the
underlying model. Thus, no remaining cost or feasibility constraint couples
decisions in different groups. The common vector
therefore optimizes over $\bigcap_{n\in\mathcal{C}_t(u)}\mathcal{X}_n$,
with objective and optimal value given by~\eqref{eq:unit-cost}.
If every group is feasible, any feasible solution has cost at least the
sum of the unit minima. Conversely, combining unit minimizers gives a
feasible solution attaining this sum. If any group is infeasible, both
sides equal $+\infty$.
\end{proof}

\commitmentdpexactness*
\begin{proof}
We first show that recursive end-stage choices represent exactly the binary
policies satisfying~\eqref{eq:sca-rev-link}. Starting at the root, consider
any such choices. For each commitment start $u\in\mathcal{N}$ with chosen
end stage $e\in[t(u),T]$, set $\rho_n=0$ for
$n\in\mathcal{C}_t(u)$ with $t(u)\le t<e$, and set $\rho_n=1$ for
$n\in\mathcal{C}_e(u)$ when $e<T$. The next commitments start at the nodes
in $\mathcal{C}_{e+1}(u)$. This construction assigns every non-leaf indicator
and gives $\sigma_n(\boldsymbol{\rho})=u$ throughout the interval.
The sets $\mathcal{C}_t(u)$ over the selected intervals partition $\mathcal{N}$.
Fix a stage $t$ with $2\le t<T$ and a pair
$\{n_1,n_2\}\in\Npairs_t$. If
$\Gamma_{n_1,n_2}(\boldsymbol{\rho})=0$,
part~(a) of Lemma~\ref{lem:commitment-structure} gives a common start
$u:=\sigma_{n_1}(\boldsymbol{\rho})=\sigma_{n_2}(\boldsymbol{\rho})$.
For its selected end stage $e\in[t(u),T]$, both nodes lie in
$\mathcal{C}_t(u)$, whose indicators are zero for $t<e$ and one for $t=e<T$.
If $\Gamma_{n_1,n_2}(\boldsymbol{\rho})>0$, the binary indicators already
satisfy both linking inequalities. Thus~\eqref{eq:sca-rev-link} holds.
Conversely, part~(c) of Lemma~\ref{lem:commitment-structure} shows that every policy satisfying these
restrictions has this form.

For \cscenmsp{}, we establish~\eqref{eq:dp-cscen} by backward induction on
$t(u)$. At $u\in\mathcal{N}_T$, the only end stage is $T$, giving
$V(u)=c(u,T)$. Now fix $u\in\mathcal{N}$ with $t(u)<T$ and an end stage
$e\in[t(u),T]$. The decomposition~\eqref{eq:interval-decomp} gives the interval
cost $\sum_{t=t(u)}^e c(u,t)$. When $e<T$, the successor subtrees rooted at
$u'\in\mathcal{C}_{e+1}(u)$ are disjoint, and the construction permits
independent choices in them. By the induction hypothesis, their cost-to-go
values sum to $\sum_{u'\in\mathcal{C}_{e+1}(u)}V(u')$.
When $e=T$, this sum is empty. Minimizing over $e$ gives~\eqref{eq:dp-cscen},
and the root value is $V(1)$.

For \bscenmsp{}, fix $u\in\mathcal{N}$ and a remaining budget
$b\in[0,\mu]$. If $b=0$, only $e=T$ is allowed. If $b\ge1$, any end stage
$e<T$ consumes one revision on every path through $u$, leaving $b-1$ at each
successor $u'\in\mathcal{C}_{e+1}(u)$. The same induction gives~\eqref{eq:dp-bscen}.

For the stage-wise variants, fix $t\in[1,T]$.
All starts $u\in\mathcal{N}_t$ share an end stage $e\in[t,T]$; when $e<T$,
the next starts are precisely $\mathcal{N}_{e+1}$. Summing their interval
costs and applying the same budget argument gives~\eqref{eq:dp-cstg}
and~\eqref{eq:dp-bstg}, with zero cost-to-go beyond stage $T$.

If the corresponding DP value is finite, following minimizing end stages
yields an admissible revision policy whose selected units all have finite
costs. By assumption, their minima in~\eqref{eq:unit-cost} are attained.
Assigning each unit's minimizing strategic vector to every node in its
commitment group therefore gives optimal strategic decisions by
Proposition~\ref{prop:interval-decomp}.
A DP value of $+\infty$ means that every admissible policy contains an
infeasible unit, so the formulation is infeasible.
\end{proof}

\revisionmonotonicity*
\begin{proof}
Fix admissible policies $\boldsymbol{\rho}$ and $\boldsymbol{\rho}'$ with
$\boldsymbol{\rho}\le\boldsymbol{\rho}'$.
If $v(\boldsymbol{\rho})=+\infty$, the claim is immediate.
Otherwise, for any $t\in[2,T]$ and
$\{n_1,n_2\}\in\Npairs_t$, the definition~\eqref{eq:sca-gamma} gives
$\Gamma_{n_1,n_2}(\boldsymbol{\rho}')\ge\Gamma_{n_1,n_2}(\boldsymbol{\rho})$.
Since $M_i\ge0$ for every $i\in\mathcal{I}$, replacing $\boldsymbol{\rho}$
by $\boldsymbol{\rho}'$ relaxes the strategic-decision linking
restrictions~\eqref{eq:sca-nac}. The objective and all other constraints
on strategic and operational decisions are unchanged.
Every feasible solution under $\boldsymbol{\rho}$ therefore remains
feasible under $\boldsymbol{\rho}'$ with the same cost, giving
$v(\boldsymbol{\rho}')\le v(\boldsymbol{\rho})$.
The stage-wise mapping preserves componentwise order and fixed-policy
values, so the result also holds for \bstgmsp{}.
\end{proof}

\eocvalidity*
\begin{proof}
(a) Fix an admissible policy $\boldsymbol{\rho}$ with
$v(\boldsymbol{\rho})<\infty$ and $n\in\mathcal{N}$.
For any nonempty node groups
$\mathcal{N}_g\subseteq\mathcal{N}_{g'}\subseteq\mathcal{N}$,
a layout feasible at every node in $\mathcal{N}_{g'}$ is feasible at every
node in $\mathcal{N}_g$. Hence
$\mathrm{SP}^*(\mathcal{N}_g)\le\mathrm{SP}^*(\mathcal{N}_{g'})$.
If $\rho_u=0$ for every $u\in\mathcal{P}(\sigma_n^j,n)$, then
$\sigma_n(\boldsymbol{\rho})$ is $\sigma_n^j$ or an ancestor, so
$\mathcal{C}_{t(n)}(\sigma_n^j)\subseteq
\mathcal{C}_{t(n)}(\sigma_n(\boldsymbol{\rho}))$. Thus
$v_n(\boldsymbol{\rho})\ge v_n^j$ in \bscensta{}. If, additionally,
$\rho_{a(\sigma_n^j)}=1$ when $\sigma_n^j\neq1$, the commitment start is
exactly $\sigma_n^j$. The group is then unchanged, and
Lemma~\ref{lem:commitment-structure}\textup{(b)} gives the same lag.
Hence $v_n(\boldsymbol{\rho})=v_n^j$ in \cscensta{}.
When the relevant conditions fail, the cut's right-hand side is nonpositive,
proving validity since $v_n(\boldsymbol{\rho})\ge0$.
At $\boldsymbol{\rho}=\bar{\boldsymbol{\rho}}^j$, the conditions hold for
every $n\in\mathcal{N}$, so each cut has right-hand side $v_n^j$.

(b) Fix an admissible policy $\boldsymbol{\rho}$ and $n\in\mathcal{N}$.
The nodes $u\in\mathcal{P}(\sigma_n^j,n)$ are a subset of the non-leaf nodes
with $\bar{\rho}_u^j=0$. The cost-based EOC additionally includes
$a(\sigma_n^j)$ when $\sigma_n^j\neq1$, whose evaluated indicator is one.
Since $\rho_u\in\{0,1\}$ for every $u\in\mathcal{N}\setminus\mathcal{N}_T$
and $v_n^j\ge0$, replacing these subsets by the full index sets weakens the bound.
Multiplying the weakened inequalities by $p_n$ and summing over
$n\in\mathcal{N}$, using $\theta=\sum_{n\in\mathcal{N}}p_n\theta_n$
and $\sum_{n\in\mathcal{N}}p_n v_n^j=v^j$, gives the corresponding MOC
or IOC with $L=0$. For stage-wise policies, each
$u\in\mathcal{P}(\sigma_n^j,n)$ maps to a distinct revision stage $t(u)+1$.
When $\sigma_n^j\neq1$, its parent's indicator maps to stage
$t(\sigma_n^j)$, outside those stages. The same subset argument therefore
applies to the stage-indexed cuts without repeating an indicator.
\end{proof}

\weightedswitchcut*
\begin{proof}
Fix a non-leaf node $n\in\Nrec$, a child $n'\in\mathcal{C}(n)$, and binary
$\boldsymbol{y}_n$.
For any child layout $\boldsymbol{s}_{n'}\in\mathcal{X}_{n'}$ satisfying
$\sum_{\ell\in\mathcal{L}}x_{n'\ell}\le\zsta{n'}$, binary assignments give,
for every $(k,\ell)\in\mathcal{K}\times\mathcal{L}$,
\[
    (1-y_{nk\ell})y_{n'k\ell}
    \ge (1-y_{nk\ell}^i)y_{n'k\ell}
       -(1-y_{nk\ell}^i)y_{nk\ell}.
\]
For a fixed child layout, the expected cost-to-go from its children is
independent of the incoming assignment, as is the feasible set of child
layouts. Multiplying by $w_{n'}(\boldsymbol{\eta}^*)$, summing over
$(k,\ell)\in\mathcal{K}\times\mathcal{L}$, adding the same expected
cost-to-go from the children of $n'$, and minimizing over the feasible
layouts gives
\[
    \mathcal{Q}_{n'}(\boldsymbol{y}_n)
    \ge v_{n'}^i-w_{n'}(\boldsymbol{\eta}^*)
    \sum_{k\in\mathcal{K}}\sum_{\ell\in\mathcal{L}}
    (1-y_{nk\ell}^i)y_{nk\ell}.
\]
Multiplying by $q_{nn'}$ and summing over $n'\in\mathcal{C}(n)$
gives the claimed lower bound, since the child weights coincide and
$\sum_{n'\in\mathcal{C}(n)}q_{nn'}=1$.
At $\boldsymbol{y}_n=\boldsymbol{y}_n^i$, the subtraction term is zero,
and the intercept is the exact expected child value, proving tightness.
\end{proof}

\commitmentgain*
\begin{proof}
  Fix $t\in[2,T]$. Under full adaptivity, every non-leaf
  revision indicator equals one and every node has lag zero. Setting
  $\rho_u=0$ for $u\in\mathcal{N}_{t-1}$ affects only its children: they
  share one layout with lag one, while all other nodes retain lag zero.
  By Proposition~\ref{prop:interval-decomp} and~\eqref{eq:unit-cost-sta},
  the resulting reduction in optimal cost is $\gamma_u$. The child sets
  $\mathcal{C}(u)$ are disjoint across $u\in\mathcal{N}_{t-1}$, so these
  reductions add when several such indicators are set to zero. These policies
  satisfy~\eqref{eq:sca-rev-link}: indicators can differ within
  $\mathcal{N}_{t-1}$ because all earlier indicators equal one, while every
  indicator outside $\mathcal{N}_{t-1}$ remains one. A stage-wise policy can set
  either all indicators at stage $t-1$ to zero or none, giving
  $\big[\sum_{u\in\mathcal{N}_{t-1}}\gamma_u\big]^+$. A scenario-adaptive
  policy can select any subset, giving
  $\sum_{u\in\mathcal{N}_{t-1}}[\gamma_u]^+$. If the optimal \cstgsta{}
  cost equals the \msp{} cost, a positive stage-wise single-stage gain would
  contradict optimality, proving the first conditional claim. The best
  scenario-adaptive single-stage policy at each $t$ is feasible for
  \cscensta{}, which proves the stated lower bound on its improvement.
\end{proof}

%%%%%%%%%%%%%%%%%%%%%%%%%%%%%%%%%%%%%%%%
%%%%%%%%%%%%%%%%%%%%%%%%%%%%%%%%%%%%%%%%

\section{Objective Values} \label{app:obj-levels}

  \begin{table}[htbp!]
      \centering
      \footnotesize
      \setlength{\tabcolsep}{3pt}
      \begin{tabular}{ll*{12}{c}}
          & & \multicolumn{4}{c}{\textbf{M-HM}} & \multicolumn{4}{c}{\textbf{M-HT}} & \multicolumn{4}{c}{\textbf{M-RG}} \\
          \cmidrule(lr){3-6} \cmidrule(lr){7-10} \cmidrule(lr){11-14} & & S4 & S6 & S9 & S12 & S4 & S6 & S9 & S12 & S4 & S6 & S9 & S12 \\ \hline
          \multicolumn{2}{l}{\msp{}} & 19.39 & 28.64 & 44.34 & 57.76 & 20.21 & 29.47 & 45.45 & 59.24 & 23.81 & 38.68 & 62.08 & 101.12 \\ \hline
          \multirow{2}{*}{$c = 0.03$} & \twosp{} & 20.83 & 29.48 & 45.06 & 54.47 & 30.51 & 40.83 & 59.10 & 70.56 & 27.62 & 43.19 & 67.61 & 87.88 \\
          & \cscensta{} & 19.37 & 28.28 & 43.44 & 54.30 & 20.20 & 29.29 & 44.99 & 58.27 & 23.40 & 36.98 & 56.58 & 87.88 \\ \hline
          \multirow{2}{*}{$c = 0.05$} & \twosp{} & 20.19 & 27.80 & 40.55 & 46.72 & 29.38 & 38.18 & 52.60 & 59.50 & 26.70 & 40.64 & 60.65 & 75.01 \\
          & \cscensta{} & 19.33 & 27.76 & 40.54 & 46.72 & 20.19 & 29.17 & 44.61 & 56.20 & 23.13 & 35.72 & 52.73 & 75.01 \\ \hline
          \multirow{2}{*}{$c = 0.10$} & \twosp{} & 18.58 & 23.57 & 29.31 & 30.02 & 26.56 & 31.53 & 36.44 & 36.03 & 24.40 & 34.26 & 43.33 & 47.21 \\
          & \cscensta{} & 18.53 & 23.57 & 29.31 & 30.02 & 20.13 & 28.73 & 36.10 & 36.03 & 22.45 & 32.55 & 41.88 & 47.21 \\ \midrule \midrule
          & & \multicolumn{4}{c}{\textbf{Q-HM}} & \multicolumn{4}{c}{\textbf{Q-HT}} & \multicolumn{4}{c}{\textbf{Q-RG}} \\
          \cmidrule(lr){3-6} \cmidrule(lr){7-10} \cmidrule(lr){11-14}  & & S4 & S6 & S8 & S12 & S4 & S6 & S8 & S12 & S4 & S6 & S8 & S12 \\ \hline
          \multicolumn{2}{l}{\msp{}} & 20.14 & 29.47 & 38.87 & 57.02 & 21.00 & 30.46 & 40.17 & 57.91 & 24.10 & 38.77 & 52.56 & 88.34 \\ \hline
          \multirow{2}{*}{$c = 0.03$} & \twosp{} & 20.95 & 30.79 & 38.46 & 53.08 & 32.02 & 42.03 & 52.13 & 66.83 & 30.32 & 45.89 & 61.19 & 76.42 \\
          & \cscensta{} & 19.90 & 28.85 & 37.85 & 53.06 & 20.97 & 30.24 & 39.76 & 56.58 & 23.67 & 36.94 & 49.01 & 76.36 \\ \hline
          \multirow{2}{*}{$c = 0.05$} & \twosp{} & 20.25 & 29.05 & 35.28 & 45.59 & 30.83 & 39.38 & 47.40 & 56.21 & 29.20 & 43.14 & 55.97 & 65.67 \\
          & \cscensta{} & 19.72 & 28.19 & 35.27 & 45.59 & 20.95 & 30.07 & 39.40 & 54.03 & 23.37 & 35.68 & 46.35 & 65.67 \\ \hline
          \multirow{2}{*}{$c = 0.10$} & \twosp{} & 18.50 & 24.67 & 27.27 & 29.52 & 27.86 & 32.73 & 35.46 & 33.84 & 26.40 & 36.26 & 42.80 & 42.39 \\
          & \cscensta{} & 18.50 & 24.67 & 27.27 & 29.52 & 20.90 & 29.55 & 34.78 & 33.84 & 22.60 & 32.48 & 39.68 & 42.39
      \end{tabular}
      \caption{Optimal values of the station problem (expected discounted station-use cost over the horizon) for \msp{}, \twosp{}, and \cscensta{} by dataset and number of stages, averaged over seeds and the three discount shapes. The value of \msp{} does not depend on $c$ or on the shape of $h$. The RG datasets are described in Section~\ref{sec:value-pa}.}
      \label{tab:obj-levels}
  \end{table}
  
\end{document}

%% file: Figures/tree_overview.tex
% Composite scenario tree (left) and zoom on n2, n4, n5 (right)
\begin{tikzpicture}[line width=0.6pt, >={Straight Barb[length=1.3mm, width=2.2mm]}]
  % ---------- left: scenario tree ----------
  \begin{scope}[sctree]
    \node[stnode] (n1) at (0,3.5) {$n_1$};
    \foreach \i/\y in {2/5.5, 3/1.5}
      \node[stnode] (n\i) at (1,\y) {$n_{\i}$};
    \foreach \i/\y in {4/6.5, 5/4.5, 6/2.5, 7/0.5}
      \node[stnode] (n\i) at (2,\y) {$n_{\i}$};
    \foreach \i/\y in {8/7, 9/6, 10/5, 11/4, 12/3, 13/2, 14/1, 15/0}
      \node[stnode] (n\i) at (3,\y) {$n_{\i}$};
    \foreach \p/\c in {1/2, 1/3, 2/4, 2/5, 3/6, 3/7, 4/8, 4/9,
                       5/10, 5/11, 6/12, 6/13, 7/14, 7/15}
      \draw[stedge] (n\p.east) -- (n\c.west);
  \end{scope}
  % zoom box around n2, n4, n5 with magnifier icon
  \coordinate (zbLL) at ($(n2.west |- n5.south)+(-0.15cm,-0.15cm)$);
  \coordinate (zbUR) at ($(n5.east |- n4.north)+(0.15cm,0.45cm)$);
  \draw[densely dotted, line width=0.5pt] (zbLL) rectangle (zbUR);
  \coordinate (mg) at ($(zbLL |- zbUR)+(0.22cm,-0.2cm)$);
  \draw[line width=0.5pt] (mg) circle (0.1cm);
  \draw[line width=1.1pt] ($(mg)+(0.07cm,-0.07cm)$) -- ++(0.12cm,-0.12cm);
  \draw[line width=0.4pt] ($(mg)+(-0.05cm,0)$) -- ++(0.1cm,0)
                          ($(mg)+(0,-0.05cm)$) -- ++(0,0.1cm);

  % ---------- right: composite nodes ----------
  % \compnode{name}{(x,y)}{label}: strategic substage (S) above,
  % operational substage (O) with demand realizations D_n below
  \def\compnode#1#2#3{%
    \begin{scope}[shift={#2}]
      \draw[rounded corners=8pt] (-1.2,-1.1) rectangle (1.2,1.1);
      \node at (0,1.42) {$#3$};
      \draw[dashed, line width=0.4pt] (-1.2,0) -- (1.2,0);
      \node[anchor=north west, font=\footnotesize, inner sep=3pt] at (-1.2,1.1) {S};
      \node[anchor=north west, font=\footnotesize, inner sep=3pt] at (-1.2,0) {O};
      \node[anchor=north east, font=\scriptsize, inner sep=3pt] at (1.2,0) {$\mathcal{D}_{#3}$};
      \node[circle, draw, minimum size=0.44cm, inner sep=0pt, fill=white] (#1) at (0,0.55) {};
      \node[circle, draw, densely dotted, minimum size=0.4cm, inner sep=0pt] (#1l) at (-0.6,-0.78) {};
      \node[circle, draw, densely dotted, minimum size=0.4cm, inner sep=0pt] (#1r) at (0.6,-0.78) {};
      \node[font=\scriptsize] at (0,-0.8) {$\cdots$};
      \draw[line width=0.4pt] (#1.south) -- (#1l.north) (#1.south) -- (#1r.north);
      \foreach \x in {-0.15,0,0.15}
        \draw[line width=0.4pt] (#1.south) -- (\x,-0.5);
    \end{scope}}
  \begin{scope}[shift={(9.2cm,2.625cm)}]
    \draw[densely dotted, line width=0.5pt] (-2.1,-3.6) rectangle (5.8,3.6);
    \begin{scope}[yshift=-0.25cm]
      \compnode{z2}{(0,0)}{n_2}
      \compnode{z4}{(3.8,1.9)}{n_4}
      \compnode{z5}{(3.8,-1.9)}{n_5}
      
        \draw[->] (-2.1,-0.4) -- (z2);
        \draw[->] (z2) -- (z4);
        \draw[->] (z2) -- (z5);
        \foreach \z in {z4,z5} {
          \draw[->] (\z) -- ($(\z)+(1.8,0.45)$);
          \draw[->] (\z) -- ($(\z)+(1.8,-0.3)$);
        }
    \end{scope}
  \end{scope}
\end{tikzpicture}

%% file: Figures/tree_mu0.tex
% Full commitment (2SP, mu = 0): all stage-t nodes merged
\begin{tikzpicture}[sctree]
  \sctreebbox
  \node[stnode]   (a) at (0,3.5) {$n_1$};
  \node[stmerged] (b) at (1,3.5) {$n_{2\text{--}3}$};
  \node[stmerged] (c) at (2,3.5) {$n_{4\text{--}7}$};
  \node[stmerged] (d) at (3,3.5) {$n_{8\text{--}15}$};
  \draw[stedge] (a) -- (b);
  \draw[stedge] (b) -- (c);
  \draw[stedge] (c) -- (d);
\end{tikzpicture}

%% file: Figures/tree_stage_mu2.tex
% Stage-wise revision, mu = 2: r_2 = 0, r_3 = r_4 = 1
\begin{tikzpicture}[sctree]
  \sctreebbox
  \node[stnode]   (n1)  at (0,3.5) {$n_1$};
  \node[stmerged] (n23) at (1,3.5) {$n_{2\text{--}3}$};
  \foreach \i/\y in {4/6.5, 5/4.5, 6/2.5, 7/0.5}
    \node[stnode] (n\i) at (2,\y) {$n_{\i}$};
  \foreach \i/\y in {8/7, 9/6, 10/5, 11/4, 12/3, 13/2, 14/1, 15/0}
    \node[stnode] (n\i) at (3,\y) {$n_{\i}$};
  \draw[stedge] (n1.east) -- (n23.west);
  \foreach \c in {4,5,6,7}
    \draw[stedge] (n23.east) -- (n\c.west);
  \foreach \p/\c in {4/8, 4/9, 5/10, 5/11, 6/12, 6/13, 7/14, 7/15}
    \draw[stedge] (n\p.east) -- (n\c.west);
\end{tikzpicture}

%% file: Figures/tree_msp.tex
% Full adaptivity (MSP): binary tree, T = 4
\begin{tikzpicture}[sctree]
  \sctreebbox
  \node[stnode] (n1) at (0,3.5) {$n_1$};
  \foreach \i/\y in {2/5.5, 3/1.5}
    \node[stnode] (n\i) at (1,\y) {$n_{\i}$};
  \foreach \i/\y in {4/6.5, 5/4.5, 6/2.5, 7/0.5}
    \node[stnode] (n\i) at (2,\y) {$n_{\i}$};
  \foreach \i/\y in {8/7, 9/6, 10/5, 11/4, 12/3, 13/2, 14/1, 15/0}
    \node[stnode] (n\i) at (3,\y) {$n_{\i}$};
  \foreach \p/\c in {1/2, 1/3, 2/4, 2/5, 3/6, 3/7, 4/8, 4/9,
                     5/10, 5/11, 6/12, 6/13, 7/14, 7/15}
    \draw[stedge] (n\p.east) -- (n\c.west);
\end{tikzpicture}

%% file: Figures/tree_scen_mu2.tex
% Scenario-adaptive revision, mu = 2: rho_{n2} = rho_{n6} = rho_{n7} = 0, others 1
\begin{tikzpicture}[sctree]
  \sctreebbox
  \node[stnode]   (n1)  at (0,3.5) {$n_1$};
  \node[stnode]   (n2)  at (1,5.5) {$n_2$};
  \node[stnode]   (n3)  at (1,1.5) {$n_3$};
  \node[stmerged] (n45) at (2,5.5) {$n_{4\text{--}5}$};
  \node[stnode]   (n6)  at (2,2.5) {$n_6$};
  \node[stnode]   (n7)  at (2,0.5) {$n_7$};
  \foreach \i/\y in {8/7, 9/6, 10/5, 11/4}
    \node[stnode] (n\i) at (3,\y) {$n_{\i}$};
  \node[stmerged] (n1213) at (3,2.5) {$n_{12\text{--}13}$};
  \node[stmerged] (n1415) at (3,0.5) {$n_{14\text{--}15}$};
  \foreach \p/\c in {1/2, 1/3, 2/45, 3/6, 3/7, 45/8, 45/9, 45/10, 45/11,
                     6/1213, 7/1415}
    \draw[stedge] (n\p.east) -- (n\c.west);
\end{tikzpicture}

%% file: Figures/SolutionApproachOverview.tex
% Editable figure fragment, included by content.tex.
% Requires tikz and arrows.meta, positioning, calc; no standalone wrapper.
% All method references follow the manuscript's subsection numbering.
\begingroup
\def\solutionArrowGap{6.5mm}
\begin{tikzpicture}[
    x=1cm,y=1cm,
    connection/.style={draw=black!70,line width=0.8pt,line join=round},
    flow/.style={connection,-{Stealth[length=3mm,width=2mm]}},
    box/.style={
        draw=black!65,fill=white,line width=0.5pt,
        rounded corners=2pt,align=center,
        inner xsep=2.2mm,inner ysep=1.2mm,outer sep=0pt,
        minimum height=1.15cm,font=\small
    },
    mainbox/.style={box,text width=5cm,fill=black!3},
    altbox/.style={box,text width=5cm},
    widebox/.style={box,text width=9cm},
    tallbox/.style={minimum height=1.45cm},
    route/.style={align=center,font=\small,inner sep=0pt},
    edgelabel/.style={font=\footnotesize,inner sep=2pt}
]
    % STA: two alternative solvers for the same joint optimization problem.
    \node[align=center,font=\small\bfseries] at (1.075,2.15)
        {\sta{}: joint revision and layout optimization};

    % Center each section reference on a separate line below its method title.
    \node[route] at (-3.5,1.25)
        {\textbf{Dynamic program (DP)}\\
         {\footnotesize(Sections~\ref{sec:structure}-\ref{sec:dp})}};
    \node[route] at (3.5,1.25)
        {\textbf{Integer L-shaped alternative}\\
         {\footnotesize(Section~\ref{sec:lshaped})}};

    \node[mainbox] (stationcosts) at (-3.5,0)
        {\textbf{Station subproblems}\\
         Precompute commitment costs};
    \node[mainbox,below=\solutionArrowGap of stationcosts] (commitmentdp)
        {\textbf{Commitment DP}\\
         Select commitment end stages};

    \node[altbox] (master) at (3.5,0)
        {\textbf{Revision master problem}\\
         Select $\boldsymbol{\eta}$};
    \node[altbox,below=\solutionArrowGap of master] (lowerlayer)
        {\textbf{Fixed-revision subproblem}\\
         Optimize layout decisions};

    \draw[flow] (stationcosts.south) --
        node[edgelabel,right,xshift=1.2mm] {Reusable costs} (commitmentdp.north);
    \draw[flow] (master.south) --
        node[edgelabel,right,xshift=1.2mm] {Candidate $\boldsymbol{\eta}$} (lowerlayer.north);
    \draw[flow] (lowerlayer.east) --
        ($(lowerlayer.east)+(0.4cm,0)$) --
        ($(master.east)+(0.4cm,0)$)
        -- (master.east);
    \node[edgelabel,align=center,text width=2cm] at
        ($(master.east)!0.54!(lowerlayer.east)+(1.4cm,0)$)
        {Strengthened \\ optimality \\ cuts};

    \draw[draw=black!50,line width=0.6pt,rounded corners=5pt]
        (-6.5,2.6) rectangle ($(lowerlayer.south)+(5.2cm,-2.5mm)$);

    % Only revision decisions, counts, and weights are passed to REC.
    % STA's activation and assignment vectors are reoptimized in REC.
    \coordinate (handoff) at
        ($(commitmentdp.south)!0.5!(lowerlayer.south)+(0,-4.5mm)$);
    \node[widebox,tallbox,below=\solutionArrowGap of handoff] (staoutputs)
        {\textbf{Selected \sta{} outputs}\\
         Revision decisions $\boldsymbol{\eta}^{*}$ and station counts $\zsta{n}$\\
         Fixed weights $w_n(\boldsymbol{\eta}^{*})$};
    \draw[connection] (commitmentdp.south) |- (handoff);
    \draw[connection] (lowerlayer.south) |- (handoff);
    \draw[flow] (handoff) -- (staoutputs.north);
    
    % The straight vertical flow arrows all use the same anchor-to-anchor gap.
    \node[widebox,below=\solutionArrowGap of staoutputs] (mergedtree)
        {\textbf{Transformed scenario tree}\\
         Merge same-stage nodes whose layouts must coincide};
    \draw[flow] (staoutputs.south) -- (mergedtree.north);

    \node[widebox,tallbox,fill=black!3,below=\solutionArrowGap of mergedtree] (recsolve)
        {\textbf{\rec{}: minimize reconfiguration cost}
         (Section~\ref{sec:lexicographic})\\
         Reoptimize layouts by nested decomposition\\
         Weighted asymmetric Lagrangian switch cuts};
    \draw[flow] (mergedtree.south) -- (recsolve.north);

    \node[widebox,tallbox,below=\solutionArrowGap of recsolve] (finallayouts)
        {\textbf{Final layout plans}\\
         \sta{}-optimal layouts under $\boldsymbol{\eta}^{*}$\\
         Minimum reconfiguration cost};
    \draw[flow] (recsolve.south) -- (finallayouts.north);
\end{tikzpicture}
\endgroup